# Improving Causal Inference with Measurement Errors in Exposures and Confounders: A New Method and Its Application to Air Pollution Exposure Assessment and Epidemiology

Honghyok Kim

Division of Environmental and Occupational Health Sciences, School of Public Health, University of Illinois Chicago, Chicago, IL, 60661, USA

Corresponding author: Honghyok Kim, honghyok@uic.edu

## Abstract

When exposure measurement error (EME), confounder measurement error (CME), or both are present, health effect estimates regarding exposure mixtures and critical exposure time-window may not represent the true effects. For example, in air pollution epidemiology, modeled estimates for multiple air pollutants and meteorological factors may serve as surrogates for exposures and confounders. Methods for simultaneously addressing EME and CME remain understudied. We developed a two-stage causal effect modeling framework to estimate average exposure/treatment effects (AEE) by addressing EME and CME. We identified conditions under which AEE is identifiable with minimal bias given linear or non-linear potential outcomes models and developed a new method, referred to as multi-dimensional regression calibration (MRC). The first stage of the framework estimates MRC models. The second stage estimates AEE by using g-computation with MR-Calibrated variables. Simulation analyses confirmed the bias-correction capability. As an application, we analyzed the association between air pollution and COVID-19 mortality in Cook County, Illinois. We developed machine learning-based 500m-gridded daily estimates of air pollutants and meteorological factors in a way for what we refer to as doubly EME&CME-correction. Using distributed lag variables, a one interquartile range (22.7ppb) increase in 3-week $O_3$ exposure below 70ppb was associated with an 135.3% (95% CI: 68.4, 233.0) increase in COVID-19 mortality risk, comparable to that for $PM_{2.5}$ exposure, which contradicts the previously reported no association for 3-week $O_3$ in Cook County. At low levels, reducing pollution may have helped prevent premature deaths from COVID-19. Our new framework can address measurement error in multiple covariates simultaneously.

The supplementary materials begin after p39 (separate page numbers)

## 1. Introduction

For environmental epidemiological studies, investigators may need to address multiple continuous exposures and confounders, such as exposure mixtures[1,2], including multiple contaminants, and different (cumulative) exposure time-windows as distributed lag variables[3-8]. In identifying the health effect of one pollutant, another may be a confounder. An exposure variable in one period may be a confounder in identifying the effect of exposure in another period[8].

Measurement error is prevalent in observational studies. Exposure measurement error (EME) and confounder measurement error (CME) may introduce bias in effect estimates and compromise causal inference and reproducibility[9-26]. When direct exposure/confounder assessment is not a practical option, modelled estimates may be linked to health data, to serve as surrogates for exposures and confounders. For example, modelled estimates for environmental exposures (e.g., air pollution, temperature, chemicals including perfluoroalkyl substances, heavy metals, and other compounds) may be used[16,27-35]. Modelled estimates are inherently subject to error, reflecting the well-known notion that “all models are wrong, but some are useful”. Claims about the direction or magnitude of bias should be made with caution when error corrections or quantitative bias analyses have not been performed[10,11,14].

EME has been considered, including modelled estimates for air pollutants and statistical simulations[18,27,29,36-41]. Correcting for EME is preferable to not correcting at all[10,11,14,42,43]. Under the potential outcomes framework[44], investigators examined EME in a single continuous or binary exposure variable[45-49]. In statistics, it is established that the extent to correction for bias due to EME may depend on modeling conditions and correction methods, especially when the true outcome model or effect of continuous covariates is non-linear and EME is not small[10,43]. Non-linear models are common (e.g., logistic, log-linear models). Effects of continuous variables may be sub/non-linear[50-52]. In fact, statistical simulations might not be comprehensive, meaning the limited error-correction ability of certain methods may not be fully revealed. Empirical studies found that effect estimates may differ across correction methods[53], while others demonstrated that they may not necessarily differ[54,55].

Addressing both EME and CME simultaneously remains understudied. The importance of addressing EME and CME has likewise been emphasized in nutritional epidemiology, where dietary intake measurements are frequently error-prone[15,24,56,57]. Bayesian methods and likelihood-based methods may be used to simultaneously adjust for EME and CME[43,58]. Addressing both EME and CME can be technically intensive, may involve subjective assumptions[58], and remains not well-established within the potential outcomes framework and when correcting for error in high-dimensional covariates.

Therefore, this study addresses two practical questions:

Q1) Under what conditions the bias due to EME can be minimized in non-linear settings, we refer to as *approximation conditions*.

Q2) How to estimate causal effects by addressing EME and CME simultaneously.

We developed a new causal effect modeling framework that extends regression calibration (RC), a technically efficient method[42,43], to accommodate multiple covariates within the potential outcomes framework, and conducted simulation analyses to evaluate its performance. As its application, we analyzed the association of short-term exposure to $PM_{2.5}$ and $O_3$ with COVID-19 mortality.

## 2. The Need for Error Correction

Due to space constraints, Appendix S1 provides an illustration of EME and CME using directed acyclic graphs (Figure S1) and error in modelled estimates. Here, we illustrate how EME and CME may affect causal effect estimation under the potential outcomes framework.

Let $X$ and $X^{ep}$ denote a true continuous exposure and an exposure measured with error, respectively. *ep* means "error-prone". $Y(X)$ denotes the potential outcome when exposed to $X$. $Y(X)$ is discrete (including binary outcomes) and may be deterministic or non-deterministic[59]. $Y(x)$ denotes the potential outcome at $X = x$. $\boldsymbol{Z}'$ denotes a set of confounders. $\boldsymbol{z}'$ is a set of values. Conditional average exposure effect (AEE) in risk difference (RD) given $\boldsymbol{z}'$ is defined as

$$AEE(X|\boldsymbol{z}') \coloneqq E[Y(X = e_T + \delta) - Y(X = e_T)|X = x, \boldsymbol{z}']$$

$AEE(X|\boldsymbol{z}')$ means that given $x$ and $\boldsymbol{z}'$, the difference between the potential outcome when exposed to $e_T + \delta$ and the potential outcome when exposed to $e_T$. Conditional AEE in risk ratio (RR) is

$$AEE(X|\boldsymbol{z}') \coloneqq \frac{E[Y(X = e_T + \delta)|X = x, \boldsymbol{z}']}{E[Y(X = e_T)|X = x, \boldsymbol{z}']}$$

Under exchangeability conditional on $\boldsymbol{z}'$, which is $Y(x) \perp\!\!\!\perp X|\boldsymbol{z}'$, and the consistency assumption[60],

$AEE(X|\boldsymbol{z}') = E[Y(X = e_T + \delta)|X = x, \boldsymbol{z}'] - E[Y(X = e_T)|X = x, \boldsymbol{z}']$

$= E[Y(X = e_T + \delta)|X = e_T + \delta, \boldsymbol{z}'] - E[Y(X = e_T)|X = e_T, \boldsymbol{z}']$

$= E\left[Y_{X=e_T+\delta}|\boldsymbol{z}'\right] - E\left[Y_{X=e_T}|\boldsymbol{z}'\right]$

$Y_X$ denotes the observed outcome of those exposed to $X$. Similarly,

$$AEE(X|\boldsymbol{z}') = \frac{E[Y(X = e_T + \delta)|X = e_T + \delta, \boldsymbol{z}']}{E[Y(X = e_T)|X = e_T, \boldsymbol{z}']} = \frac{E\left[Y_{X=e_T+\delta}|\boldsymbol{z}'\right]}{E\left[Y_{X=e_T}|\boldsymbol{z}'\right]}$$

Using $X^{ep}$, $AEE(X|\boldsymbol{z}')$ is not necessarily identifiable. Nevertheless, when the effect of $X$ on $Y$ is strictly linear and the error is purely Berksonian, $AEE(X|\boldsymbol{z}')$ can be identifiable with minimal or no bias. The error is said to be purely Berksonian when this follows the form

$$X = X^{ep} + U$$

where $U$ is completely random with $E[U] = 0$. This implies $Y(x) \perp\!\!\!\perp X^{ep}|x, \boldsymbol{z}'$, which is the non-differential error assumption[45,48]. The effect of $X$ is said to be strictly linear when the true potential outcomes model, which may be referred to as the true (data-generating) outcome model, follows the form: $Y(X) = \alpha + \beta X + j(\boldsymbol{Z}') + W$. $\alpha$ is the intercept. $W$ is an independent (unknown) random variable that can subsume the effect of a set of non-confounder but risk factors, describing non-deterministic potential outcomes[59]. When the effect is strictly linear and EME is purely Berksonian,

$$E\left[E_U[Y(x)|X^{ep} = x, \boldsymbol{z}']\right] = E[Y(x)|\boldsymbol{z}']$$

implying that $AEE(X|\boldsymbol{z}')$ is identifiable.

Even if $X^{ep}$ is not purely Berksonian, RC with validation data may be used to estimate causal effects[43,45,61]. For continuous or categorical $X$[10,45,58], an RC model is fit such that

$$X = \gamma_{RC} + s(X^{ep}) + k(\boldsymbol{Z}') + \epsilon$$

$X_{RC}^{ep} := \hat{\gamma}_{RC} + \hat{s}(X^{ep}) + \hat{k}(\boldsymbol{Z}')$, referred to as predicted/calibrated variable, may be used in place of $X^{ep}$ for estimating causal effects. $\epsilon$ is the residual, independent of $X_{RC}^{ep}$, and $E[\epsilon] = 0$. An intuition is that $X_{RC}^{ep}$ is purely Berksonian[43,61]. Confounders should be included in a RC model like $k(\boldsymbol{Z}')$[10,42,43,58,61-63], which we refer to as the confounder inclusion condition.

However, even if $X^{ep}$, including $X_{RC}^{ep}$, is purely Berksonian, the expected value of effect estimates may not necessarily equal or approximate the true exposure effect. Each of the following three may introduce bias, or collectively, inflate overall bias:

1) The potential outcomes/data-generating model is unknown. Pure Berkson error may introduce non-negligible bias in estimating the model or the effect if the error is not small and the effect is not strictly linear[10,43,64,65];
2) When CME exists, exchangeability may not hold (See Appendix S1 for cases where CME may not necessarily introduce residual confounding);
3) Even if it is assumed that some part of the potential outcomes model is known, when the effect is not strictly linear, bias in estimating AEE may not necessarily be negligible (Examples 1–3 in Appendix S2).

## 3. A New Causal Effect Modeling Framework

To address Q1 and Q2, we propose a two-stage causal effect modeling framework for estimating average exposure/treatment effect (AEE) by simultaneously addressing EME and CME:

Stage 1) Develop calibration models that we refer to as multi-dimensional regression calibration (MRC) models and obtain calibrated exposure and confounder variables, referred to as MR-Calibrated variables.

Stage 2) Estimate AEE using g-computation[66,67] with MR-Calibrated variables.

### 3.1. Effect Estimands, Interpretability, and Conditions for Minimal Bias: A Solution to Q1

We define an effect estimand based on $X^{ep}$ and find conditions under which this effect estimand may equal or approximate $AEE(X|\mathbf{z}')$.

Let $Y(X=?\,|X^{ep})$ denote the potential outcome of those with $X^{ep}$ and $Y(X^{ep})$ denote $Y(X=?\,|X^{ep})$ for notational simplicity. Let $Y_{X^{ep}}$ denote the observed outcome of those with $X^{ep}$, which can be seen as observed outcomes labeled using or linked to $X^{ep}$ in an observational dataset. With the consistency assumption (Assumption 1), the value of each individual's $Y(X=e_T)$, $Y(X^{ep}=x^{ep})$, $Y_{X=e_T}$, and $Y_{X^{ep}=x^{ep}}$ is identical because they all indicate their potential outcomes, where $Y(X^{ep}=x^{ep})$ denotes $Y(X=?\,|X^{ep}=x^{ep})$. $x^{ep}$ may not necessarily be $e_T$. However, investigators cannot identify $Y(X=e_T)$ and $Y_{X=e_T}$ when $X$ is unknown. $Y_{X^{ep}=e_T}$ is identifiable but may not equal $Y_{X=e_T}$ and $Y(X=e_T)$ (See Table S1). Therefore, investigators who use $X^{ep}$ (including $X_{RC}^{ep}$) cannot directly target $AEE(X|\boldsymbol{z}')$ but may target

$$AEE(X^{ep}|\boldsymbol{z}) \coloneqq E[Y(X^{ep} = e_T + \delta) - Y(X^{ep} = e_T)|X^{ep} = x^{ep}, \boldsymbol{z}]$$

by hoping that $AEE(X^{ep}|\boldsymbol{z})$ may equal or approximate $AEE(X|\boldsymbol{z})$ or $AEE(X|\boldsymbol{z}')$ assuming that $V$ is not an effect modifier, where $\boldsymbol{z} \coloneqq \{\boldsymbol{z}', v\}$ and $v$ is a value of $V$ that denotes an error component that makes EME differential (See Appendix S1 for details). For non-differential EME, $AEE(X^{ep}|\boldsymbol{z}')$ may be targeted. The expectation for $Y(X^{ep})$ is taken regarding the (joint) distribution of $X$ given $x^{ep}$ and $\boldsymbol{z}$ and (of $Y(X)$ if $Y(X)$ is non-deterministic). This equality or approximation may sometimes fail because, for example, $E[Y(X^{ep})|X^{ep} = e_T, \boldsymbol{z}]$ may not equal $E[Y(X)|X = e_T, \boldsymbol{z}]$. See Appendix S3 with Table S1 for illustration.

Nevertheless, $AEE(X^{ep}|\boldsymbol{z}\ or\ \boldsymbol{z}')$ may be endowed with a causal interpretation, by assuming that this differs from $AEE(X|\boldsymbol{z}')$ within a reasonable/acceptable range (e.g., minimal bias). For Berkson error (Figures S1E–H), intervening $X^{ep}$ directs the change in $X$, consequently $Y(X)$, such as intervening area-level $O_3$ may change personal exposure to $O_3$ and thereby impact $Y$. For non-Berkson error (Figures S1A–D), the assumption that the change in $X^{ep}$ is subject to only or at least reflects the change in $X$ is needed. This assumption appears implicit in public health practice. Investigators may measure exposure levels (e.g., pollutant, nutrition intake, or biomarkers) using error-prone instruments and determine whether an intervention is desirable based on the measured values, which is not uncommon.

To estimate $AEE(X^{ep}|\boldsymbol{z} \text{ or } \boldsymbol{z}')$ using observational data, exchangeability under EME (EUEME) is needed.

**Definition (EUEME).** For discrete outcomes including binary outcomes, EUEME is defined as,

$$Y(X^{ep} = x^{ep}) \perp\!\!\!\perp X^{ep}|\boldsymbol{z}$$

implying

$$E[Y(X^{ep} = x^{ep})|X^{ep} = x^{ep} + \delta, \boldsymbol{z}] = E[Y(X^{ep} = x^{ep})|X^{ep} = x^{ep}, \boldsymbol{z}]$$

Note that EUEME requires being conditional on $v$ for differential EME. Otherwise, $Y(X^{ep} = x^{ep})$ would not be independent of $X^{ep}$. For non-differential EME, EUEME may be

$$Y(X^{ep} = x^{ep}) \perp\!\!\!\perp X^{ep}|\boldsymbol{z}'$$

Consequently, in the case of $AEE(X^{ep}|\boldsymbol{z})$, for example,

$AEE(X^{ep}|\boldsymbol{z}) = E[Y(X^{ep} = e_T + \delta)|X^{ep} = x^{ep}, \boldsymbol{z}] - E[Y(X^{ep} = e_T)|X^{ep} = x^{ep}, \boldsymbol{z}]$

$= E[Y(X^{ep} = e_T + \delta)|X^{ep} = e_T + \delta, \boldsymbol{z}] - E[Y(X^{ep} = e_T)|X^{ep} = e_T, \boldsymbol{z}]$ with **Assumption 2 (EUEME)**

$= E\left[Y_{X^{ep}=e_T+\delta}|\boldsymbol{z}\right] - E\left[Y_{X^{ep}=e_T}|\boldsymbol{z}\right]$ with **Assumption 1 (Consistency)**

To use $Y_{X^{ep}}$, the following two assumptions are also needed:

**Assumption 3 (Positivity).** There is a non-zero probability that every individual is exposed to $X$[60].

**Assumption 4 (Consistency in Exposure Measurement)**. $X$ of every individual is measured using the same instrument and the error model is identical for all individuals.

Finally, $AEE(X|\boldsymbol{z} \text{ or } \boldsymbol{z}')$ may be identifiable under certain circumstances because

**Theorem 1 (Minimal bias when EME is Berksonian).** When the effect of $X$ is approximately strictly linear, for purely Berkson EME,

$$E\left[\widehat{AEE}(X^{ep}|\boldsymbol{z}')\right] \approx AEE(X|\boldsymbol{z}')$$

and for $X = X^{ep} + V + U$ or $X = X^{ep} + U$ when $X^{ep}$ is correlated with $V$ that is correlated with $Y$,

$$E\left[\widehat{AEE}(X^{ep}|\boldsymbol{z})\right] \approx AEE(X|\boldsymbol{z})$$

**Theorem 2 (Minimal bias when EME is Berksonian under non-linearity).** When the effect of $X$ is not strictly linear, there are approximation conditions such that for purely Berkson EME,

$$E\left[\widehat{AEE}(X^{ep}|\boldsymbol{z}')\right] \approx AEE(X|\boldsymbol{z}')$$

and for $X = X^{ep} + V + U$ or $X = X^{ep} + U$ when $X^{ep}$ is correlated with $V$ that is correlated with $Y$,

$$E\left[\widehat{AEE}(X^{ep}|\boldsymbol{z})\right] \approx AEE(X|\boldsymbol{z})$$

When $V$ is not an effect modifier and not a confounder, $AEE(X|\boldsymbol{z}) = AEE(X|\boldsymbol{z}')$. Appendix S4 provides the proofs and details approximation conditions under non-linearity. Generally, approximation conditions may include

1) EME should be small or modest,

2) $q\big(E(Y)\big)$ should lie within a domain that this can be approximated as a linear, where $q$ is the link function, and/or

3) $h(X)$ as the exposure-response relationship function is approximately linear in the range of $X$, $(e_T + v - u_P, e_T + v + u_P)$ where $\pm u_P$ indicates the primary mass of the probability density function of $U$.

Due to non-linearity, the specific conditions vary depending on the form of the outcome model, $h(X)$, and possibly estimation methods.

### 3.2. MRC: A Solution to Q2

Let $C$ and $C^{ep}$ denote a true continuous confounder and a confounder measured with error, respectively. $V$ may also make $C^{ep}$ differential. Let $V^{ep}$ denote $V$ measured with error. When $C^{ep}$ and/or $V^{ep}$ are used in place of $C$ and/or $V$, error-prone $\boldsymbol{z}$, $\boldsymbol{z^{ep}}$, arises. EUEME may not hold under $\boldsymbol{z^{ep}}$. To address this, consider,

**Condition 1.** The main study is not subject to unmeasured confounding, selection bias, and outcome misclassification, although EME and CME may be present.

**Condition 2.** Using the validation study, the following series of models, referred to as MRC models, can be estimated and transportable to the main study. For *P* variables for exposures of interest $\{X_1 \ldots, X_P\}$, *Q* confounder variables $\{C_1 \ldots, C_Q\}$ that may include exposures not of interest, while $\boldsymbol{Z}'^{-\boldsymbol{C}}$ denotes a set of all confounders excluding $\boldsymbol{C} = \{C_1 \ldots, C_Q\}$,

$$X_1 = \gamma^*_{0,1} + \sum_{p=1}^{P} s_{X_1,X_p^{ep}}\left(X_p^{ep}\right) + \sum_{q=1}^{Q} s_{X_1,C_q^{ep}}\left(C_q^{ep}\right) + s_{X_1,V^{ep}}\left(V^{ep}\right) + s_{X_1,\boldsymbol{Z}'^{-C}}\left(\boldsymbol{Z}'^{-\boldsymbol{C}}\right) + U_1^*$$

$$\ldots$$

$$X_P = \gamma^*_{0,P} + \sum_{p=1}^{P} s_{X_P,X_p^{ep}}\left(X_p^{ep}\right) + \sum_{q=1}^{Q} s_{X_P,C_q^{ep}}\left(C_q^{ep}\right) + s_{X_P,V^{ep}}\left(V^{ep}\right) + s_{X_P,\boldsymbol{Z}'^{-C}}\left(\boldsymbol{Z}'^{-\boldsymbol{C}}\right) + U_P^*$$

$$C_1 = \alpha^*_{0,1} + \sum_{p=1}^{P} s_{C_1,X_p^{ep}}\left(X_p^{ep}\right) + \sum_{q=1}^{Q} s_{C_1,C_q^{ep}}\left(C_q^{ep}\right) + s_{C_1,V^{ep}}\left(V^{ep}\right) + s_{C_1,\boldsymbol{Z}'^{-C}}\left(\boldsymbol{Z}'^{-\boldsymbol{C}}\right) + U_1^{C^*}$$

$$\ldots$$

$$C_Q = \alpha^*_{0,Q} + \sum_{p=1}^{P} s_{C_Q,X_p^{ep}}\left(X_p^{ep}\right) + \sum_{q=1}^{Q} s_{C_Q,C_q^{ep}}\left(C_q^{ep}\right) + s_{C_Q,V^{ep}}\left(V^{ep}\right) + s_{C_Q,\boldsymbol{Z}'^{-C}}\left(\boldsymbol{Z}'^{-\boldsymbol{C}}\right) + U_Q^{C^*}$$

$$V = \alpha_0^{**} + \sum_{p=1}^{P} s_{V,X_p^{ep}}(X_p^{ep}) + \sum_{q=1}^{Q} s_{V,C_q^{ep}}(C_q^{ep}) + s_{V,V^{ep}}(V^{ep}) + s_{V,\boldsymbol{Z}'^{-c}}(\boldsymbol{Z}'^{-\boldsymbol{C}}) + U^{V^*}$$

where $U_1^*, \ldots, U_P^*$, $U_1^{C^*}, \ldots, U_Q^{C^*}$ and $U^{V^*}$ are residuals with mean zero and a function, $s$, is considered in a way that these residuals are independent of the original covariates in the right-hand side, which can be typical in standard regression. MR-Calibrated variables are defined as

$$X_{1,MRC}^{ep} \coloneqq \hat{\gamma}_{0,1}^* + \sum_{p=1}^{P} \hat{s}_{X_1,X_p^{ep}}(X_p^{ep}) + \sum_{q=1}^{Q} \hat{s}_{X_1,C_q^{ep}}(C_q^{ep}) + \hat{s}_{X_1,V^{ep}}(V^{ep}) + \hat{s}_{X_1,\boldsymbol{Z}'^{-c}}(\boldsymbol{Z}'^{-\boldsymbol{C}})$$

$$\ldots$$

$$X_{P,MRC}^{ep} \coloneqq \hat{\gamma}_{0,P}^* + \sum_{p=1}^{P} \hat{s}_{X_P,X_p^{ep}}(X_p^{ep}) + \sum_{q=1}^{Q} \hat{s}_{X_P,C_q^{ep}}(C_q^{ep}) + \hat{s}_{X_P,V^{ep}}(V^{ep}) + \hat{s}_{X_P,\boldsymbol{Z}'^{-c}}(\boldsymbol{Z}'^{-\boldsymbol{C}})$$

$$C_{1,MRC}^{ep} \coloneqq \hat{\alpha}_{0,1}^* + \sum_{p=1}^{P} \hat{s}_{C_1,X_p^{ep}}(X_p^{ep}) + \sum_{q=1}^{Q} \hat{s}_{C_1,C_q^{ep}}(C_q^{ep}) + \hat{s}_{C_1,V^{ep}}(V^{ep}) + \hat{s}_{C_1,\boldsymbol{Z}'^{-c}}(\boldsymbol{Z}'^{-\boldsymbol{C}})$$

$$\ldots$$

$$C_{Q,MRC}^{ep} \coloneqq \hat{\alpha}_{0,Q}^* + \sum_{p=1}^{P} \hat{s}_{C_Q,X_p^{ep}}(X_p^{ep}) + \sum_{q=1}^{Q} \hat{s}_{C_Q,C_q^{ep}}(C_q^{ep}) + \hat{s}_{C_Q,V^{ep}}(V^{ep}) + \hat{s}_{C_Q,\boldsymbol{Z}'^{-c}}(\boldsymbol{Z}'^{-\boldsymbol{C}})$$

$$V_{MRC}^{ep} \coloneqq \hat{\alpha}_0^{**} + \sum_{p=1}^{P} \hat{s}_{V,X_p^{ep}}(X_p^{ep}) + \sum_{q=1}^{Q} \hat{s}_{V,C_q^{ep}}(C_q^{ep}) + \hat{s}_{V,V^{ep}}(V^{ep}) + \hat{s}_{V,\boldsymbol{Z}'^{-c}}(\boldsymbol{Z}'^{-\boldsymbol{C}})$$

To estimate an outcome model that serves as the potential outcomes model, referred to as a Q-model and then $AEE(X|\boldsymbol{z} \text{ or } \boldsymbol{z}')$ using g-computation[66,67], MR-Calibrated variables can be used because,

**Theorem 3 (Bias removal using MRC).** With Conditions 1 and 2, for *p*-th exposure, the error model is $X_p = X_{p,MRC}^{ep} + U_p^*$ and the confounder inclusion condition holds in a way that the following EUEME holds

$$Y(X_{p,MRC}^{ep} = x_{p,MRC}^{ep}) \perp\!\!\!\perp X_{p,MRC}^{ep} | \boldsymbol{z}^{**}$$

where $\boldsymbol{z}^{**} \coloneqq \{\boldsymbol{z}'^{-\boldsymbol{C}}, \boldsymbol{c}_{\boldsymbol{MRC}}^{\boldsymbol{ep}}, v_{MRC}^{ep}, \boldsymbol{x}_{\boldsymbol{MRC}}^{\boldsymbol{ep-P}}\}$ where $\boldsymbol{x}_{\boldsymbol{MRC}}^{\boldsymbol{ep-P}}$ denotes a set of values of the MR-Calibrated exposures excluding $X_{p,MRC}^{ep}$, $\boldsymbol{c}_{\boldsymbol{MRC}}^{\boldsymbol{ep}}$ denotes a set of values of the MR-Calibrated confounders, $v_{MRC}^{ep}$ denotes a value of $V_{MRC}^{ep}$, and $\boldsymbol{z}'^{-\boldsymbol{C}}$ denotes a set of values of $\boldsymbol{Z}'^{-\boldsymbol{C}}$. By applying Theorems 1 and 2,

$$E\left[\widehat{AEE}\left(X_{p,MRC}^{ep} \middle| \boldsymbol{z}^{**}\right)\right] \approx AEE(X_p | \boldsymbol{z}')$$

Appendix S5 provides the proof and remarks about $s$, multiplicative error, and marginal AEE.

## 4. Simulation Analyses

Due to space constraints, Appendix S6 provides details with Figures S2–S9, confirming the bias-reduction capability of the developed framework.

## 5. Application

### 5.1. Population and Design

During the COVID-19 pandemic, the Medical Examiner's Office Case Archive in Cook County captured almost all confirmed COVID-19 deaths that occurred in Cook County, at least from 2020 to Feb. 2021, which was used previously[68], indicating that the source population likely encompassed nearly all residents. Following the previous work[68], time-stratified case-crossover design that compares deaths and their self-control days at prespecified time points was used, which is a matched case-control design. Time-invariant factors or slowly time-varying factors within a short-time period such as age, sex, race/ethnicity, socioeconomic status, and healthy behaviors are controlled by selecting self-control days as in the same calendar year, month, and day of the week of death[68]. Different sets of COVID-19 deaths based on the type of geocoded locations at the time of an incident (e.g., home, nursing homes) and the definition of the date of death were analyzed. Table 1 provides details and previous findings[68].

### 5.2. EME and CME-Corrected Causal Effect Estimation

#### 5.2.1. Machine Learning-Based Exposure/Confounder Modeling

The previous work used less accurate modelled estimates of environmental conditions[68]. Here, we developed more spatially granular and accurate estimates for $PM_{2.5}$, $O_3$, temperature, and humidity using eXtreme Gradient Boosting (XGBoost): 500m-gridded estimates of daily mean $PM_{2.5}$, daily 8-hour moving average maximum $O_3$, daily mean temperature and relative humidity. XGBoost is a machine learning (ML) algorithm based on ensemble regression trees that address non-linear relationships among predictors to predict a dependent variable. The training dataset was constructed using ground-based monitoring data, excluding ones used to construct the validation dataset that included six locations, distributed across Cook County (Figure S10). Appendix S7 provides details.

Prediction models may be developed as a special case of MRC models. For illustration, let $X_1$, $X_2$, $X_3$, and $X_4$ denote ambient $PM_{2.5}$, $O_3$, temperature, and relative humidity. We use $X_p$ ($p = \{1,2,3,4\}$) for convenience. If models are trained in a way that

$$X_p = \theta_{1,p} + \theta_{2,p}\hat{X}_p + residuals_p$$

$$\hat{X}_p \coloneqq \hat{s}_{ML.p}(predictors)$$

where $\theta_{2,1} = \theta_{2,2} = \theta_{2,3} = \theta_{2,4} = 1$, a set of predictors include other environmental conditions and potential confounders, and $\hat{X}_p$ is independent of $residuals_{p'}$ ($p' = \{1,2,3,4\}$), then Condition 2 may hold. We refer to such $\hat{X}_p$s as MRC-Like ML-based estimates. ML algorithms can be used to obtain $\hat{s}_{ML.p}$. To increase the plausibility of Condition 2, we used several predictors including existing estimates of $PM_{2.5}$, $O_3$ (US EPA's the Community Multiscale Air Quality Modeling (CMAQ), other models, and satellite retrievals) and of temperature, humidity (PRISM group's model), day-of-the-week, and other various spatiotemporal predictors that exhibit temporal and spatial variations (e.g., greenness, traffic volumes, built-environment variation) (Table S2). The use of variables that increase predictive power in a RC model may improve efficiency[61], which applies to MRC models. Models were trained to avoid overfitting (Appendix S7).

Figures S11–14 present the agreement between MRC-Like ML-based estimates and measurements in the validation data. For $PM_{2.5}$, $O_3$, temperature, and relative humidity, overall $R^2$ was 0.73, 0.91, 0.99, and 0.83, respectively (Figure S11) with the range of site-specific $R^2$ of 0.60–0.80, 0.87–0.95, 0.95–0.99, or 0.78–0.95, respectively (Figure S12), with $\hat{\theta}_{2,p} \approx 1$. The agreement between CMAQ/PRISM-based estimates and measurements was lower than that for MRC-Like ML-based estimates (Figures S13–S14).

### 5.2.2. MRC Modeling

It would be prudent to fit MRC models using the validation data to augment EME and CME correction regarding MRC-Like ML-based estimates, we refer to as doubly

EME&CME correction. If MRC-Like ML-based estimates still deviate from Condition 2 (i.e., $\hat{X}_1 \not\approx X_{1,MRC}^{ep}$, $\hat{X}_2 \not\approx X_{2,MRC}^{ep}$, $\hat{X}_3 \not\approx C_{1,MRC}^{ep}$, and/or $\hat{X}_4 \not\approx C_{2,MRC}^{ep}$), fitting MRC models using the validation dataset would be necessary to correct for EME and CME.

For notational convenience, let $X_p^{ep}$ denote MRC-Like ML-based estimates for $X_p$. Henceforth, we overuse $\hat{X}_p$ to denote an MR-Calibrated variable: $\hat{X}_1 \equiv X_{1,MRC}^{ep}$; $\hat{X}_2 \equiv X_{2,MRC}^{ep}$, $\hat{X}_3 \equiv C_{1,MRC}^{ep}$, $\hat{X}_4 \equiv C_{2,MRC}^{ep}$. We use $L$ to denote a location within an area, $a$, that is a 500m-grid, and $t$ to denote time (i.e., day) as subscripts. MRC models were fit using the validation dataset as generalized additive models with smoothers:

$$X_{L,t,p} = \gamma_{0,p}^{*} + \sum_{l=0}^{l''} \gamma_{X_{L,t,p},X_{L,t-l,p}^{ep}}^{*} X_{L,t-l,p}^{ep} + \sum_{l=0}^{l''} \sum_{k\neq p} s_{X_{L,t,p},X_{L,t-l,k}^{ep}}\left(X_{L,t-l,k}^{ep}\right) + s_{X_{L,t,p},\boldsymbol{Z}_{L,t}^{-C}}\left(\boldsymbol{Z}_{L,t}^{-C}\right) + U_{L,t,p} \quad \text{(Eq.1)}$$

$$\hat{X}_{L,t,p} \coloneqq \hat{\gamma}_{0,p}^{*} + \sum_{l=0}^{l''} \hat{\gamma}_{X_{L,t,p},X_{L,t-l,p}^{ep}}^{*} X_{L,t-l,p}^{ep} + \sum_{l=0}^{l''} \sum_{k\neq p} \hat{s}_{X_{L,t,p},X_{L,t-l,k}^{ep}}\left(X_{a,t-l,k}^{ep}\right) + \hat{s}_{X_{L,t,p},\boldsymbol{Z}_{L,t}^{-C}}\left(\boldsymbol{Z}_{L,t}^{-C}\right).$$

Eq.1 considers potential non-linear relationships between $X_{L,t,p}$ and $X_{a,t-l,k\neq p}^{ep}$. $\boldsymbol{Z}_{L,t}^{-C}$ is a set of variables or smoothing terms that may adjust for the remaining confounders and potential differential error (See below). $l''$ indicates the maximum time lag, 21days.

For time-stratified case-crossover design, time-invariant or slowly time-varying confounders are adjusted for by matching. The matched strata can be statistically expressed as dummy variables indicating location–calendar year–month–day-of-the-week strata because all individuals at the identical location at day *t* share the same value of environmental conditions, which can be reasonably approximated to the combination of a spatial smoothing function based on coordinates of $L$, dummy variables for day-of-the-week, and a temporal smoothing function based on day-of-the-year variable in our setting, based on the mathematical relationship between time-stratified case-crossover analyses and time-series analyses[69,70], including case-time-series analysis[71]. These smoothing functions well-captured temporal and spatial patterns in $X_{L,t,p}$. Following the previous work[68], COVID-19 infectivity as a potential confounder adjusted for in an outcome model using lag-structure of ZIP Code-level COVID-19 confirmed cases and adjacent ZIP codes was also included. Consequently, the fitted MRC models did not have temporal autocorrelation of residuals. No correlations of $\hat{X}_{L,t,p}$ with the residuals were found.

#### 5.2.3. Estimation of Conditional AEE and Premature Deaths

Let $\hat{X}_{i,L,t,p}$ denote $\hat{X}_{L,t,p}$ linked to the health data using the geocoded location for individuals. We used conditional logistic regression, which is standard for case-crossover analyses. Covariates included constrained distributed lags of $\hat{X}_{i,L,t,1}$ and $\hat{X}_{i,L,t,2}$ and constrained distributed lag non-linear terms[5] of $\hat{X}_{i,L,t,3}$ and $\hat{X}_{i,L,t,4}$ and constrained distributed non-linear terms of ZIP Code-level COVID-19 infection rates, following the previous work[68].

Our health data did not include the mobility patterns of individuals, meaning that the linkage between MR-Calibrated variables and the health data may have introduced error if individuals had moved far around the geocoded location. Additionally, as we used MR-Calibrated variables for ambient $PM_{2.5}$, $O_3$, temperature, and humidity, the difference between personal exposure and ambient levels may have existed.

Collectively, we acknowledge that even after MRC modeling, EME and CME may not be fully eliminated. To address these, we

1) factorized potential error sources when modelled estimates for ambient environmental conditions are used (See Appendix S8 and Figure S15), how these were addressed collectively, and how potential differential error from sources were adjusted for in our analysis (See Appendix S8).
2) sub-optimally identified potential bias due to unmeasured mobility (See Appendix S9 with Figures S16–17), although many short-term air pollution epidemiological studies have not accounted for mobility.

We estimated the conditional AEE using the OR (odds ratio), RD, and RR scales and the number of premature deaths from COVID-19 from Mar. 2020 to Feb. 2021 that would have been prevented had $PM_{2.5}$ and $O_3$ levels been reduced to average background levels over the US[72] (~5μg/$m^3$, ~51ppb, respectively), which also corresponds to WHO air quality recommendations, and to observed low levels in Cook County from 2020 to Feb. 2021 (~2μg/$m^3$ for $PM_{2.5}$, monthly minimum for $O_3$). Appendix S10 explains methods. The term, prematurity, is used to indicate that the timing of death is shifted earlier, as death is an inevitable health endpoint.

The estimated conditional probability of $Y$=1 in each matched stratum that is related to conditional likelihood for each stratum was concentrated between 0.2 and 0.8, implying good approximation condition. This does not mean the probability of $Y$=1 for each individual (See the difference between these two in Appendix S10).

Bootstrapping[10,43] was used to obtain central estimates and 95% CIs. Parametric bootstrapping was chosen over non-parametric bootstrapping to generate bootstrapped MR-Calibrated $PM_{2.5}$ and $O_3$ variables from MRC models because MRC models include smoothing functions. Conditional logistic regression was separately fit over 500 bootstrapping samples with bootstrapped MR-Calibrated variables for $PM_{2.5}$ (or $O_3$), while the other MR-Calibrated variables for temperature, humidity, and $O_3$ (or $PM_{2.5}$) were co-adjusted for.

For comparison, AEE was also estimated using MRC-Like ML-based estimates or CMAQ/PRISM-based estimates of $PM_{2.5}$, $O_3$, temperature, and humidity.

### 5.3. Results

Figure 1 shows spatial distributions of MRC-Like ML-based estimates of $PM_{2.5}$, $O_3$, temperature, and relative humidity. Figure 2 presents scatter plots for MRC-Like ML-based estimates of $PM_{2.5}$ and $O_3$ levels linked to case and self-controls. Many of the linked exposures were below the National Ambient Air Quality Standards for 24-hour $PM_{2.5}$ and 8-hour $O_3$ levels.

Figure 3 presents OR estimates based on different exposure/confounder variables. Table S3 provides RR and RD estimates based on MR-Calibrated variables; RR estimates nearly equal OR estimates. AEE estimates based on MR-Calibrated exposure/confounder estimates, with the most precise 95% CIs, were comparable to those based on MRC-Like ML-based estimates, showing positive associations of both $PM_{2.5}$ and $O_3$ with COVID-19 mortality, while some of the CMAQ-based AEE estimates did not present positive associations. For example, for one interquartile range increase in 3-week $PM_{2.5}$ and $O_3$ exposures (5.6μg/$m^3$ and 22.7ppb, respectively) were associated with a 133.9% (95% CI: 75.8, 197.2) increase and a 135.3% (68.4, 233.0) increase in the mortality risk, respectively (Home/Work & Diff. ≤7d in Figure 3). The $PM_{2.5}$ results are consistent but the $O_3$ results contradict the previous findings (Table 1). The effect of shorter-term exposures adjusted for the remaining exposure period within 3-week was also estimated: 1-day (Lag0), 3-day (Lag0-2), 1-week (Lag0-7), 2-week (Lag0-14) in Figure 3.

Among 7,507 deaths from COVID-19, the number of preventable premature deaths by reducing $PM_{2.5}$ and $O_3$ levels to observed low levels or US-average natural background levels was 6036 (5326, 6489) or 3366 (2538, 3969), respectively (Table 2).

## 6. Discussion

We developed a new framework for estimating causal effects by addressing both EME and CME. Our applied analysis with advanced prediction modeling for air pollution and meteorological factors and error correction found positive association between 3-week exposure to $O_3$ and COVID-19 mortality, which was not found in the previous work. We could not consider individual mobility in linking environmental data and health data; however, if unmeasured mobility introduced substantial bias, it likely biased the results toward the null (See Appendix S9).

The estimates of preventable premature mortality may not necessarily imply that these individuals would have ultimately survived but implies that the timing of deaths would have been delayed had $PM_{2.5}$ and $O_3$ exposures been eliminated. The extent of this delay would have varied across individuals due to their inherent heterogeneity[73-79].

This research opens avenues for future research. First, studies could examine whether approximation conditions can be relaxed for research settings that may deviate from the conditions identified here. For example, other error correction and effect estimation methods[45,58,80,81] may merit investigation due to differing statistical properties. Developing efficient computational algorithms could be beneficial. Second, studies could investigate whether the confounder inclusion condition can be relaxed. Unless bias amplification in mixture/multi-pollutant analyses is concerning[1], which may occur when unmeasured confounding exists (a violation of Condition 1), all confounder variables should be considered. In practice, constructing validation data for high-dimensional exposures and confounders, may be challenging for reasons. Furthermore, a minimal sufficient set of covariates for confounding adjustment[82] and how to consider this may differ. For example, effect estimation may rely on spatial/between-individual contrast, temporal/within-individual contrast, or both. Depending on study designs, accordingly techniques to address unmeasured confounders may also differ. Future work should address other real-world constraints.

**Table 1. Percent increase (95% CI) in the risk of deaths from COVID-19 by one interquartile range (IQR) increase in cumulative 3-week exposure to $PM_{2.5}$ (per 5.2μg/$m^3$) and $O_3$ (per 16 ppb) reported by Kim et al. (2022)[68].**

| Type of analysis† | Pollutants | *n*‡ | All Locations§ | *n*‡ | Home/Work§ | *n*‡ | Home/Work & Diff. ≤7d§ |
|---|---|---|---|---|---|---|---|
| Date of Death | $PM_{2.5}$ | 7462 | 69.6 (34.6, 113.8) | 4881 | 60.4 (21.6, 111.7) | 1941 | 78.5 (14.1, 179.4) |
| | $O_3$ | 4434 | −9.0 (−49.6, 64.3) | 2591 | −20.7 (−61.7, 64.4) | 1077 | −16.2 (−75.7, 188.3) |
| Date of Incident | $PM_{2.5}$ | 7507 | 123.9 (78.3, 181.1) | 4923 | 109.4 (58.8, 176.3) | 1943 | 138.6 (52.6, 273.2) |
| | $O_3$ | 4647 | −16.2 (−54.2, 53.3) | 2767 | −21.0 (−63.1, 69.3) | 1089 | −17.4 (−76.5, 190.7) |

†"Date of Death" is the date the deceased was pronounced dead. "Date of Incident" is the date of a relevant event that preceded the date that the deceased was pronounced dead. For illustration, suppose that a person started experiencing severe symptoms at home on 16 March ("Date of Incident") and was hospitalized. He/she tested positive on 19 March. His/her heart stopped beating at a hospital on 24 March, and he/she was eventually pronounced dead on the same day ("Date of Death"). See Kim et al (2022)[68].

‡The number of deaths from COVID-19 used for analyses.

§"All locations" refers to findings from analyses using all types of the geocoded locations; "Home/Work" refers to findings from analyses restricted to only deaths where the geocoded locations at the "Date of Incident" were identified as residential or occupational addresses; "Home/Work & Diff. ≤7d" refers to findings from analyses restricted to only deaths where the geocoded locations at the "Date of Incident" were identified as residential or occupational addresses with the difference between "Date of Incident" and "Date of Death" ≤7days, which was used to exclude deaths from COVID-19 due to adverse outcomes other than air pollution exposures. See Kim et al (2022)[68].

**Table 2. Estimated number of preventable premature deaths from COVID-19**

| Target† | Number (95% CI) |
|---|---|
| (1) Reduction to observed low levels within and around Cook County | 6,036 (5,326, 6,489) |
| (1-1) Only $PM_{2.5}$ reduction (to 2μg/$m^3$) | 4,692 (3,675, 5,389) |
| (1-2) Only $O_3$ reduction (to monthly low levels) | 3,348 (2,397, 4,188) |
| (2) Reduction to US nationwide average natural background levels | 3,366 (2,538, 3,969) |
| (2-1) Only $PM_{2.5}$ reduction (to 5μg/$m^3$) | 3,277 (2,462, 3,915) |
| (2-2) Only $O_3$ reduction (to 51ppb) | 131 (96, 165) |

†US nationwide average natural background levels are consistent with World Health Organization's Air Quality Guidelines 2021 ($PM_{2.5}$ annual: 5μg/$m^3$; $O_3$ 8-hour: 100μg/$m^3$). Because the contribution of natural sources to air pollutant levels varies regionally, observed low levels may differ. Monthly observed low levels for $O_3$ include 7ppb (Jan.); 18ppb (Feb.); 8ppb (Mar.); 16ppb (Apr.); 25ppb (May); 27ppb (Jun.); 26ppb (Jul.); 17ppb (Aug.); 12ppb (Sep.); 5ppb (Oct.); 10ppb (Nov.); and 3ppb (Dec.)

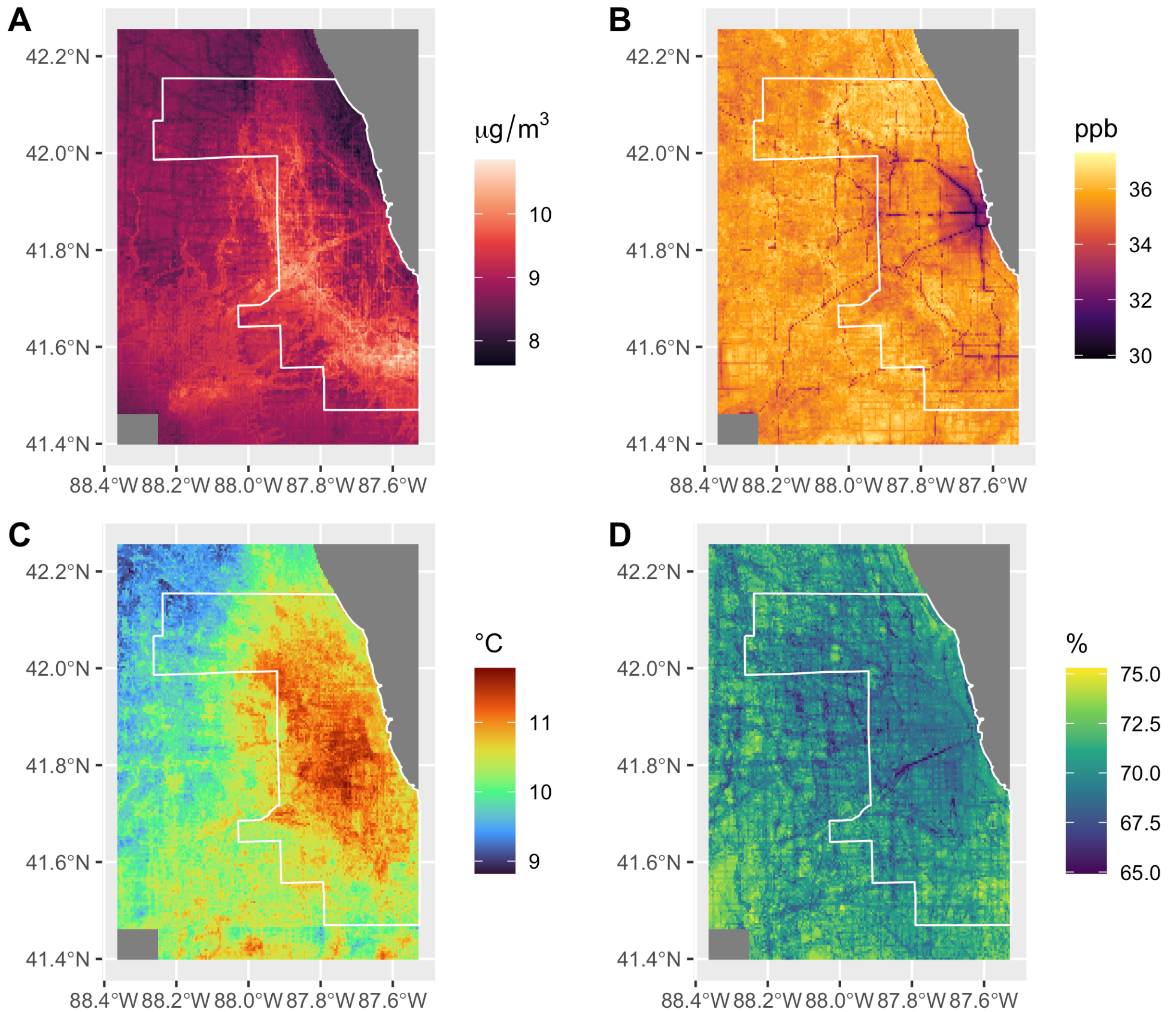

**Figure 1. MRC-Like ML-based 500m gridded estimates of $PM_{2.5}$ (A), $O_3$ (B), temperature (C), and relative humidity (D):**

Note. These averages are from Mar. 2020 to Feb. 2021. The white solid line indicates the boundary of Cook County, IL. The gray area in the top right is Lake Michigan. $PM_{2.5}$ levels were high around major highways, while $O_3$ levels were low, possibly due to $NO_x$ scavenging. In a large metropolitan city, $O_3$ levels may be lower in downtown areas with high traffic volumes compared to outskirts of the city or suburban areas due to the difference between $NO_x$-saturated conditions and $NO_x$-limited conditions with $NO_x$ scavenging and the difference in biogenic volatile organic compound emissions[83,84]. Temperature levels are higher in central areas of Cook County, possibly due to urban heat island effects. Humidity levels may be lower in roads than residential areas that include soils and trees because paved surfaces lack moisture storage and evapotranspiration and because heat radiated from roads may increase evaporation.

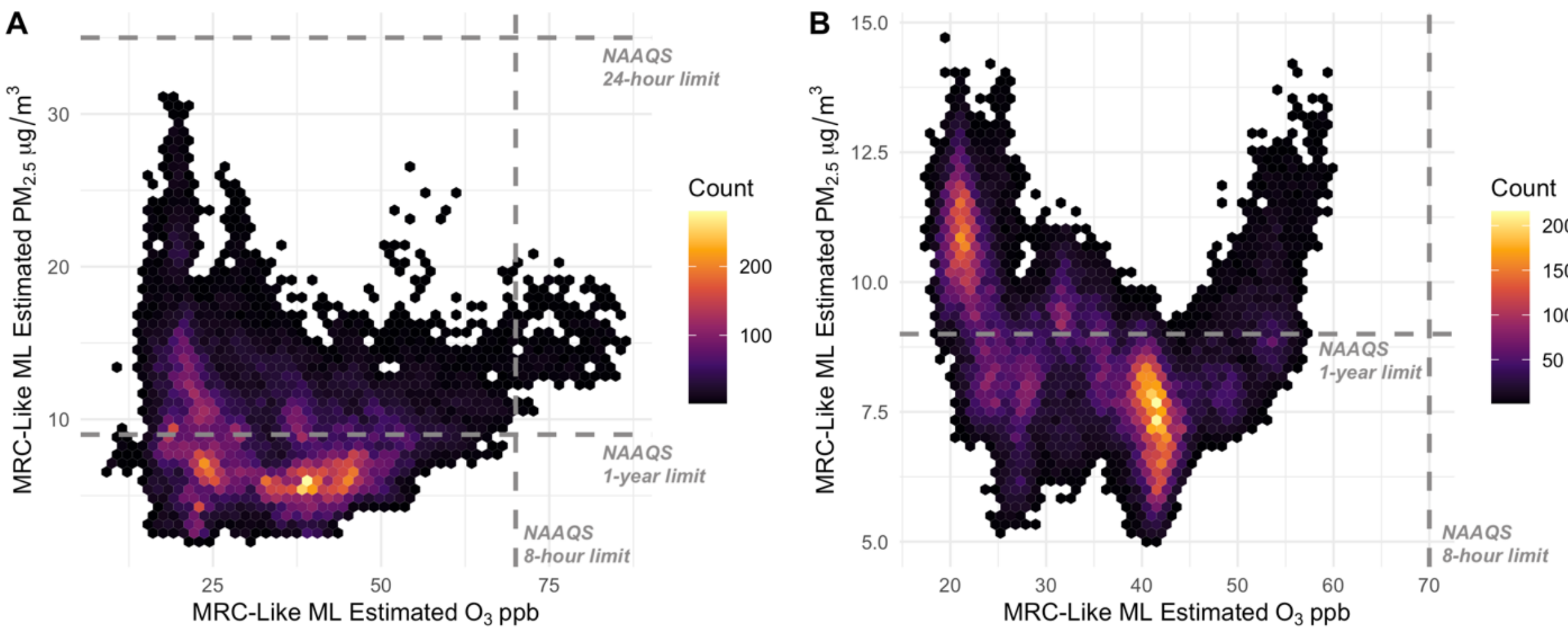

**Figure 2. Scatter plots of MRC-Like ML-based estimates of $PM_{2.5}$ and $O_3$ levels:**

**Lag0 (1-day) exposure (A); Cumulative Lag0-21 (3-week) exposure (B)**

Note. NAAQS stands for the National Ambient Air Quality Standards in the US (For $PM_{2.5}$, 24-hour limit is 35μg/$m^3$ and 1-year limit is 9μg/$m^3$; For $O_3$, 8-hour limit is 70ppb). A total number of cases–self-controls–days is 32,914.

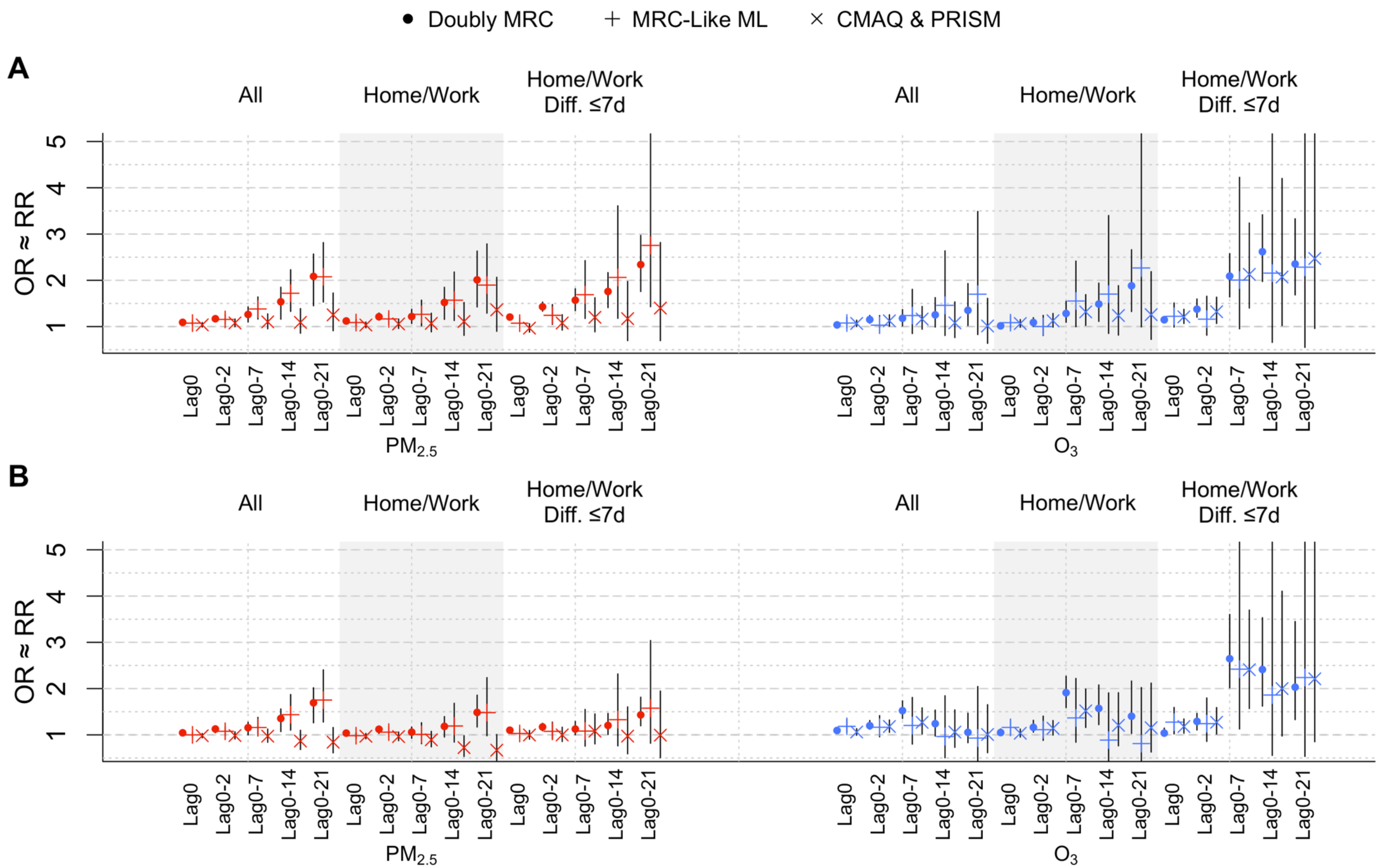

**Figure 3. Conditional AEE estimates in the OR scale for one interquartile increase in cumulative lags of $PM_{2.5}$ and $O_3$ levels by MR-Calibrated variables, MRC-Like ML-based estimates, and CMAQ & PRISM-based estimates. A. 'Date of Incident' based analyses. B. 'Date of Death' based analyses.**

Note 1. See a footnote of Table 1 for the definition of "Date of Incident" and "Date of Death".

Note 2. "Doubly MRC" refers to findings from analyses based on MR-Calibrated exposure/confounder estimates ($PM_{2.5}$, $O_3$, temperature, and relative humidity); "MRC-Like ML" refers to findings from analyses based on MRC-Like ML-based exposure/confounder estimates (i.e., no additional correction by MRC models); and "CMAQ & PRISM" refers to findings from analyses based on CMAQ estimates of $PM_{2.5}$ and $O_3$ and PRISM estimates of temperature and relative humidity.

**Supplementary Materials**

**Improving Causal Inference with Measurement Errors in Exposures and Confounders: A New Method and Its Application to Air Pollution Exposure Assessment and Epidemiology**

**List of Contents**

**Appendix S1. Illustrations of EME and CME and Errors from Prediction Models.**

**S1.1. EME and Error Types**

We begin with introducing Non-Berkson (NB) error and Berkson error[1] with directed acyclic graphs (DAGs). Let $X$ denote the true continuous exposure and $X^{ep}$ denote an error-prone $X$ meaning that $X$ is measured or estimated with error. NB error may arise when an instrument applied to measure or estimate $X$ is error-prone for various reasons: random noise, sensor error, recall bias, the fate of $X^{ep}$ (e.g., biomarkers and recent measurements as a surrogate for past exposures or past (baseline) measurements as a surrogate for cumulative exposures (over a follow-up)). NB error (Figure S1A) may be encoded as

$$X^{ep} = g(X) + U$$

where $U$ is completely random with $E[U] = 0$; $g$ is a transformation function. Examples include classical or linear errors ($X^{ep} = X + U$ or $X^{ep} = \eta_0 + \eta_1 X + U$). Error may be multiplicative, such as

$$m(X^{ep}) = g(m(X)) + U$$

where $m$ is other transformation function (e.g., logarithm).

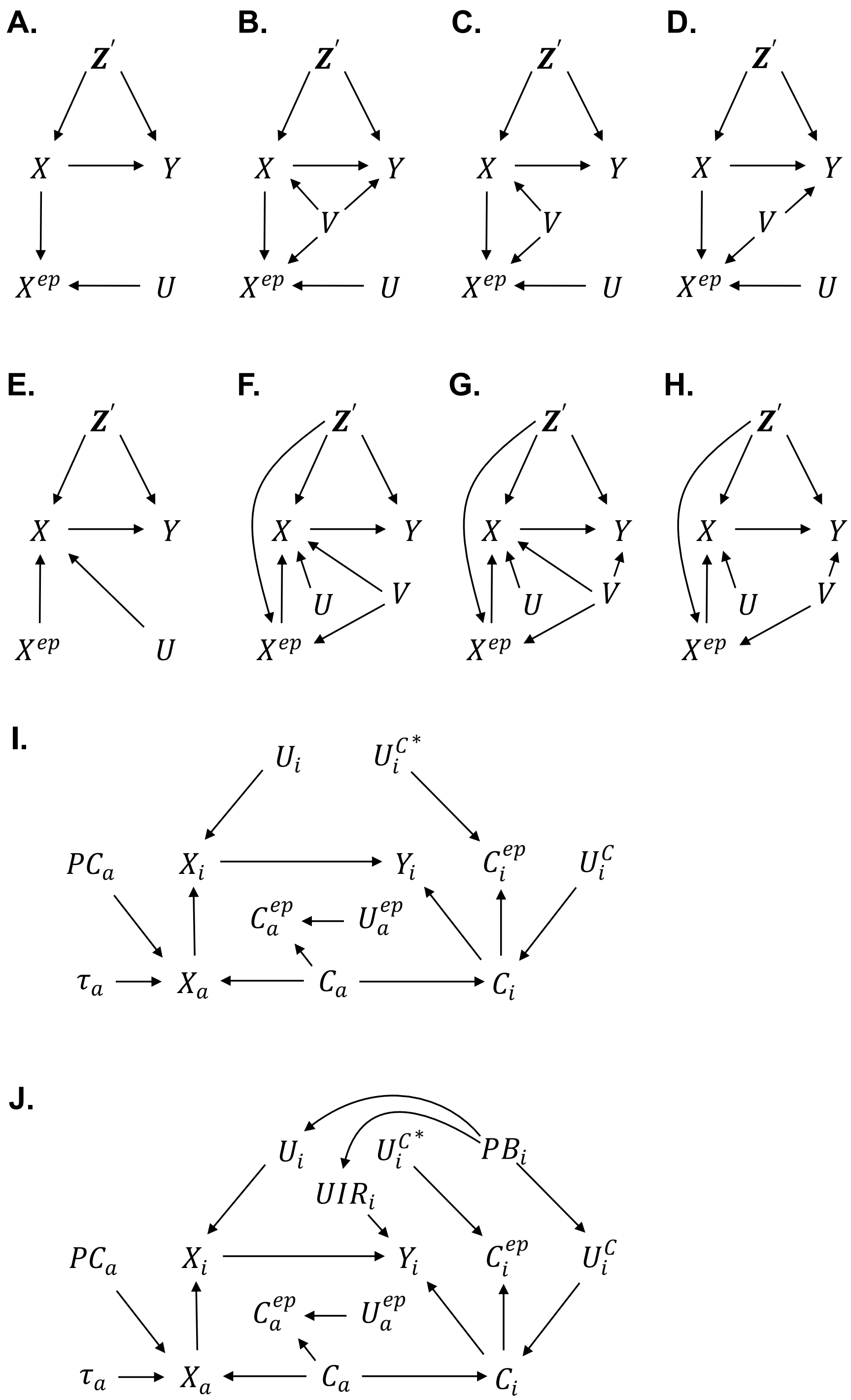

**Figure S1. Directed acyclic graphs for non-Berkson error (A–D), for Berkson error (E–H), for other examples (I–J).**

EME may be differential. Differential error means that the error is not independent of a health outcome, $Y$. When EME varies systematically by individual, group, geography, or time, differential error may arise if these characteristics are related to the distribution of $Y$. This is unlikely uncommon in environmental exposures, as they are generally characterized by spatiotemporal patterns, with exposure pathways, routes, and environmental fate varying according to individuals' behavioral patterns or biological and ecological processes, while measurement quality can differ by individual, location, and over time. Let $V$ be a variable indicating what makes the error differential (Figures S1B–D).

$$X^{ep} = g(X) + V + U$$

$$m(X^{ep}) = m\big(g(X)\big) + V + U$$

The total error $\mathcal{E} = V + U$ may be differential. $V$ may render $\mathcal{E}$ correlated with $Y$ directly and/or through $X$, including autocorrelated errors that are correlated with $X$ or $Y$. For example, $X^{ep}$ may be a biomarker that decays over time, $g(X)$ is inhalation to an environmental chemical (in the past), $V$ denotes chronological age or time that encodes the decay[2], and $Y$ may represent developmental outcomes (e.g., lung function, neurodevelopment) (Figure S1D). If $V$ also affects time-varying exposure, $X$, then Figure S1B is possible; and if the risk of $Y$ does not significantly change over time, then Figure S1C is possible. $X^{ep}$ may also be a measurement from a low-cost sensor and $V$ may denote operation time reflecting quality degradation (Figure S1D) or calendar or follow-up time, which may capture both the operation time and the seasonality of an actual environmental pollutant level (Figure S1C) and of $Y$ (Figure S1B). $X^{ep}$ may also refer to a recent measurement, where $X$ may represent the past exposure of interest and the recent exposure is temporally autocorrelated with the past exposure of interest, and $V$ may represent a causal factor that affects this time-varying exposure, (e.g., human behaviors, environmental polluting patterns). In Figure 1B, $V$ is also a confounder[3]. In Figure S1D, $V$ is a confounder when using $X^{ep}$ for effect estimation but not a confounder when using $X$ for effect estimation.

Berkson error may arise when using aggregated or ecological variables. For example, an aggregated/area-level exposure variable may be measured or estimated in large population studies, including measurements at a fixed monitor, modelled estimates for pollutants, or metrics based on distance to polluting sources. To describe this, for example,

$$X_{i,L,t} = g(X_{a,t}) + U_{i,L,t}$$

where $a$ denotes an area, $L$ denotes a location, at which an individual, $i$, stays, within $a$, and $t$ denotes time $t$. Other disciplines with Berkson error can be found elsewhere[4-9]. We refer to

$$X = X^{ep} + U$$

as the pure Berkson error. The term *pure* is used to emphasize that $U$ is completely random, implying the independence between $U$ and $X^{ep}$, and the linear coefficient of 1. Berkson error models (Figures S1E–G) may be generally expressed as

$$X = g(X^{ep}) + V + U$$

$$m(X) = g\left(m(X^{ep})\right) + V + U$$

The total error $\mathcal{E} = V + U$ may be differential if $V$ affects $X^{ep}$ and $Y$ directly or through $X$ ($X^{ep} - V - X - Y$ or $X^{ep} - V - Y$ pathways) (Figures S1F-G). Figure 1H is an exception that the total error is $U$, such as $X = g(X^{ep}) + U$, but $X^{ep}$ is affected by $V$ such that $X^{ep}$ is subject to differential error. In Figure 1H, $V$ may alternatively be regarded as a confounder only when using $X^{ep}$. Generally, $V$ can be seen as a confounder when using $X$ or $X^{ep}$ ($X - X^{ep} - V - Y$), depending on a DAG. For example, $V$ may represent area-level ambient temperature that may affect area-level ambient air pollution[10], $X^{ep}$, that may affect personal exposure, $X$. This example can be extended to Figures S1F, or S1G, depending on the study. Temperature may indicate behavioral pattern changes (e.g., people tend to spend more time outdoors in warm or cool weather) that affect $X$ (Figure S1F) and for certain health outcomes, temperature may be a risk factor (Figure S1G). $X^{ep}$ may also represent area-level pollutant concentration that can be increased by emissions from sources or ecological processes in environment. $V$ may encode these.

**S1.2. CME**

In addition to EME, studies may rely on confounder variables that are also subject to measurement error. Using error-prone confounder variables for confounder adjustment may result in residual confounding. The mechanisms of residual confounding are various. Even when EME is purely Berksonian, bias may arise due to unmeasured confounding[11], which can be extended to residual confounding due to CME. Some investigators have noted that the net bias due to EME and unmeasured confounding may cancel out in certain conditions[12]. Extending this concept to net bias from both EME and CME may be tempting, but conditions would be specific but not necessarily practical. Specifically, if a non-linear covariate-response relationship exists or multiple covariates are of interest, the identification of such conditions could be tedious, and the conditions would not be generalizable but study-specific.

As shown below, CME may introduce residual confounding but may not do so, which may be relevant in the context of air pollution epidemiology. Suppose a DAG in Figure S1I that may describe

$$X_i = X_a + U_i$$

$$X_a = f(C_a, PC_a) + \tau_a$$

$$C_a^{ep} = C_a + U_a^C$$

$$C_i = C_a + U_i^C$$

$$C_i^{ep} = C_i + U_i^{C^*}$$

$$Y_i = \alpha + \beta_X X_i + \beta_C C_i$$

where $X$ is $O_3$, $C$ is high temperature as a confounder, $PC$ includes volatile organic compounds and nitrogen oxides as precursors and $f(C_a, PC_a)$ describes $O_3$ formation through photochemical reactions. $a$ denotes an area and $i$ denotes an individual. When estimating the health effect of $X$ using $X_a$, adjustment for $C_a^{ep}$ or $C_i^{ep}$ cannot fully block the backdoor path ($X_a - C_a - C_i - Y_i$), meaning residual confounding. However, in certain circumstances, the use of variables subject to (pure) Berkson EME and (pure) CME may be preferred over error-free variables to avoid bias from EME and CME, and unmeasured confounding. The forementioned example demonstrates that adjustment for $C_a$ is sufficient when using $X_a$ for effect estimation, meaning no residual confounding. Also, by extending the point made by Weisskopf and Webster (2017)[13] — that “the set of possible confounding variables for the ambient concentration estimation cannot be larger than that for personal exposure" — to the EME perspective, we can consider a situation like Figure S1J, where $U_i$ and $U_i^C$ depend on personal behaviors, $PB_i$ that also affects unmeasured individual risk factors, $UIR_i$, that are confounders when using $X_i$ for effect estimation. Here, adjustment for $C_a$, would still be sufficient if $X_a$ is used to estimate the effect of $X$ because $UIR_i$ are not confounders when using $X_a$. In the case of Weisskopf and Webster (2017)[13], $X_a$ may denote the ambient concentration and $X_i$ may denote personal exposure. This example demonstrates that CME may not necessarily result in bias.

### S1.3. Modelled/Predicted Exposure and Confounder Estimates

The resulting error from modelled estimates in health outcome modeling may not equal to the error identified in exposure/confounder modeling studies. For illustration, let $\boldsymbol{Z'}$ denote a set of confounders. Suppose that modeled exposure or confounder ($C$) estimates are derived from an exposure and confounder prediction models without $\boldsymbol{Z'}$,

$$X = \vartheta_0 + s(predictors\ without\ \boldsymbol{Z'}) + \epsilon$$

$$X^{ep} \coloneqq \hat{\vartheta}_0 + \hat{s}(predictors\ without\ \boldsymbol{Z'})$$

$$C = \vartheta_0^C + s^C(predictors\ without\ \boldsymbol{Z'}) + \epsilon^C$$

$$C^{ep} \coloneqq \hat{\vartheta}_0^C + \hat{s}^C(predictors\ without\ \boldsymbol{Z'})$$

where $\epsilon$ and $\epsilon^C$ are residuals with mean zero.

So, the error in exposure/confounder modeling may *appear* to be purely Berksonian as

$$X = X^{ep} + \epsilon$$

$$C = C^{ep} + \epsilon^C$$

However, epidemiological analyses are often conditional on confounders, meaning that the resulting error structure in outcome modeling that comes from these prediction models is dependent on confounders, as noted by Haber et al (2021)[11]. For illustration, let $\boldsymbol{Z'}^{-C}$ denote a set of confounders without $C$. The resulting EME and CME in an outcome regression model with $X^{ep}$, $C$ or $C^{ep}$, and $\boldsymbol{Z'}^{-C}$ as covariates would follow:

$X|\boldsymbol{Z'} = g_1(X^{ep}|\boldsymbol{Z'}) + \epsilon|\boldsymbol{Z'}$ (Eq. S1-1) (if $C$ is measured without error)

$X|\boldsymbol{Z'}^{-C}, C^{ep} = g_2(X^{ep}|\boldsymbol{Z'}^{-C}, C^{ep}) + \epsilon|\boldsymbol{Z'}^{-C}, C^{ep}$ (Eq. S1-2) (if $C^{ep}$ is used)

$C|\boldsymbol{Z'}^{-C}, X^{ep} = g_3(C^{ep}|\boldsymbol{Z'}^{-C}, X^{ep}) + \epsilon^{C}|\boldsymbol{Z'}^{-C}, X^{ep}$ (Eq. S1-3) (if $C^{ep}$ is used)

where $variable1|variable2$ means the residuals of $variable1$ after adjusting for $variable2$. All these equations are in the form of Berkson type but may not be purely Berskonian, if model predictors are correlated with $\boldsymbol{Z'}$ for Eq. S1-1 and Eq. S1-2, or correlated with $\boldsymbol{Z'}^{-C}, X$ for Eq. 1-3. $g_1$, $g_2$, and $g_3$ may not necessarily be the identity function. If $g_1$ and $g_2$ are not the identity function, then bias due to EME would arise. Even if they are the identity function, $X^{ep}|\boldsymbol{Z'}$ may be correlated with $\epsilon|\boldsymbol{Z'}$ or $X^{ep}|\boldsymbol{Z'}^{-C}, C^{ep}$ may be correlated with $\epsilon|\boldsymbol{Z'}^{-C}, C^{ep}$, meaning that bias due to EME may arise. $X^{ep}|\boldsymbol{Z'}^{-C}, C^{ep}$ may be correlated with $\epsilon^{C}|\boldsymbol{Z'}^{-C}, X^{ep}$, meaning that residual confounding due to CME may arise. This statistical property makes it complicated to understand the magnitude and direction of bias when EME and CME as model errors are not addressed—even if epidemiological studies are based on the same exposure and confounder models—because the list of confounders and the joint distribution of confounders in each epidemiological study may differ. To avoid this issue, confounders should be included in the prediction models, which is analogous to the confounder inclusion condition in RC.

**Appendix S2. Substantiation of Point 3 in Section 2 of the Main Manuscript.**

*Example 1*

Suppose that the true potential outcomes model/outcome data generating model is

$$q\big(Y(X)\big) = \alpha + \beta X + j(\boldsymbol{Z}') + W$$

where $q$ is the non-linear link function. For example, log-linear models may be used for count data. Even if investigators correctly assumed this model form and $X^{ep}$ is purely Berksonian,

$$E_W\big[E_U[\,Y(X = e_T + \delta + u_1)|X^{ep} = e_T + \delta, \boldsymbol{z}']\big] - E_W\big[E_U[Y(X = e_T + u_2)|X^{ep} = e_T, \boldsymbol{z}']\big]$$

$$= E_W\big[E_U[\,q^{-1}(\alpha + \beta(e_T + \delta + u_1) + j(\boldsymbol{z}') + w)|X^{ep} = e_T + \delta, \boldsymbol{z}']\big]$$

$$- E_W\big[E_U[q^{-1}(\alpha + \beta(e_T + u_2) + j(\boldsymbol{z}') + w)|X^{ep} = e_T, \boldsymbol{z}']\big]$$

$$\neq E_W[q^{-1}(\alpha + \beta(e_T + \delta) + j(\boldsymbol{z}') + w)|\boldsymbol{z}'] - E_W[q^{-1}(\alpha + \beta(e_T) + j(\boldsymbol{z}') + w)|\boldsymbol{z}']$$

$$= AEE(X|\boldsymbol{z}')$$

because $X$ is unknown and the effect is not strictly linear. The inequality arises due to $E[q^{-1}(\beta x + \beta u)] \neq E[q^{-1}(\beta x) + q^{-1}(\beta u)]$.

*Example 2*

Even for linear outcome models if $\beta X$ is replaced with $h(X)$ where $h$ is a non-linear function such as

$$Y(X) = \alpha + h(X) + j(\boldsymbol{Z}') + W$$

Suppose that investigators correctly assumed this model form and $X^{ep}$ is purely Berksonian. Then,

$$E_W\big[E_U[\alpha + h(e_T + \delta + u_1) + j(\boldsymbol{z}') + w|X^{ep} = e_T + \delta, \boldsymbol{z}']\big]$$

$$- E_W\big[E_U[\alpha + h(e_T + u_2) + j(\boldsymbol{z}') + w|X^{ep} = e_T, \boldsymbol{z}']\big]$$

$= E_U[h(e_T + \delta + u_1) - h(e_T + u_2)]$

$\neq h(X = e_T + \delta) - h(X = e_T)$

$= E_W[\alpha + h(e_T + \delta) + j(\boldsymbol{z}') + w|\boldsymbol{z}'] - E_W[\alpha + h(e_T) + j(\boldsymbol{z}') + w|\boldsymbol{z}']$

$= AEE(X|\boldsymbol{z}')$

The inequality arises due to $E[h(e_T + u)] \neq E[h(e_T) + h(u)] = h(e_T) + E[h(u)]$

*Example 3*

Suppose the investigator correctly assumed the model form as

$$logit(P(Y(X))) = \alpha + \beta X + j(\boldsymbol{Z}') + W$$

which may be used for binary outcomes. The inequality arises because the logistic function is non-linear:

$$E_W\left[E_U\left[\frac{1}{1 + \exp\left(-(\alpha + \beta(e_T + \delta + u_1) + j(\boldsymbol{Z}') + W)\right)}\middle| X^{ep} = e_T + \delta, \boldsymbol{z}'\right]\right]$$

$$\neq E_W\left[\frac{1}{1 + \exp\left(-(\alpha + \beta(e_T + \delta) + j(\boldsymbol{Z}') + W)\right)}\middle|\boldsymbol{z}'\right]$$

Note that for non-linear functions small numerical differences in a variable may have significant consequences unless the function's range falls within a domain where such differences have minimal impact. This consequence depends on a functional form.

**Appendix S3. An Illustration of An Extension of the Potential Outcomes Framework For EME**

For illustration, Table S1 shows 20 individuals in a hypothetical population, where, for illustrational simplicity, the potential outcomes model is $Y(x) = 0.1x$, without $W$ and $\boldsymbol{Z}'$, given $x \in X = \{8,9,10\}$ and the error model is $X^{ep} = X + U'$ where $u' \in U' = \{-1,0,1\}$ as random error. Now, suppose that $X$ is unknown and investigators want to identify the outcome over exposure levels using $X^{ep}$. For the value of $X^{ep} = 9$, the value of $Y_{X^{ep}=9}$ and $Y(X^{ep} = 9)$ can be 0.8, 0.9, or 1.0, because in this population, a group of individuals with $X^{ep} = 9$ includes some individuals having $Y(X = 8)$, some having $Y(X = 9)$, and some having $Y(X = 10)$. Generally, $E[Y(X^{ep} = 9)]$ may not equal $E[Y(X = 9)] = 0.9$. These two are equal if and only if the joint distribution of $X$ and $U'$ makes them so. With the consistency assumption, investigators can identify $E[Y(X^{ep} = 9)]$, not $E[Y(X = 9)]$ from the observed data without $X$. The causal effect of 0.1 may not be identifiable through, $E[Y(X^{ep} = 10) - Y(X^{ep} = 9)] = E[Y_{X^{ep}=10}] - E[Y_{X^{ep}=9}]$. In this case, underestimation may arise as in classical EME.

**Table S1. A hypothetical data that illustrates identifiability issues if $X$ is measured with error. See Note.**

| ID | | | $Y$ | | $Y(X)$ | |
|---|---|---|---|---|---|---|
| | $X$ (a) | $X^{ep}$ (b) | $Y_X$ (c: if $X$ is measured) | $Y_{X^{ep}}$ (d: if $X$ is measured with error) | $Y(X)$ (e: if $X$ is measured) | $Y(X^{ep})$ (f: if $X$ is measured with error) |
| 1 | 8 | 7 | 0.8 | 0.8 | 0.8 | 0.8 |
| 2 | 8 | 8 | 0.8 | 0.8 | 0.8 | 0.8 |
| 3 | 9 | 8 | 0.9 | 0.9 | 0.9 | 0.9 |
| 4 | 9 | 8 | 0.9 | 0.9 | 0.9 | 0.9 |
| 5 | 8 | 8 | 0.8 | 0.8 | 0.8 | 0.8 |
| 6 | 10 | 9 | 1.0 | 1.0 | 1.0 | 1.0 |
| 7 | 9 | 9 | 0.9 | 0.9 | 0.9 | 0.9 |
| 8 | 9 | 9 | 0.9 | 0.9 | 0.9 | 0.9 |
| 9 | 9 | 9 | 0.9 | 0.9 | 0.9 | 0.9 |
| 10 | 9 | 9 | 0.9 | 0.9 | 0.9 | 0.9 |
| 11 | 9 | 9 | 0.9 | 0.9 | 0.9 | 0.9 |
| 12 | 8 | 9 | 0.8 | 0.8 | 0.8 | 0.8 |
| 13 | 9 | 9 | 0.9 | 0.9 | 0.9 | 0.9 |
| 14 | 9 | 9 | 0.9 | 0.9 | 0.9 | 0.9 |
| 15 | 10 | 10 | 1.0 | 1.0 | 1.0 | 1.0 |
| 16 | 10 | 10 | 1.0 | 1.0 | 1.0 | 1.0 |
| 17 | 9 | 10 | 0.9 | 0.9 | 0.9 | 0.9 |
| 18 | 10 | 10 | 1.0 | 1.0 | 1.0 | 1.0 |
| 19 | 9 | 10 | 0.9 | 0.9 | 0.9 | 0.9 |
| 20 | 10 | 11 | 1.0 | 1.0 | 1.0 | 1.0 |

**Note.** Values in columns *c*–*f* for each individual are identical as they indicate the same outcome of each individual. In practice, they may be labeled differently depending on the exposure variable in use. This different labeling has consequences in data analysis as shown earlier.

**Appendix S4. Conditions for Minimal Bias (Proof for Theorems 1 and 2).**

Note. Readers who wish to find approximation conditions directly can refer to the end of each subsection (S4.2, S4.3) in this appendix. This appendix follows the order in which the conditions were identified.

**S4.1. Setting**

We define

$$P_{X^{ep}=e_T,X=e_T,Y(X^{ep})=y_1,Y(X)=y_2|\boldsymbol{z}} := \mathrm{P}(Y(X^{ep} = e_T) = y_1, Y(X = e_T) = y_2|\boldsymbol{z})$$

for $x \in \mathcal{X} \subset \mathbb{R}$. This probability is defined as $Y$ is discrete.

Given $\boldsymbol{z}$, $Y(X = e_T)$ is a single value if the outcome is deterministic but is a distribution if the outcome is non-deterministic. As in the manuscript, let $W$ denote a random variable that indicates a set of non-confounder risk factors that makes the outcome non-deterministic. $Y(X^{ep} = e_T)$ is not a single value but a distribution because this may refer to $Y(X = x_1)$, $Y(X = x_2)$,..., and so on. Furthermore, if the outcome is non-deterministic, each $Y(X = x)$ is not a single value but a distribution. Thus, to identify conditions for minimal bias in estimating *AEE*, these pluralistic characteristics need to be addressed. $P_{X^{ep}=e_T,X=e_T,Y(X^{ep})=y_1,Y(X)=y_2|\boldsymbol{z}}$ is introduced to address these.

We first want to find out when the following equality holds:

$$E[Y(X^{ep} = e_T)|\boldsymbol{z}] = E[Y(X = e_T)|\boldsymbol{z}]$$

for $e_T \in \mathbb{R}$ and $\delta \in \mathbb{R}$. This implies $E[Y(X^{ep} = e_T + \delta)|\boldsymbol{z}] = E[Y(X = e_T + \delta)|\boldsymbol{z}]$. Note that under EUEME, $E[Y(X^{ep})|X^{ep} = x^{ep}, \boldsymbol{z}] = E[Y(X^{ep})|\boldsymbol{z}]$.

Then,

$$AEE(X^{ep}|\boldsymbol{z}) = E[Y(X^{ep} = e_T + \delta)|\boldsymbol{z}] - E[Y(X^{ep} = e_T)|\boldsymbol{z}]$$

$$= E[Y(X = e_T + \delta)|\boldsymbol{z}] - E[Y(X = e_T)|\boldsymbol{z}] = AEE(X|\boldsymbol{z})$$

Note that under (conditional) exchangeability, $E[Y(X)|X = x, \boldsymbol{z}] = E[Y(X)|\boldsymbol{z}]$.

The expectation in $E[Y(X^{ep} = e_T)|\boldsymbol{z}]$ is taken with respect to the (joint) probability distribution regarding $X$ given $X^{ep} = e_T$ and $\boldsymbol{z}$ (and $W$ if the outcome is non-deterministic). The expectation in $E[Y(X = e_T)|\boldsymbol{z}]$ is taken with respect to $W$.

As $V$ is not an effect modifier and not a confounder, $AEE(X|\boldsymbol{z}) = AEE(X|\boldsymbol{z}')$. To find out the conditions for $E[Y(X^{ep} = e_T)|\boldsymbol{z}] = E[Y(X = e_T)|\boldsymbol{z}]$, we rewrite this as

$$E[Y(X^{ep} = e_T)|\boldsymbol{z}] - E[Y(X = e_T)|\boldsymbol{z}] = 0$$

$$E[Y(X^{ep} = e_T)|\boldsymbol{z}] - E[Y(X = e_T)|\boldsymbol{z}] = E[Y(X^{ep} = e_T) - Y(X = e_T)|\boldsymbol{z}]$$

$$= \sum_{y_1 \in \mathcal{Y}, y_2 \in \mathcal{Y}} (y_1 - y_2) P_{X^{ep}=e_T, X=e_T, Y(X^{ep})=y_1, Y(X)=y_2|\boldsymbol{z}} = 0$$

If $P_{X^{ep}=e_T, X=e_T, Y(X^{ep})=y_1, Y(X)=y_2|\boldsymbol{z}}$ is symmetric around $y_1 = y_2$, then $E[Y(X^{ep} = e_T)|\boldsymbol{z}] = E[Y(X = e_T)|\boldsymbol{z}]$ holds. We can also consider approximation, $E[Y(X^{ep} = e_T)|\boldsymbol{z}] - E[Y(X = e_T)|\boldsymbol{z}] \ll \varepsilon$ where $\varepsilon$ is a certain real number.

$P_{X^{ep}=e_T, X=e_T, Y(X^{ep})=y_1, Y(X)=y_2|\boldsymbol{z}}$ can be decomposed as

$$P_{X^{ep}=e_T, X=e_T, Y(X^{ep})=y_1, Y(X)=y_2|\boldsymbol{z}} = \sum_{W} P(Y(X = e_T) = y_2|w, \boldsymbol{z}) \times P(Y(X^{ep} = e_T) = y_1|w, \boldsymbol{z})$$

if W is discrete; and

$P_{X^{ep}=e_T,X=e_T,Y(X^{ep})=y_1,Y(X)=y_2|\boldsymbol{z}}$

$$= \int_{-\infty}^{\infty} P(Y(X = e_T) = y_2|w, \boldsymbol{z}) \times P(Y(X^{ep} = e_T) = y_1|w, \boldsymbol{z}) f_W(w) dw$$

if W is continuous.

### S4.2. When the Exposure Effect is Approximately Strictly Linear

Suppose the true potential outcomes model/outcome data-generating model is

$$Y(X) = \alpha + h(X) + j(\boldsymbol{z}) + W$$

where $h$ *and* $j$ are functions that yield discrete outputs and $W$ is a discrete random variable. Recall that the outcome is discrete. The estimation of the model will be discussed later. With this model,

$$P(Y(X^{ep} = e_T) = y_1|w, \boldsymbol{z}) = \int_{-\infty}^{\infty} \mathbf{1}\{y_1 = \alpha + h(x) + j(\boldsymbol{z}) + w\} f_{X|X^{ep},\boldsymbol{Z}}(x|e_T, \boldsymbol{z}) dx$$

$$P(Y(X = e_T) = y_2|w, \boldsymbol{z}) = \mathbf{1}\{y_2 = \alpha + h(e_T) + j(\boldsymbol{z}) + w\}$$

where $f_{X|X^{ep},\boldsymbol{Z}}(x|e_T, \boldsymbol{z})$ is the probability density function of $X$ given $X^{ep} = e_T, \boldsymbol{Z} = \boldsymbol{z}$. We use $f$ to denote a density function. Accordingly, $P_{X^{ep}=e_T,X=e_T,Y(X^{ep})=y_1,Y(X)=y_2|\boldsymbol{z}} =$

$$\int_{-\infty}^{\infty} \mathbf{1}\{y_2 = \alpha + h(e_T) + j(\boldsymbol{z}) + w\} \mathbf{1}\{y_1 = \alpha + h(x) + j(\boldsymbol{z}) + w\} f_{X|X^{ep},\boldsymbol{Z}}(x|e_T, \boldsymbol{z}) dx$$

$f_{X|X^{ep},\boldsymbol{Z}}(x|e_T, \boldsymbol{z})$ depends on the form of EME. For NB error, where $m(X^{ep}) = g\big(m(X)\big) + V + U$,

$$f_{X|X^{ep},\boldsymbol{Z}}(x|e_T, \boldsymbol{z}) = \frac{f_{X|\boldsymbol{Z}}(x|\boldsymbol{z}) f_U\big(u = m(e_T) - g\big(m(x)\big) - v\big)}{f_{X^{ep}|\boldsymbol{Z}}(e_T|\boldsymbol{z})}$$

by the convolution of probability density functions. For Berkson error, where $m(X) = g\big(m(X^{ep})\big) + V + U$,

$$f_{X|X^{ep},\mathbf{Z}}(x|e_T,\mathbf{z}) = f_U\big(u = m(x) - g\big(m(e_T)\big) - v\big)\left|\frac{d}{dx}m(x)\right|$$

If EME is purely Berksonian (i.e., $X = X^{ep} + U$), then

$$f_{X|X^{ep},\mathbf{Z}}(x|e_T,\mathbf{z}) = f_U(u = x - e_T)$$

which is symmetric around $x = e_T$.

Suppose $h(X) \coloneqq round(\beta X, d)$. For simplicity, we use $[\beta X]$ to denote this. Then,

$$P_{X^{ep}=e_T,X=e_T,Y(X^{ep})=y_1,Y(X)=y_2|\mathbf{z}} =$$

$$\int_{-\infty}^{\infty}\mathbf{1}\{y_2 = \alpha + [\beta_X e_T] + j(\mathbf{z}) + w\}\mathbf{1}\{y_1 = \alpha + [\beta_X x] + j(\mathbf{z}) + w\}f_U(u = x - e_T)dx$$

which would be symmetric around $y_1 = y_2$ if the rounding function satisfies

$$[\beta_X(e_T + u)] - [\beta_X e_T] = -([\beta_X(e_T - u)] - [\beta_X e_T])$$

This condition would approximately hold if $h(X) \approx \beta X$. Thus, the following equalities would approximately hold:

$$E[Y(X^{ep} = e_T)|\mathbf{z}'] \approx E[Y(X = e_T)|\mathbf{z}'] \text{ (Approximate Equality 1-1)}$$

and as EME is purely Berksonian, instead of $\mathbf{z}'$, we can use

$$E[Y(X^{ep} = e_T)|\mathbf{z}] \approx E[Y(X = e_T)|\mathbf{z}] \text{ (Approximate Equality 1-2)}$$

We can also use Approximate Equality 1-2 when $X = X^{ep} + U$ but $X^{ep}$ is correlated with $V$ that is correlated with $Y$ like when $V$ becomes a confounder when using $X^{ep}$ (Figure S1H).

Also, when $V$ is not an effect modifier (and not a confounder when using $X$),

$$E[Y(X^{ep} = e_T)|\mathbf{z}] = E[Y(X^{ep} = e_T)|\mathbf{z}'] \approx E[Y(X = e_T)|\mathbf{z}'] = E[Y(X = e_T)|\mathbf{z}]$$

(Approximate Equality 2)

also holds.

If EME is Berksonian but follows this form, $X = X^{ep} + V + U$, in which $V$ makes EME differential,

$$f_{X|X^{ep},V,\boldsymbol{Z}}(x|e_T, v, \boldsymbol{z}) = f_U(u = x - e_T - v)$$

This would be symmetric around $x = e_T + v$. So, $P_{X^{ep}=e_T,X=e_T,Y(X^{ep})=y_1,Y(X)=y_2|\boldsymbol{z}}$ would be symmetric around $y_1 = y_2 - [\beta_X v]$. Thus,

$$E[Y(X^{ep} = e_T)|\boldsymbol{z}] \approx E[Y(X = e_T)|\boldsymbol{z}] + \beta_X v$$

$$E[Y(X^{ep} = e_T + \delta)|\boldsymbol{z}] \approx E[Y(X = e_T + \delta)|\boldsymbol{z}] + \beta_X v$$

Approximate Equalities 1-2 and 2 do not hold. However, for conditional AEE that is

$$E[Y(X^{ep} = e_T + \delta)|\boldsymbol{z}] - E[Y(X^{ep} = e_T)|\boldsymbol{z}]$$

$$\approx E[Y(X = e_T + \delta)|\boldsymbol{z}] + \beta_X v - (E[Y(X = e_T)|\boldsymbol{z}] + \beta_X v)$$

$$= E[Y(X = e_T + \delta)|\boldsymbol{z}] - E[Y(X = e_T)|\boldsymbol{z}] = AEE(X|\boldsymbol{z})$$

Thus, the use of $X^{ep}$ would not introduce bias or may introduce negligible bias.

Now, we consider estimation. In practice, investigators need to estimate the outcome model. $\widehat{W_1}$, $\hat{\alpha}$, $\widehat{\beta_X}$, and $\hat{\jmath}$ are used to estimate $AEE(X^{ep}|\boldsymbol{z})$. Being similar to Approximate Equalities 1-1, 1-2, and 2, the following needs to be satisfied:

$$E\left[E\left[\hat{Y}(X^{ep} = e_T)\middle|\boldsymbol{z}'\right]\right] \approx E[Y(X = e_T)|\boldsymbol{z}'] \text{ (Approximate Equality 3) or}$$

$$E\left[E\left[\hat{Y}(X^{ep} = e_T)\middle|\boldsymbol{z}\right]\right] \approx E[Y(X = e_T)|\boldsymbol{z}] \text{ (Approximate Equality 4)}$$

For simplicity, with the notation $\boldsymbol{z}$, we should show

$$E\left[E\left[\sum\nolimits_{y_1\in\mathcal{Y},y_2\in\mathcal{Y}}(y_1 - y_2)\hat{P}_{X^{ep}=e_T,X=e_T,Y(X^{ep})=y_1,Y(X)=y_2|\boldsymbol{z}}\right]\right] \cong 0$$

The outer expectation is taken with respect to $\hat{P}_{X^{ep}=e_T,X=e_T,Y(X^{ep})=y_1,Y(X)=y_2|\boldsymbol{z}}$.

For purely Berkson error,

$$\hat{P}_{X^{ep}=e_T, X=e_T, Y(X^{ep})=y_1, Y(X)=y_2|\mathbf{z}}$$

$$= \int_{-\infty}^{\infty} \mathbf{1}\{y_2 = \alpha + [\beta_X e_T] + j(\mathbf{z}) + w\}\mathbf{1}\{y_1$$

$$= \hat{\alpha} + [\hat{\beta}_X x] + \hat{j}(\mathbf{z}) + w\} f_U(u = x - e_T) dx$$

When EME is purely Berksonian and the effect of $X$ is strictly linear, $E\left[\widehat{\beta_X}\right] = \beta_X$ as known in the exposure measurement error literature[14-16]. For $h(x) \coloneqq [\beta_X X]$, $E\left[\widehat{\beta_X}\right] \approx \beta_X$ when $d$ is such that $h(X) \approx \beta X$. If $E[\hat{\alpha}]$ and $E[\hat{j}]$ equal $\alpha$ and $j$, respectively, $\hat{P}_{X^{ep}=e_T, X=e_T, Y(X^{ep})=y_1, Y(X)=y_2|\mathbf{z}}$ would be symmetric around $y_1 = y_2$ and Approximate Equalities 3 and 4 would hold. Otherwise, the symmetry would be shifted by either the difference between $E[\hat{\alpha}]$ and $\alpha$, the difference between $E[\hat{j}]$ and $j$, or both. We refer to either the difference–or both collectively, as *diff.* So,

$$E\left[E\left[\hat{Y}(X^{ep} = e_T)\middle|\mathbf{z}'\right]\right] \approx E[Y(X = e_T)|\mathbf{z}] + diff$$

$$E\left[E\left[\hat{Y}(X^{ep} = e_T + \delta)\middle|\mathbf{z}'\right]\right] \approx E[Y(X = e_T + \delta)|\mathbf{z}] + diff$$

Then,

$$E\left[\widehat{AEE}(X^{ep}|\mathbf{z})\right] = E\left[E\left[\hat{Y}(X^{ep} = e_T + \delta)\middle|\mathbf{z}'\right]\right] - E\left[E\left[\hat{Y}(X^{ep} = e_T)\middle|\mathbf{z}'\right]\right]$$

$$\approx E[Y(X = e_T)|\mathbf{z}] + diff - (E[Y(X = e_T + \delta)|\mathbf{z}] + diff)$$

$$= E[Y(X = e_T + \delta)|\mathbf{z}] - E[Y(X = e_T)|\mathbf{z}] = AEE(X|\mathbf{z})$$

For $X = X^{ep} + V + U$, similarly,

$$E\left[E\left[\hat{Y}(X^{ep} = e_T)\middle|\mathbf{z}\right]\right] \approx E[Y(X = e_T)|\mathbf{z}] + \beta_X v, \text{ or}$$

$$E\left[E\left[\hat{Y}(X^{ep} = e_T)\middle|\mathbf{z}\right]\right] \approx E[Y(X = e_T)|\mathbf{z}] + \beta_X v + diff$$

$$E\left[E\left[\hat{Y}(X^{ep} = e_T + \delta)\middle|\boldsymbol{z}\right]\right] \approx E[Y(X = e_T + \delta)|\boldsymbol{z}] + \beta_X v, \text{ or}$$

$$E\left[E\left[\hat{Y}(X^{ep} = e_T + \delta)\middle|\boldsymbol{z}\right]\right] \approx E[Y(X = e_T + \delta)|\boldsymbol{z}] + \beta_X v + diff$$

So, Approximate Equalities 3 and 4 do not hold. Nevertheless,

$$E\left[\widehat{AEE}(X^{ep}|\boldsymbol{z})\right] = E\left[E\left[\hat{Y}(X^{ep} = e_T + \delta)\middle|\boldsymbol{z}'\right]\right] - E\left[E\left[\hat{Y}(X^{ep} = e_T)\middle|\boldsymbol{z}'\right]\right]$$

$$\approx E[Y(X = e_T + \delta)|\boldsymbol{z}] - E[Y(X = e_T)|\boldsymbol{z}] = AEE(X|\boldsymbol{z})$$

when $V$ is not an effect modifier (and not a confounder when using $X$), $AEE(X|\boldsymbol{z}) = AEE(X|\boldsymbol{z}')$

This completes the proof of Theorem 1.

### S4.3. When the Exposure Effect is Not Strictly Linear

Different forms of the outcome model will be examined in this sub-section.

$P_{X^{ep}=e_T, X=e_T, Y(X^{ep})=y_1, Y(X)=y_2|\boldsymbol{z}}$ may depend on the form of the true potential outcomes model.

First, when the potential outcomes model is

$$q\big(Y(X)\big) = \alpha + h(X) + j(\boldsymbol{z}) + W$$

and EME is purely Berksonian,

$$P_{X^{ep}=e_T, X=e_T, Y(X^{ep})=y_1, Y(X)=y_2|\boldsymbol{z}}$$

$$= \int_{-\infty}^{\infty} \mathbf{1}\{q(y_2) = \alpha + h(e_T) + j(\boldsymbol{z}) + w\}\mathbf{1}\{q(y_1) = \alpha + h(x) + j(\boldsymbol{z}) + w\} f_U(u$$

$$= x - e_T)dx$$

If $q$ is a monotonic function such as the logarithm, $q(y_2)$ and $q(y_1)$ could be seen as $y_2^*$ and $y_1^*$, respectively. If $h(x) \approx \beta x$ (may not need to be discrete depending on $q$), then, the consequence would be similar to the case for approximately strictly linear exposure effect. We will discuss the logistic link later. If $h(x)$ is non-linear, the following approximation may hold under the following conditions, which may not be exhaustive,

- The error is purely Berksonian, AND
- $h(x)$ is approximately linear in the range of $X$, $[e_T - u_P, e_T + u_P]$, where $\pm u_P$ indicates the primary mass of $f_U(u)$ (e.g., $h(x)$ may approximate the first order of a Taylor expansion)

The second condition is because $h(x) \approx \beta x$ in that range. This does not mean that $h(x)$ should be approximately linear within all possible values of $X$. For example, if $f_U(u)$ has the most of their mass near zero, which is the case when EME is small, $P_{X^{ep}=e_T, X=e_T, Y_1=y_1, Y_2=y_2|\mathbf{z}}$ may be approximately symmetric around $y_1^* - y_2^*$. If $h(x)$ is not strongly non-linear, EME may not need to be small for the approximation. Also, if the effect size is small, because it is more likely that $h(e_T + u) \approx h(e_T)$, the error may not need to be small for the approximation (e.g., modest error).

As in Appendix S4.2, these can be extended for $X = X^{ep} + V + U$. In this case,

- $h(x)$ is approximately linear in the range of $X$, $[e_T + v - u_P, e_T + v + u_P]$, where $\pm u_P$ indicates the primary mass of $f_U(u)$ (e.g., $h(x)$ may approximate the first order of a Taylor expansion)

Second, when the potential outcomes model/data-generating model takes the following logistic model form for binary outcomes,

$$P(Y(X = x) = 1|\mathbf{z}, w) = \frac{1}{1 + \exp(-(h(x) + j(\mathbf{z}) + w))}$$

where $w$, $h(x)$, and $j(\mathbf{z})$ can be continuous.

$$P(Y(X = e_T) = 1|w, \boldsymbol{z}) = \frac{1}{1 + \exp\big(-(h(e_T) + j(\boldsymbol{z}) + w)\big)}$$

$$P(Y(X^{ep} = e_T) = 1|w, \boldsymbol{z}) = \int_{-\infty}^{\infty} \frac{1}{1 + \exp\big(-(h(x) + j(\boldsymbol{z}) + w)\big)} f_{X|X^{ep},\boldsymbol{Z}}(x|e_T, \boldsymbol{z})dx$$

For simplicity, let $W$ be continuous. For binary outcomes, we want to find out $P_{X^{ep}=e_T,X=e_T,Y(X^{ep})=1,Y(X)=0|\boldsymbol{z}}$ and $P_{X^{ep}=e_T,X=e_T,Y(X^{ep})=0,Y(X)=1|\boldsymbol{z}}$ such that

$$E[Y(X = e_T)|\boldsymbol{z}] - E[Y(X^{ep} = e_T)|\boldsymbol{z}] = E[Y(X = e_T) - Y(X^{ep} = e_T)|\boldsymbol{z}]$$

$$= \int_{-\infty}^{\infty} (1 - 0) \times \left(P_{X^{ep}=e_T,X=e_T,Y(X^{ep})=1,Y(X)=0|\boldsymbol{z}}\right) + (0$$

$$- 1) \times \left(P_{X^{ep}=e_T,X=e_T,Y(X^{ep})=0,Y(X)=1|\boldsymbol{z}}\right) f_W(w)dw \cong 0$$

Let $D \coloneqq P_{X^{ep}=e_T,X=e_T,Y(X^{ep})=1,Y(X)=0|\boldsymbol{z}} - P_{X^{ep}=e_T,X=e_T,Y(X^{ep})=0,Y(X)=1|\boldsymbol{z}}$. When the error follows the form $X = X^{ep} + V + U$,

$P_{X^{ep}=e_T,X=e_T,Y(X^{ep})=0,Y(X)=1|\boldsymbol{z}} =$

$$= \frac{1}{1 + \exp\big(-(\alpha + h(e_T) + j(\boldsymbol{z}) + w)\big)} \int_{-\infty}^{\infty} \frac{\exp\big(-(\alpha + h(x) + j(\boldsymbol{z}) + w)\big)}{1 + \exp\big(-(\alpha + h(x) + j(\boldsymbol{z}) + w)\big)} f_U(u = x - e_T - v)dx =$$

$$= \frac{1}{1 + \exp\big(-(\alpha + h(e_T) + j(\boldsymbol{z}) + w)\big)} \int_{-\infty}^{\infty} \frac{\exp\big(-(\alpha + h(e_T + v + u) + j(\boldsymbol{z}) + w)\big)}{1 + \exp\big(-(\alpha + h(e_T + v + u) + j(\boldsymbol{z}) + w)\big)} f_U(u)du$$

$P_{X^{ep}=e_T,X=e_T,Y(X^{ep})=1,Y(X)=0|\boldsymbol{z}} =$

$$= \frac{\exp\big(-(\alpha + h(e_T) + j(\boldsymbol{z}) + w)\big)}{1 + \exp\big(-(\alpha + h(e_T) + j(\boldsymbol{z}) + w)\big)} \int_{-\infty}^{\infty} \frac{1}{1 + \exp\big(-(\alpha + h(x) + j(\boldsymbol{z}) + w)\big)} f_U(u = x - e_T - v)dx =$$

$$= \frac{\exp\big(-(\alpha + h(e_T) + j(\boldsymbol{z}) + w)\big)}{1 + \exp\big(-(\alpha + h(e_T) + j(\boldsymbol{z}) + w)\big)} \int_{-\infty}^{\infty} \frac{1}{1 + \exp\big(-(\alpha + h(e_T + v + u) + j(\boldsymbol{z}) + w)\big)} f_U(u)du$$

When the error is purely Berksonian, then $v = 0$.

*D* can be approximately zero (then $E[Y(X=e_T)|\boldsymbol{z}]-E[Y(X^{ep}=e_T)|\boldsymbol{z}]\approx 0$) when $\exp\left(-(\alpha+h(e_T)+j(\boldsymbol{z})+w)\right)\approx\exp\left(-(\alpha+h(e_T+v+u)+j(\boldsymbol{z})+w)\right)\approx k_1$, implying

$$AEE(X^{ep}|\boldsymbol{z})\approx AEE(X|\boldsymbol{z})$$

or

$$AEE(X|\boldsymbol{z})=AEE(X|\boldsymbol{z}')$$

when $V$ is not an effect modifier (and not a confounder when using $X$). To demonstrate,

$P_{X^{ep}=e_T,X=e_T,Y(X^{ep})=0,Y(X)=1|\boldsymbol{z}}$

$$\approx\frac{1}{1+k_1}\int_{-\infty}^{\infty}\frac{k_1}{1+k_1}f_U(u)du=\frac{k_1}{(1+k_1)^2}\int_{-\infty}^{\infty}f_U(u)\,du$$

$P_{X^{ep}=e_T,X=e_T,Y(X^{ep})=1,Y(X)=0|\boldsymbol{z}}$

$$\approx\frac{k_1}{1+k_1}\int_{-\infty}^{\infty}\frac{1}{1+k_1}f_U(u)du=\frac{k_1}{(1+k_1)^2}\int_{-\infty}^{\infty}f_U(u)\,du$$

Thus, $D\approx 0$ when $\exp\left(-(\alpha+h(e_T)+j(\boldsymbol{z})+w_2)\right)\approx\exp\left(-(\alpha+h(e_T+v+u)+j(\boldsymbol{z})+w_1)\right)\approx k_1$. For this approximation, $v$ should be zero, meaning that EME is purely Berskonian and $u$ should be modest, meaning that the error should be modest. If the exposure effect is small, the error may not need to be small because it is likely $h(e_T+u)\approx h(e_T)$.

There may exist other approximation conditions. Suppose $D\approx\pm\,k$, even if $k\neq 0$,

$$E[Y(X=e_T)-Y(X^{ep}=e_T)|\boldsymbol{z}]\approx\pm k$$

$$E[Y(X=e_T+\delta)-Y(X^{ep}=e_T+\delta)|\boldsymbol{z}]\approx\pm k$$

Subtracting the former from the latter implies

$$AEE(X^{ep}|\boldsymbol{z})\approx AEE(X|\boldsymbol{z})$$

If the probability of $Y$=1 is small (i.e., rare disease assumption), this approximation may likely occur. To demonstrate,

let $A_0 = \exp\big(-(\alpha + h(e_T) + j(\boldsymbol{z}) + w)\big)$ and $B_0 = \exp\big(-(\alpha + h(e_T + v + u) + j(\boldsymbol{z}) + w)\big)$,

let $A_1 = \exp\big(-(\alpha + h(e_T + \delta) + j(\boldsymbol{z}) + w)\big)$ and $B_1 = \exp\big(-(\alpha + h(e_T + \delta + v + u) + j(\boldsymbol{z}) + w)\big)$,

$P_{X^{ep}=e_T, X=e_T, Y(X^{ep})=0, Y(X)=1|\boldsymbol{z}} =$

$$= \frac{1}{1 + A_0} \int_{-\infty}^{\infty} \frac{B_0}{1 + B_0} f_U(u) du$$

$P_{X^{ep}=e_T, X=e_T, Y(X^{ep})=1, Y(X)=0|\boldsymbol{z}} =$

$$= \frac{A_0}{1 + A_0} \int_{-\infty}^{\infty} \frac{1}{1 + B_0} f_U(u) du$$

$P_{X^{ep}=e_T+\delta, X=e_T+\delta, Y(X^{ep})=0, Y(X)=1|\boldsymbol{z}} =$

$$= \frac{1}{1 + A_1} \int_{-\infty}^{\infty} \frac{B_1}{1 + B_1} f_U(u) du$$

$P_{X^{ep}=e_T+\delta, X=e_T+\delta, Y(X^{ep})=1, Y(X)=0|\boldsymbol{z}} =$

$$= \frac{A_1}{1 + A_1} \int_{-\infty}^{\infty} \frac{1}{1 + B_1} f_U(u) du$$

Suppose the probability of $Y$=1 is small such that $\frac{1}{1+A_0} \approx \frac{1}{1+A_1} \approx k_2 \ll 0.1$ and $\frac{A_0}{1+A_0} \approx \frac{A_1}{1+A_1} \approx 1 - k_2 \gg 0.9$ and $\frac{1}{1+B_0} \approx \frac{1}{1+B_1} \approx k_3 \ll 0.1$ and $\frac{B_0}{1+B_0} \approx \frac{B_1}{1+B_1} \approx 1 - k_3 \gg 0.9$. Then

$E[Y(X = e_T) - Y(X^{ep} = e_T)|\boldsymbol{z}]$

$$= P_{X^{ep}=e_T, X=e_T, Y(X^{ep})=1, Y(X)=0|\boldsymbol{z}} - P_{X^{ep}=e_T, X=e_T, Y(X^{ep})=0, Y(X)=1|\boldsymbol{z}}$$

$$= \frac{A_0}{1+A_0}\int_{-\infty}^{\infty}\frac{1}{1+B_0}f_U(u)\,du - \frac{1}{1+A_0}\int_{-\infty}^{\infty}\frac{B_0}{1+B_0}f_U(u)\,du$$

$$\approx (1-k_2)\int_{-\infty}^{\infty}\frac{1}{1+B_0}f_U(u)\,du - k_2\int_{-\infty}^{\infty}\frac{B_0}{1+B_0}f_U(u)\,du$$

$$\approx (1-k_2)k_3 - k_2(1-k_3)$$

$E[Y(X = e_T + \delta) - Y(X^{ep} = e_T + \delta)|\boldsymbol{z}]$

$$\approx (1-k_2)\int_{-\infty}^{\infty}\frac{1}{1+B_1}f_U(u)\,du - k_2\int_{-\infty}^{\infty}\frac{B_1}{1+B_1}f_U(u)\,du$$

$$\approx (1-k_2)k_3 - k_2(1-k_3)$$

Thus, $E[Y(X = e_T) - Y(X^{ep} = e_T)|\boldsymbol{z}] - E[Y(X = e_T + \delta) - Y(X^{ep} = e_T + \delta)|\boldsymbol{z}] \approx 0$.

Now, we consider estimation. In practice, $\hat{\alpha}$, $\hat{h}$, and $\hat{j}$ will be used in place of $\alpha$, $h$, and $j$ in $P(Y(X^{ep} = e_T) = 1|w, \boldsymbol{z})$ (this is regarding $B_0$ or $B_1$ above) because investigators need to estimate the model. $\hat{h}$ may not necessarily be unbiased when the error is purely Berksonian when the outcome model is non-linear. As known in the exposure measurement error literature[14-16], there exist conditions for minimal or negligible bias. When $E[\hat{j}]$ is significantly different from $j$, then $AEE(X|\boldsymbol{z} = 0)$ may be targeted instead, in which $\boldsymbol{z}$ may be defined as a set of centered variables. We summarize possible approximation conditions for $E[\widehat{AEE}(X^{ep}|\boldsymbol{z})] \approx AEE(X|\boldsymbol{z})$.

1. The error is purely Berksonian or $X = X^{ep} + V + U$ (when conditional on $V$), AND
2. $q(y)$ can be linearly approximated within a certain range of $y$, OR

    (e.g., log-linear models when Poisson count ranges are ≥100)

3. $h(x)$ can be linearly approximated in the range of $X$, $[e_T + v - u_P, e_T + v + u_P]$, OR
4. For logistic functions and binary outcomes, the probability of *Y*=1 ranges from 0.2 to 0.8 (because the logistic function is almost linear in this range) or the probability of *Y*=1 is low (e.g., $\ll 0.1$), OR
5. The exposure effect is small (compared to the entire range of the probability of *Y=1* or $q(y)$), OR
6. The error is modest.

This completes the proof of Theorem 2. Due to the nonlinear nature of the problem, there may exist perturbations. We may not rule out the possibility that other conditions may exist.

**Appendix S5. Proof of Theorem 3.**

Under Condition 2, $U_1^*,\ldots,U_P^*$, $U_1^{C^*},\ldots,$ $U_Q^{C^*}$, and $U^{V^*}$ are independent of $X_{p,MRC}^{ep}$ because $X_{p,MRC}^{ep}$ is a linear combination of all the variables in all MRC models (including the original covariates and potentially transformed covariates as part of $s$); each of this variable is independent of the residuals. The residuals can be seen as unmeasured variables in the main epidemiological study; they do not act as confounders because they are independent of $X_{p,MRC}^{ep}$. Also, the error model for the exposure of interest is $X_p = X_{p,MRC}^{ep} + U_p^*$ where $E\left[U_p^*\right] = 0$ and $U_p^*$ is independent of $X_{1,MRC}^{ep},\ldots,$ $X_{P,MRC}^{ep}$, $C_{1,MRC}^{ep},\ldots,$ $C_{Q,MRC}^{ep}$, $V_{MRC}^{ep}$, and $\boldsymbol{Z}'^{-\boldsymbol{C}}$, which is Berksonian. Adjusting for the other MR-Calibrated variables suffice for adjustment for $X_1,\ldots,$ $X_{p-1}$, $X_{p+1}$, $X_P$, $C_1,\ldots,$ $C_Q$, $V$ because the other MR-Calibrated variables correspond to the portions of these true variables after subtracting the residuals. Consequently, conditioning on the other MR-Calibrated variables and $\boldsymbol{Z}'^{-\boldsymbol{C}}$ is sufficient for EUEME under Conditions 1 and 2.

**Remark 1 (Independence Depending on $s$).** For example, if $s\left(X_p^{ep}\right) = \eta {X_p^{ep}}^2$, the independence between ${X_p^{ep}}^2$ and $U_p^*$ that would be typical in standard regression may not imply the independence between $X_p^{ep}$ and $U_p^*$. The independence between $X_p^{ep}$ and $U_p^*$ may be preserved by making $s\left(X_p^{ep}\right) = \eta_1 X_p^{ep} + \eta_2 {X_p^{ep}}^2$ or ensuring that $s\left(X_p^{ep}\right)$ lies within a domain that the monotonicity of $s\left(X_p^{ep}\right)$ is preserved. Many spline and smoothing techniques include the original variable itself as one of the basis functions.

**Remark 2 (For multiplicative error)**. The transformation of the dependent variable in MRC models may not be considered. Such transformation is proven to be able to reduce the error correction ability in the context of regression calibration[16]. When error is multiplicative, additive MRC models with non-linear functions and interaction terms could be considered. If the residuals may be heteroskedastic, weighted regression or quasi-likelihood/variance model estimation could be used[16,17].

**Remark 3 (Conditional AEE vs. Marginal AEE).** Investigators may wish to estimate $AEE(X|\boldsymbol{z}')$ using $E\left[\widehat{AEE}\left(X_{p,MRC}^{ep}\middle|\boldsymbol{z}^{**}\right)\right]$ or $AEE(X)$ using $E_{\boldsymbol{Z}^{**}}\left[E\left[\widehat{AEE}\left(X_{p,MRC}^{ep}\middle|\boldsymbol{z}^{**}\right)\right]\right]$. When it is unlikely that the effect of exposures and confounders are (approximately) strictly linear and the error is small, we suggest estimating $AEE(X|\boldsymbol{z}')$ using $E\left[\widehat{AEE}\left(X_{p,MRC}^{ep}\middle|\boldsymbol{Z}^{**}=0\right)\right]$ where $\boldsymbol{Z}^{**}$ may be defined as centered variables. For $AEE(X)$, bias should be minimal for the regression coefficient estimates of all variables in a Q-model, not just $X_{p,MRC}^{ep}$.

**Appendix S6. Simulation Analyses.**

Simulation analyses were conducted to identify the bias-correction capability of the developed framework when the exposure effect is not strictly linear. Several simulation scenarios were considered:

- Error magnitude: large error ($R^2 \approx 0.5$ for $X_1/X_2/X_3/C/\ V$ and $X_1^{ep}/X_2^{ep}/X_3^{ep}/C^{ep}/\ V^{ep}$), moderate error ($R^2 \approx 0.75$), low error ($R^2 \approx 0.90$)
- Error type: multiplicative or additive errors
- Potential outcomes models: logistic or log-linear (Poisson).
- $h(x)$: linear terms for $X_1$ and $X_2$; linear and quadratic terms for $X_3$, linear and cubic terms for $C$.
- Baseline risk, which affects the extent of linear approximation of the outcome model: low, medium, and high
- Correlation structures: Pearson correlations among $X_1$, $X_2$, $X_3$, $C$ are

$X_1 - X_2$: ~0.67;

$X_1 - X_3$: ~0.54;

$X_1 - C$: ~0.57;

$X_2 - X_3$: ~0.68;

$X_2 - C$: ~0.68;

$X_3 - C$: ~0.56

Even for high error scenario, correlations among $X_1^{ep}$, $X_2^{ep}$, $X_3^{ep}$, $C^{ep}$ remained at least moderate (~0.3). See details below.

### S6.1. Simulation Samples Generation

$X_1$ was generated using

$$Gamma(7,5) \times 5 + 0.5C + 1.1A + 0.7B$$

$X_2$ was generated using

$$Gamma(2,4) \times 3 + 0.4C + 0.9A + 0.7B$$

$X_3$ was generated using

$$Gamma(4,3) \times 4 + 0.5C + 1.2A + 0.7B$$

$A$ was generated using

$$Uniform(0,5)$$

$B$ was generated using

$$Discrete\ Uniform(0,4)$$

$C$ was generated using

$$Gamma(5,5) \times 10$$

$X_1$, $X_2$, and $X_3$ mimick empirical distributions of air pollutants in different cities in the world. $C$ mimicks an empirical distribution of daily temperature in a moderately warm climate region. In this simulation, $0.5C$ and $0.4C$ were used to generate $X_1$, $X_2$, and $X_3$, as some air pollutants (e.g., some constituents in $PM_{2.5}$) may be correlated with temperature. High temperature increases $O_3$ by facilitating the photochemical formation of $O_3$ in the air; high temperatures may increase $PM_{2.5}$ due to various processes (e.g., photochemical formations of secondary particles, increased $PM_{2.5}$ emissions that is correlated with increased temperature such as energy demand increases). We acknowledge that the correlation between an air pollutant and temperature can be non-linear (e.g., $O_3$ formation), but the focus of this simulation is not to perfectly mimic real-world circumstances but to evaluate the bias-correction capabilities of the developed modeling framework when multiple covariates are error-prone. $B$ mimicks categorical variables such as socioeconomic status. $A$ mimicks a set of unknown or unmeasured continuous factors that may be related to air pollution levels (e.g., emission activities at areas where individuals reside, atmospheric process). For example, distributions of air pollution may differ by socioeconomic status because socioeconomically deprived individuals may reside in areas with high pollution emission sources.

$V$ mimicking a set of factors that make EME differential, which is not a confounder in this simulation, but a risk factor of $Y$, was generated using

$$V \sim Gamma(2,2)$$

For environmental modelled estimates, uncertainties may vary across individuals if environmental exposure prediction models exhibit different uncertainties across locations or regions or over time, which is not uncommon, and individuals live in different locations and regions. Also, the risk of health outcomes may vary across locations or regions or over time due to contextual factors or other individual-level non-confounding risk factors, as demonstrated by geographically varying health risks across the United States. $X_1^{ep}$, $X_2^{ep}$, $X_3^{ep}$, $C^{ep}$, $V^{ep}$ were generated using the following error models:

1. Multiplicative differential linear error scenarios

$$\log\left(X_1^{ep}\right) = 1.13 \times \log(X_1) + 0.032 \times V + U_1$$

$$\log\left(X_2^{ep}\right) = 0.85 \times \log(X_2) + 0.045 \times V + U_2$$

$$\log(X_3^{ep}) = 0.90 \times \log(X_3) + 0.061 \times V + U_3$$

$$\log(C^{ep}) = 0.89 \times \log(C) + U_4$$

$$V^{ep} = V + U_5 \text{ (additive)}$$

$$U_1 \sim N(0, \tau_1^2)$$

$$U_2 \sim N(0, \tau_2^2)$$

$$U_3 \sim N(0, \tau_3^2)$$

$$U_4 \sim N(0, \tau_4^2)$$

$$U_5 \sim N(0, \tau_5^2)$$

- High error scenario ($R^2 \approx 0.50$ for all $X_1^{ep} - X_1$, $X_2^{ep} - X_2$ , $X_3^{ep} - X_3$, $C^{ep} - C$, $V^{ep} - V$ pairs): $\tau_1^2 = 0.25$; $\tau_2^2 = 0.23$; $\tau_3^2 = 0.23$; $\tau_4^2 = 0.36$; $\tau_5^2 = 0.71$.
- Moderate error scenario ($R^2 \approx 0.75$ for all $X_1^{ep} - X_1$, $X_2^{ep} - X_2$ , $X_3^{ep} - X_3$, $C^{ep} - C$, $V^{ep} - V$ pairs): $\tau_1^2 = 0.145$; $\tau_2^2 = 0.132$; $\tau_3^2 = 0.130$; $\tau_4^2 = 0.21$; $\tau_5^2 = 0.41$
- Small error scenario ($R^2 \approx 0.90$ for all $X_1^{ep} - X_1$, $X_2^{ep} - X_2$ , $X_3^{ep} - X_3$, $C^{ep} - C$, $V^{ep} - V$ pairs): $\tau_1^2 = 0.083$; $\tau_2^2 = 0.071$; $\tau_3^2 = 0.067$; $\tau_4^2 = 0.12$; $\tau_5^2 = 0.24$

Note. Even for high error scenario, correlations among $X_1^{ep}$, $X_2^{ep}$, $X_3^{ep}$, $C^{ep}$ remained moderate (~0.3)

2. Additive differential linear error scenarios

$$X_1^{ep} = 1.13 \times X_1 + \phi_1 \times V + U_1$$

$$X_2^{ep} = 0.85 \times X_2 + \phi_2 \times V + U_2$$

$$X_3^{ep} = 0.90 \times X_3 + \phi_3 \times V + U_3$$

$$C^{ep} = 0.89 \times \log(C) + U_4 \text{ (multiplicative)}$$

$$V^{ep} = V + U_5 \text{ (additive)}$$

$$U_1 \sim N(0, \tau_1^2)$$

$$U_2 \sim N(0, \tau_2^2)$$

$$U_3 \sim N(0, \tau_3^2)$$

$$U_4 \sim N(0, \tau_4^2)$$

$$U_5 \sim N(0, \tau_5^2)$$

- High error scenario ($R^2 \approx 0.50$ for all $X_1^{ep} - X_1$, $X_2^{ep} - X_2$ , $X_3^{ep} - X_3$, $C^{ep} - C$ pairs): $\phi_1 = 8$, $\phi_2 = 4$, $\phi_3 = 6.5$; $\tau_1^2 = 7$; $\tau_2^2 = 3.4$; $\tau_3^2 = 5.5$; $\tau_4^2 = 0.36$; $\tau_5^2 = 0.71$
- Moderate error scenario ($R^2 \approx 0.75$ for all $X_1^{ep} - X_1$, $X_2^{ep} - X_2$ , $X_3^{ep} - X_3$, $C^{ep} - C$ pairs): $\phi_1 = 4$, $\phi_2 = 2$, $\phi_3 = 3.4$; $\tau_1^2 = 4.2$; $\tau_2^2 = 2.1$; $\tau_3^2 = 3.4$; $\tau_4^2 = 0.21$; $\tau_5^2 = 0.41$
- Small error scenario ($R^2 \approx 0.90$ for all $X_1^{ep} - X_1$, $X_2^{ep} - X_2$ , $X_3^{ep} - X_3$, $C^{ep} - C$ pairs): $\phi_1 = 3$, $\phi_2 = 1.4$, $\phi_3 = 2.3$; $\tau_1^2 = 2$; $\tau_2^2 = 1$; $\tau_3^2 = 1.7$; $\tau_4^2 = 0.12$; $\tau_5^2 = 0.24$

Note. Even for high error scenario, correlations among $X_1^{ep}$, $X_2^{ep}$, $X_3^{ep}$, $C^{ep}$ remained moderate to high (0.3–0.6).

We used two types of the outcome generating models: logistic (binary outcome), and log-linear (Poisson count). Figure S2 presents dose-response curves for the effect of $X_1$, $X_2$, $X_3$, and $C$ on the outcome used for simulations. For $X_1$ and $X_2$, linear terms were used to mimick the effects of some air pollutants; for $X_3$, quadratic terms were used to mimick a sub-linear form, that is sometimes observed in air pollution epidemiology, and for $C$, cubic terms were used to mimick a *U* or *J*-shaped form that is often observed in temperature epidemiology. The following items present functional forms to generate outcomes in simulations.

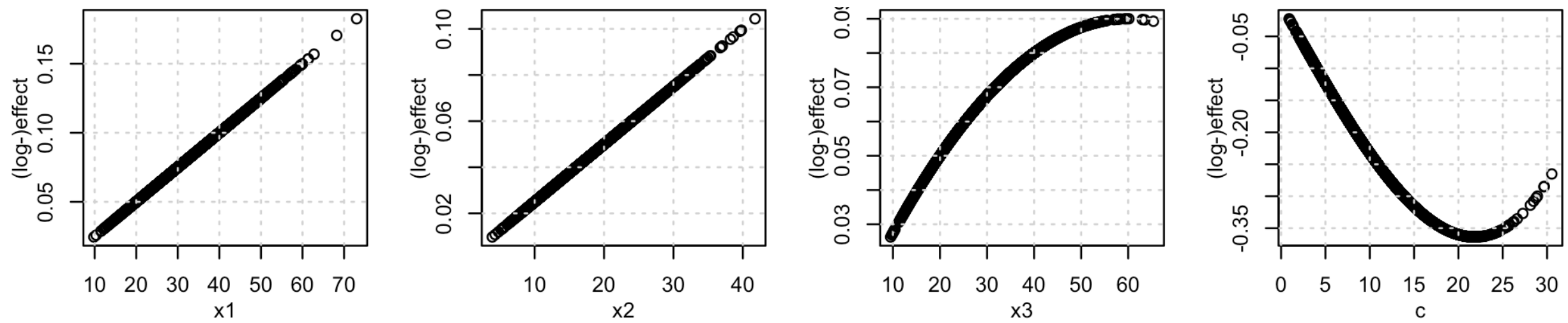

**Figure S2. Effects of $X_1$, $X_2$, $X_3$, and $C$ on the outcome used for simulations, represented as dose-response curves, based on one simulation sample.**

Regarding outcome model types,

1. Binary outcome generating models,

$$P(Y = 1) = \frac{1}{1 + \exp\left(-(\alpha + 0.01X_1 + 0.01X_2 + 0.012X_3 - 0.0001X_3^2 - 0.1C + 0.00007C^3 + 0.03V + 0.5B)\right)}$$

- Low baseline risk scenario: $\alpha = -5$, ensuring the majority of $P_i(Y = 1)$ as ≤0.1 (less favorable approximation condition, See Appendix S4.3)
- Moderate baseline risk scenario: $\alpha = -3.5$, ensuring the majority of $P_i(Y = 1)$ as $0.1 \leq P_i(Y = 1) \leq 0.2$ (the least favorable approximation condition compared to the other two baseline risk scenarios)
- High baseline risk scenario: $\alpha = -1$, ensuring the majority of $P_i(Y = 1)$ as $0.2 \leq P_i(Y = 1) \leq 0.8$ (favorable approximation condition, See Appendix S4.3)

2. Poisson count generating models,

$$\log(E[Y]) = \alpha + 0.025X_1 + 0.025X_2 + 0.03X_3 - 0.000025X_3^2 - 0.25C - 0.00000175C^3$$

- Low baseline risk scenario: $\alpha = 0$, ensuring the majority of $Y$ as <10

- Moderate baseline risk scenario: $\alpha = 3$, ensuring the majority of $P_i(Y = 1)$ as $10 \leq P_i(Y = 1) \leq 50$
- High baseline risk scenario: $\alpha = 6$, ensuring the majority of $P_i(Y = 1)$ as $350 \leq P_i(Y = 1) \leq 500$ (the log-linear function can be linearly approximated)

**S6.2. Conditional AEE Estimation**

MRC models were fitted using natural cubic splines for $X_1^{ep}$, $X_2^{ep}$, $X_3^{ep}$, $X_3^{ep^2}$, $C^{ep}$, $C^{ep^3}$, and $V^{ep}$ and dummy variables for $B$. AEE estimation was conducted as follows:

1) Fit a Q-model, using $X_{1,RC}^{ep}$, $X_{2,RC}^{ep}$, $X_{3,RC}^{ep}$, ${X_{3,RC}^{ep}}^2$, $C_{RC}^{ep}$, $C_{RC}^{ep^3}$, $V_{RC}^{ep}$, and $B$ as covariates.
2) Estimate conditional AEE.
    a. For logistic models,

For $X_1$,

$$\hat{P}\left(Y = 1 \middle| X_{1,RC}^{ep} = x_{1,RC}^{ep}, X_{2,RC}^{ep} = 0, X_{3,RC}^{ep} = 0, \boldsymbol{z}^{**} = 0\right) = \frac{1}{1 + \exp\left(-\left(\hat{\alpha} + \hat{\beta}_1 x_{1,RC}^{ep}\right)\right)}$$

$$RD = \hat{P}\left(Y = 1 \middle| X_{1,RC}^{ep} = x_{1,RC}^{ep}, X_{2,RC}^{ep} = 0, X_{3,RC}^{ep} = 0, \boldsymbol{z}^{**} = 0\right) - \hat{P}\left(Y = 1 \middle| X_{1,RC}^{ep} = mean\left(X_{1,RC}^{ep}\right), X_{2,RC}^{ep} = 0, X_{3,RC}^{ep} = 0, \boldsymbol{z}^{**} = 0\right)$$

$$RR = \frac{\hat{P}\left(Y = 1 \middle| X_{1,RC}^{ep} = x_{1,RC}^{ep}, X_{2,RC}^{ep} = 0, X_{3,RC}^{ep} = 0, \boldsymbol{z}^{**} = 0\right)}{\hat{P}\left(Y = 1 \middle| X_{1,RC}^{ep} = mean\left(X_{1,RC}^{ep}\right), X_{2,RC}^{ep} = 0, X_{3,RC}^{ep} = 0, \boldsymbol{z}^{**} = 0\right)}$$

$$OR = \exp(\hat{\beta}_1)$$

For $X_2$,

$$\hat{P}\left(Y = 1 \middle| X_{1,RC}^{ep} = 0, X_{2,RC}^{ep} = x_{2,RC}^{ep}, X_{3,RC}^{ep} = 0, \boldsymbol{z}^{**} = 0\right) = \frac{1}{1 + \exp\left(-\left(\hat{\alpha} + \hat{\beta}_2 x_{2,RC}^{ep}\right)\right)}$$

$$RD = \hat{P}\left(Y = 1 \middle| X_{1,RC}^{ep} = 0, X_{2,RC}^{ep} = x_{2,RC}^{ep}, X_{3,RC}^{ep} = 0, \boldsymbol{z}^{**} = 0\right)$$

$$- \hat{P}\left(Y = 1 \middle| X_{1,RC}^{ep} = 0, X_{2,RC}^{ep} = mean(X_{2,RC}^{ep}), X_{3,RC}^{ep} = 0, \boldsymbol{z}^{**} = 0\right)$$

$$RR = \hat{P}\left(Y = 1 \middle| X_{1,RC}^{ep} = 0, X_{2,RC}^{ep} = mean(X_{2,RC}^{ep}), X_{3,RC}^{ep} = 0, \boldsymbol{z}^{**} = 0\right)$$

$$OR = \exp(\hat{\beta}_2)$$

For $X_3$,

$$\hat{P}\left(Y = 1 \middle| X_{1,RC}^{ep} = 0, X_{2,RC}^{ep} = 0, X_{3,RC}^{ep} = x_{3,RC}^{ep}, \boldsymbol{z}^{**} = 0\right)$$

$$= \frac{1}{1 + \exp\left(-\left(\hat{\alpha} + \hat{\beta}_3 x_{3,RC}^{ep} + \hat{\beta}_4 {x_{3,RC}^{ep}}^2\right)\right)}$$

$$RD = \hat{P}\left(Y = 1 \middle| X_{1,RC}^{ep} = 0, X_{2,RC}^{ep} = 0, X_{3,RC}^{ep} = x_{3,RC}^{ep}, \boldsymbol{z}^{**} = 0\right)$$

$$- \hat{P}\left(Y = 1 \middle| X_{1,RC}^{ep} = 0, X_{2,RC}^{ep} = 0, X_{3,RC}^{ep} = mean(X_{3,RC}^{ep}), \boldsymbol{z}^{**} = 0\right)$$

$$RR = \frac{\hat{P}\left(Y = 1 \middle| X_{1,RC}^{ep} = 0, X_{2,RC}^{ep} = 0, X_{3,RC}^{ep} = x_{3,RC}^{ep}, \boldsymbol{z}^{**} = 0\right)}{\hat{P}\left(Y = 1 \middle| X_{1,RC}^{ep} = 0, X_{2,RC}^{ep} = 0, X_{3,RC}^{ep} = mean(X_{3,RC}^{ep}), \boldsymbol{z}^{**} = 0\right)}$$

$$OR_1 = \exp(\hat{\beta}_3)$$

$$OR_2 = \exp(\hat{\beta}_4)$$

For $C$,

$$\hat{P}\left(Y = 1 \middle| X_{1,RC}^{ep} = 0, X_{2,RC}^{ep} = 0, X_{3,RC}^{ep} = 0, C_{RC}^{ep} = c_{RC}^{ep}, \boldsymbol{z}^{**-C} = 0\right)$$

$$= \frac{1}{1 + \exp\left(-\left(\hat{\alpha} + \hat{\beta}_5 c_{RC}^{ep} + \hat{\beta}_6 {c_{RC}^{ep}}^3\right)\right)}$$

$$RD = \hat{P}\left(Y = 1 \middle| X_{1,RC}^{ep} = 0, X_{2,RC}^{ep} = 0, X_{3,RC}^{ep} = 0, C_{RC}^{ep} = c_{RC}^{ep}, \boldsymbol{z}^{**-C} = 0\right)$$

$$- \hat{P}\left(Y = 1 \middle| X_{1,RC}^{ep} = 0, X_{2,RC}^{ep} = 0, X_{3,RC}^{ep} = 0, C_{RC}^{ep} = mean\left(C_{RC}^{ep}\right), \boldsymbol{z}^{**-C} = 0\right)$$

$$RR = \frac{\hat{P}\left(Y = 1 \middle| X_{1,RC}^{ep} = 0, X_{2,RC}^{ep} = 0, X_{3,RC}^{ep} = 0, C_{RC}^{ep} = c_{RC}^{ep}, \boldsymbol{z}^{**-C} = 0\right)}{\hat{P}\left(Y = 1 \middle| X_{1,RC}^{ep} = 0, X_{2,RC}^{ep} = 0, X_{3,RC}^{ep} = 0, C_{RC}^{ep} = mean\left(C_{RC}^{ep}\right), \boldsymbol{z}^{**-C} = 0\right)}$$

$$OR_1 = \exp(\hat{\beta}_5)$$

$$OR_2 = \exp(\hat{\beta}_6)$$

b. For log-linear models,

For $X_1$,

$$E\left[Y \middle| X_{1,RC}^{ep} = x_{1,RC}^{ep}, X_{2,RC}^{ep} = 0, X_{3,RC}^{ep} = 0, \boldsymbol{z}^{**} = 0\right] = \exp\left(\hat{\alpha} + \hat{\beta}_1 x_{1,RC}^{ep}\right)$$

$$RD = E\left[Y \middle| X_{1,RC}^{ep} = x_{1,RC}^{ep}, X_{2,RC}^{ep} = 0, X_{3,RC}^{ep} = 0, \boldsymbol{z}^{**} = 0\right]$$
$$- E\left[Y \middle| X_{1,RC}^{ep} = mean\left(X_{1,RC}^{ep}\right), X_{2,RC}^{ep} = 0, X_{3,RC}^{ep} = 0, \boldsymbol{z}^{**} = 0\right]$$

$$RR = \exp\left(\hat{\beta}_1\right)$$

For $X_2$,

$$E\left[Y \middle| X_{1,RC}^{ep} = 0, X_{2,RC}^{ep} = x_{2,RC}^{ep}, X_{3,RC}^{ep} = 0, \boldsymbol{z}^{**} = 0\right] = \exp\left(\hat{\alpha} + \hat{\beta}_2 x_{2,RC}^{ep}\right)$$

$$RD = E\left[Y \middle| X_{1,RC}^{ep} = 0, X_{2,RC}^{ep} = x_{2,RC}^{ep}, X_{3,RC}^{ep} = 0, \boldsymbol{z}^{**} = 0\right]$$
$$- E\left[Y \middle| X_{1,RC}^{ep} = 0, X_{2,RC}^{ep} = mean\left(X_{2,RC}^{ep}\right), X_{3,RC}^{ep} = 0, \boldsymbol{z}^{**} = 0\right]$$

$$RR = \exp\left(\hat{\beta}_2\right)$$

For $X_3$,

$$E\left[Y \middle| X_{1,RC}^{ep} = 0, X_{2,RC}^{ep} = 0, X_{3,RC}^{ep} = x_{3,RC}^{ep}, \boldsymbol{z}^{**} = 0\right] = \exp\left(\hat{\alpha} + \hat{\beta}_3 x_{3,RC}^{ep} + \hat{\beta}_4 {x_{3,RC}^{ep}}^2\right)$$

$$RD = E\left[Y \middle| X_{1,RC}^{ep} = 0, X_{2,RC}^{ep} = 0, X_{3,RC}^{ep} = x_{3,RC}^{ep}, \boldsymbol{z}^{**} = 0\right]$$
$$- E\left[Y \middle| X_{1,RC}^{ep} = 0, X_{2,RC}^{ep} = 0, X_{3,RC}^{ep} = mean\left(X_{3,RC}^{ep}\right), \boldsymbol{z}^{**} = 0\right]$$

$$RR_1 = \exp\left(\hat{\beta}_3\right)$$

$$RR_2 = \exp\left(\hat{\beta}_4\right)$$

For $C$,

$$E\left[Y \middle| X_{1,RC}^{ep} = 0, X_{2,RC}^{ep} = 0, X_{3,RC}^{ep} = 0, C_{RC}^{ep} = c_{RC}^{ep}, \boldsymbol{z}^{**-C} = 0\right]$$

$$= \exp\left(\hat{\alpha} + \hat{\beta}_5 c_{RC}^{ep} + \hat{\beta}_6 {c_{RC}^{ep}}^3\right)$$

$$RD = E\left[Y \middle| X_{1,RC}^{ep} = 0, X_{2,RC}^{ep} = 0, X_{3,RC}^{ep} = 0, C_{RC}^{ep} = c_{RC}^{ep}, \boldsymbol{z}^{**-C} = 0\right]$$

$$- E\left[Y \middle| X_{1,RC}^{ep} = 0, X_{2,RC}^{ep} = 0, X_{3,RC}^{ep} = 0, C_{RC}^{ep} = mean\left(C_{RC}^{ep}\right), \boldsymbol{z}^{**-C} = 0\right]$$

$$RR_1 = \exp\left(\hat{\beta}_5\right)$$

$$RR_2 = \exp\left(\hat{\beta}_6\right)$$

**S6.3. Bias Analyses**

We tested the degree of bias in effect estimates from our new modeling framework.

Bias is defined as

$$RD - Bias(\%) \coloneqq \left(\frac{E\left[\widehat{RD}\right] - RD_{TRUE}}{RD_{TRUE}}\right) \times 100$$

$$RR - Bias(\%) \coloneqq \left(\frac{\mathrm{E}\left[\log\left(\widehat{RR}\right)\right] - \log(RR_{TRUE})}{\log(RR_{TRUE})}\right) \times 100$$

$$OR - Bias(\%) \coloneqq \left(\frac{\mathrm{E}\left[\log\left(\widehat{OR}\right)\right] - \log(OR_{TRUE})}{\log(OR_{TRUE})}\right) \times 100$$

As benchmarks, we also tested effect estimates based on true variables (i.e., $X_1$, $X_2$, $X_3$, $C$ and $V$) and effect estimates that were not adjusted for EME and CME (i.e., $X_1^{ep}$, $X_2^{ep}$, $X_3^{ep}$ , $C^{ep}$). For each scenario, 5,000 simulation samples were used.

The percentage of samples in which the 95% confidence intervals (CIs) included the true values, referred to as the coverage of the 95% confidence interval, was also calculated.

### S6.4. Results

Figure S3 shows the bias in regression coefficients of variables in logistic models as Q-models. Thus, the exponentiation of these coefficients corresponds to OR. When $X_1^{ep}$, $X_2^{ep}$, $X_3^{ep}$, and $C^{ep}$ (including ${X_3^{ep}}^2$ and ${C^{ep}}^3$) were used, biases were substantial in the simulation scenarios. For example, the biases for $X_1^{ep}$, $X_2^{ep}$, $X_3^{ep}$ were nearly -75% to -100% when the error was moderate or high, meaning that conditional AEE estimates were severely underestimated. In some circumstances, biases exceeded -100% meaning that the direction of conditional AEE estimates was reversed. Biases could also be substantial in low error scenarios. For example, the biases regarding $C^3$ were >100%, meaning that if $C$ was the main exposure of interest, conditional AEE estimates could be severely overestimated. Figure S4 shows the bias in the regression coefficients of variables in log-linear models as Q-models, exhibiting similar patterns.

For both outcome models, when MRC were applied, the biases were greatly reduced (Figures S3 and S4). Particularly, if the logistic and log-linear functions lie within the range where they can be favorably approximated linearly, the biases were nearly eliminated. In simulation settings, this was the case in the high-baseline risk scenario, followed by the low-baseline risk scenario. See Figures S5–7 for conditional RR and RD estimates from logistic models and RD from log-linear models, respectively. The consequences were similar. The coverages of 95% CIs were ~95% (Figures S8–9) when the bias was minimal.

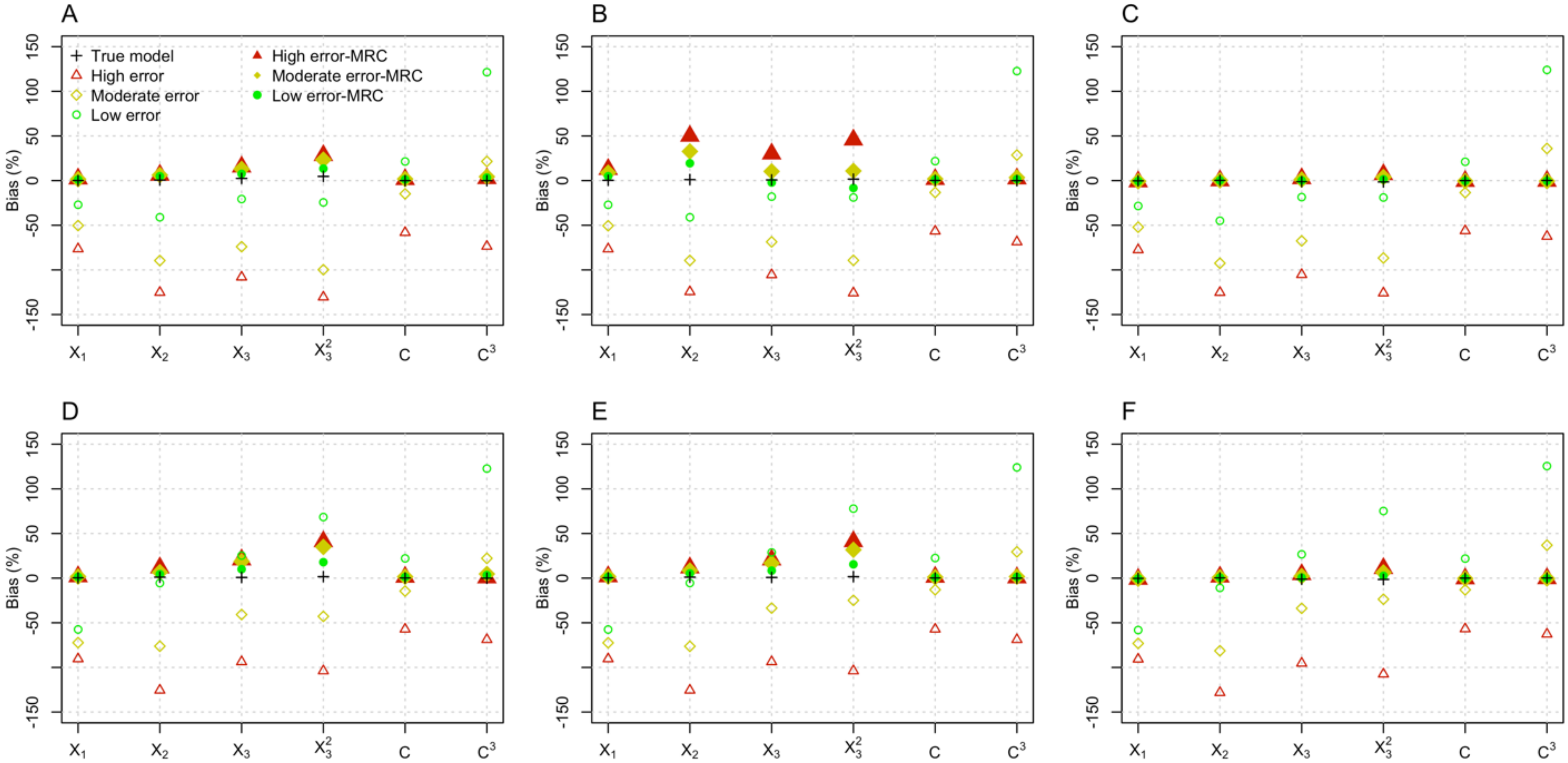

**Figure S3. Biases in regression coefficient estimates from logistic models as Q-models (i.e., conditional log-OR) with MRC.**

Note 1. The top panels (A, B, C) are for additive errors scenarios. The bottom panels (D, E, F) are for multiplicative error scenarios. Panels A and D are for 'Low baseline risk' scenarios; Panels B and E are for 'Moderate baseline risk' scenarios; and Panels C and F are for 'High baseline risk' scenarios. 'High error', 'Moderate error', and 'Low error' indicate scenarios for different magnitude of errors, when $X_1^{ep}$, $X_2^{ep}$, $X_3^{ep}$, and $C^{ep}$ (including ${X_3^{ep}}^2$ and ${C^{ep}}^3$) were used to estimate Q-models. 'True model" indicates when $X$, $X_2$, $X_3$, and $C$ (including $X_3^2$ and $C^3$) were used to estimate Q-models.

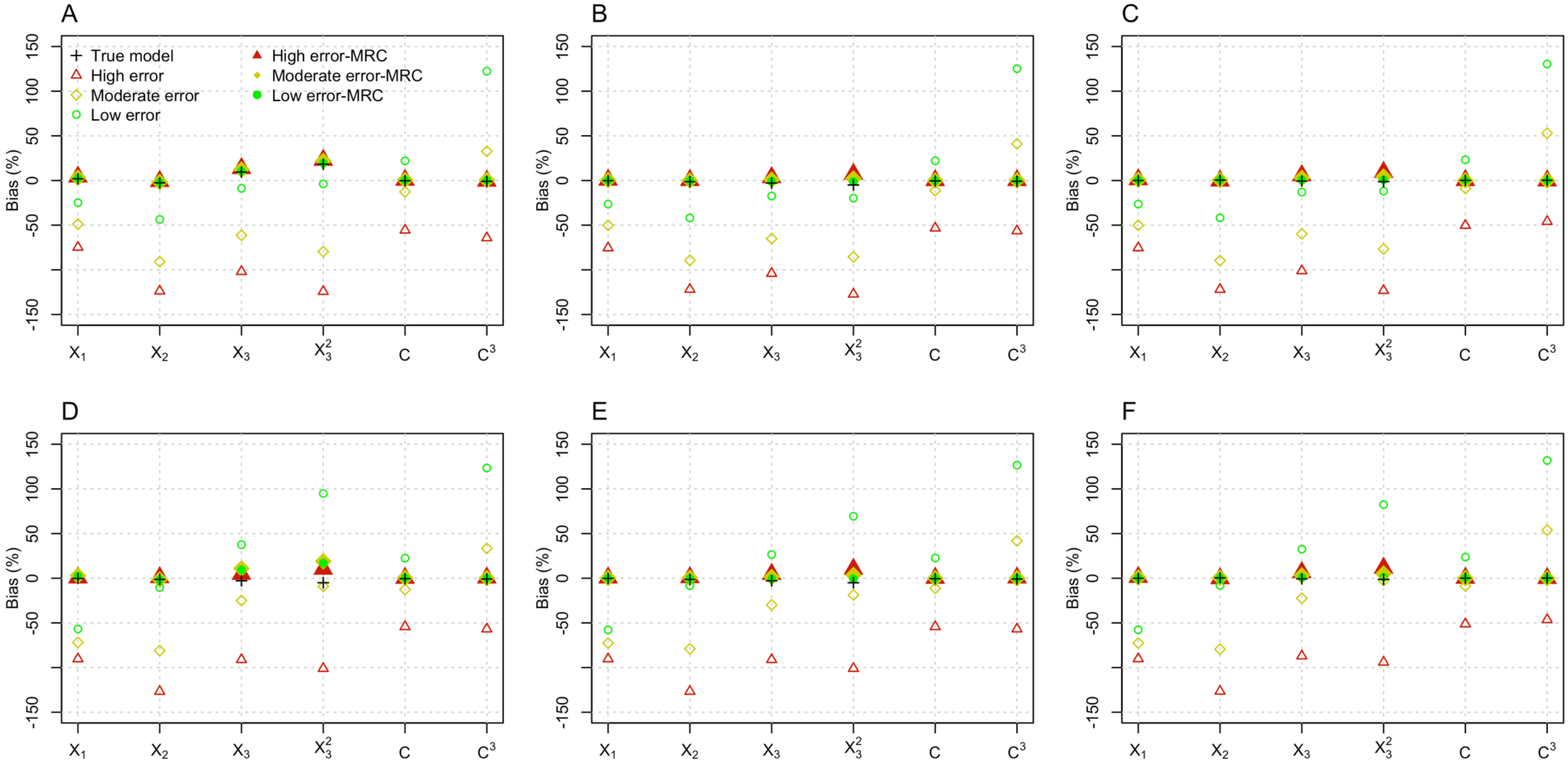

**Figure S4. Biases in regression coefficients estimates from log-linear models as Q-models (i.e., log-RR (rate ratio)) with MRC.**

Note 1. The top panels (A, B, C) are for additive errors scenarios. The bottom panels (D, E, F) are for multiplicative error scenarios. Panels A and D are for 'Low baseline risk' scenarios; Panels B and E are for 'Moderate baseline risk' scenarios; and Panels C and F are for 'High baseline risk' scenarios. 'High error', 'Moderate error', and 'Low error' indicates scenarios when $X_1^{ep}$, $X_2^{ep}$, $X_3^{ep}$, and $C^{ep}$ (including ${X_3^{ep}}^2$ and ${C^{ep}}^3$) were used. 'True model" indicates when $X$, $X_2$, $X_3$, and $C$ (including $X_3^2$ and $C^3$) were used.
Note 2. As the outcome model is log-linear (Poisson), regression coefficients indicate log-rate ratio (RR).

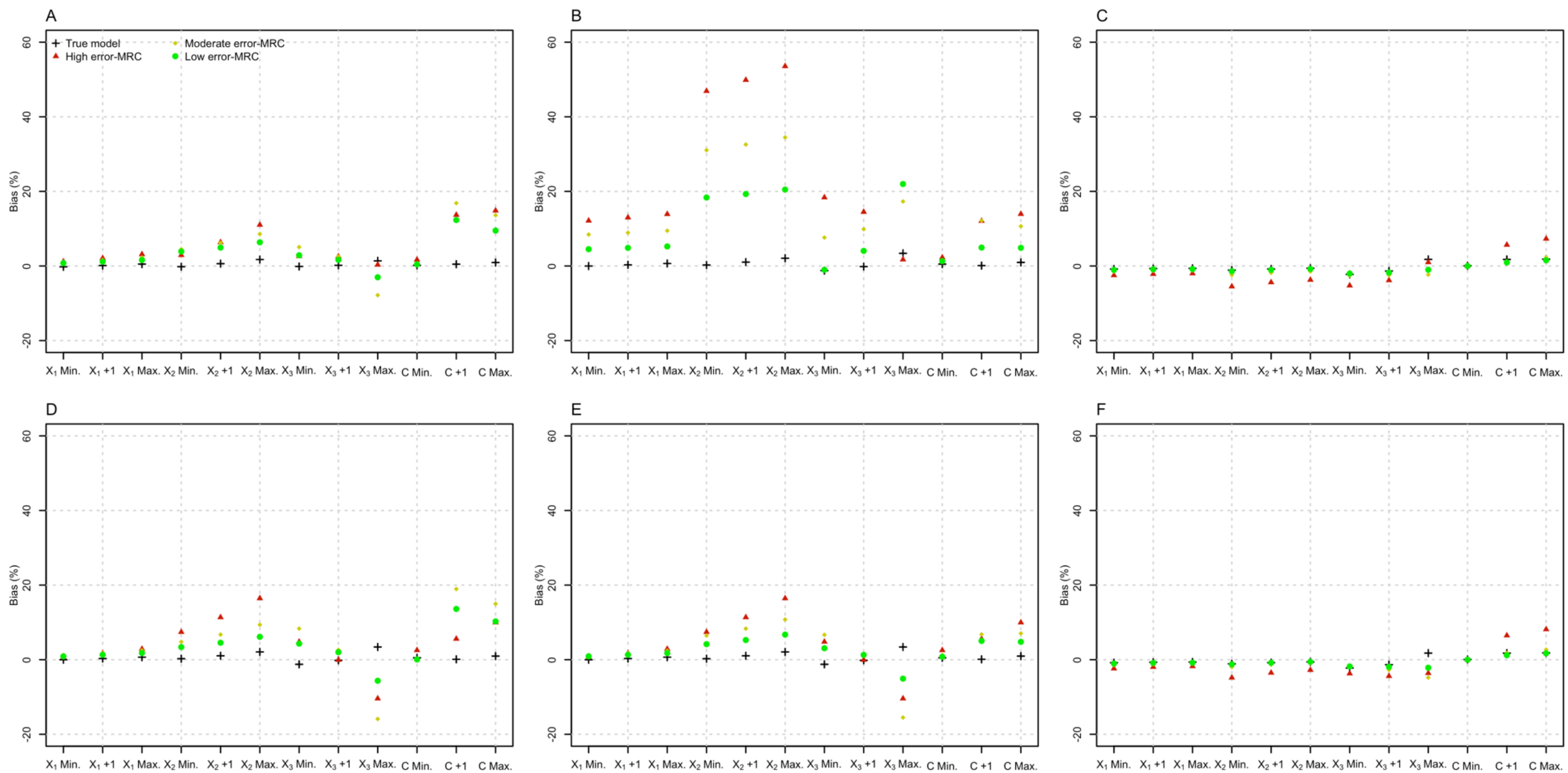

**Figure S5. Biases in AEE estimates (log-RR scale) from logistic models as Q-models with MRC.**

Note 1. The top panels (A, B, C) are for additive errors scenarios. The bottom panels (D, E, F) are for multiplicative error scenarios. Panels A and D are for 'Low baseline risk' scenarios; Panels B and E are for 'Moderate baseline risk' scenarios; and Panels C and F are for 'High baseline risk' scenarios. 'High error', 'Moderate error', and 'Low error' indicates scenarios when $X_1^{ep}$, $X_2^{ep}$, $X_3^{ep}$, and $C^{ep}$ (including ${X_3^{ep}}^2$ and ${C^{ep}}^3$) were used. 'True model" indicates when $X$, $X_2$, $X_3$, and $C$ (including $X_3^2$ and $C^3$) were used.
Note 2. As the outcome model is logistic, the effect of all covariates on the outcome in RR scale becomes non-linear Therefore, biases are presented for estimates of the effect of a covariate at its near-minimum level (Min.) and maximum level (Max.), using the mean level as a reference. Biases for estimates of the effect of a covariate at the mean level +1 (+1) are also presented.

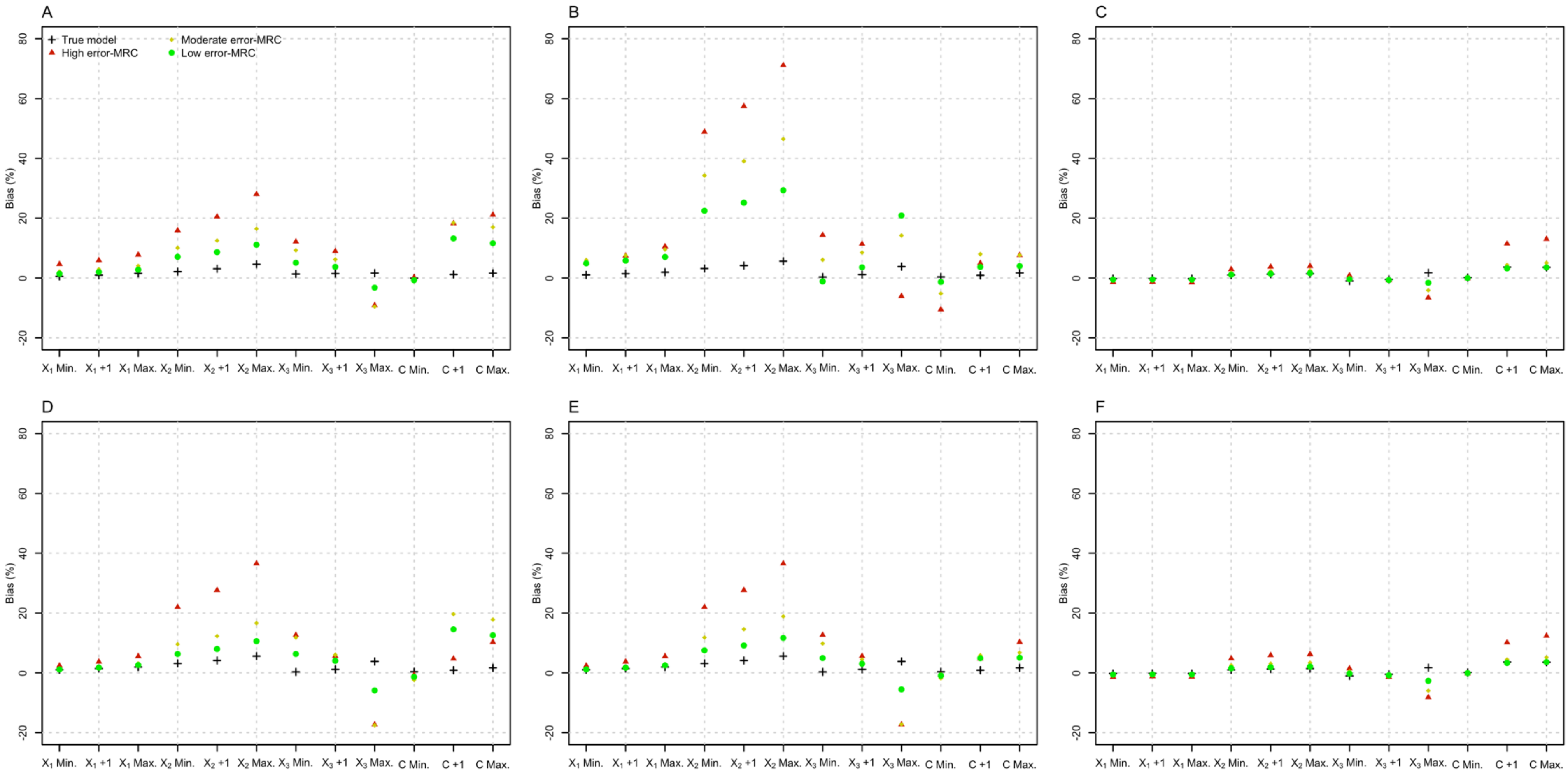

**Figure S6. Biases in AEE estimates (RD scale) from logistic models as Q-models with MRC.**

Note 1. The top panels (A, B, C) are for additive errors scenarios. The bottom panels (D, E, F) are for multiplicative error scenarios. Panels A and D are for ‘Low baseline risk’ scenarios; Panels B and E are for ‘Moderate baseline risk’ scenarios; and Panels C and F are for ‘High baseline risk’ scenarios. ‘High error’, ‘Moderate error’, and ‘Low error’ indicates scenarios when $X_1^{ep}$, $X_2^{ep}$, $X_3^{ep}$, and $C^{ep}$ (including ${X_3^{ep}}^2$ and $C^{ep^3}$) were used. ‘True model” indicates when $X$, $X_2$, $X_3$, and $C$ (including $X_3^2$ and $C^3$) were used.
Note 2. As the outcome model is logistic, the effect of all covariates on the outcome in RD scale becomes non-linear. Therefore, biases are presented for estimates of the effect of a covariate at its near-minimum level (Min.) and maximum level (Max.), using the mean level as a reference. Biases for estimates of the effect of a covariate at the mean level +1 (+1) are also presented.

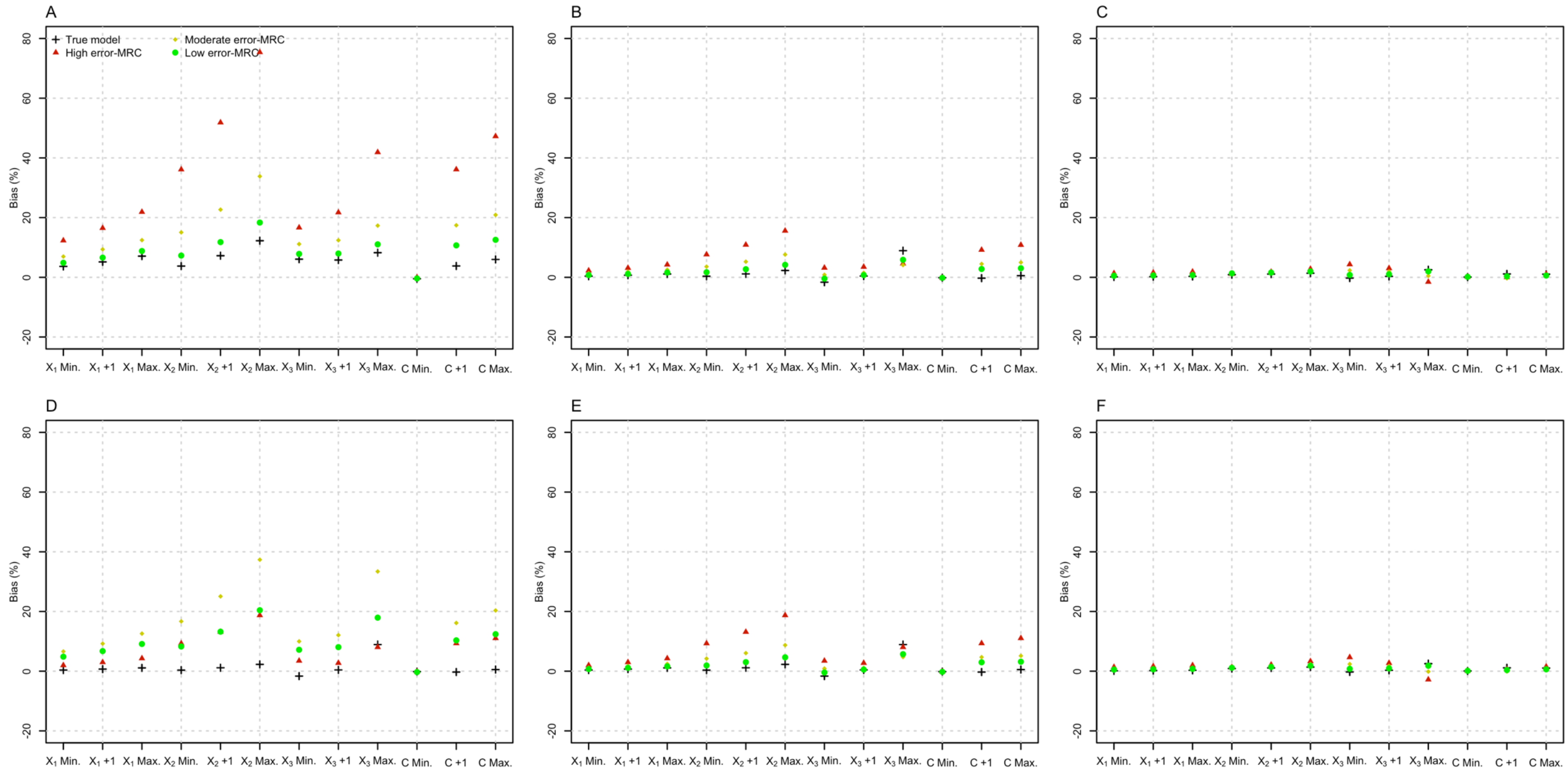

**Figure S7. Biases in AEE estimates (RD scale) from log-linear models as Q-models with MRC.**

Note 1. The top panels (A, B, C) are for additive errors scenarios. The bottom panels (D, E, F) are for multiplicative error scenarios. Panels A and D are for ‘Low baseline risk’ scenarios; Panels B and E are for ‘Moderate baseline risk’ scenarios; and Panels C and F are for ‘High baseline risk’ scenarios. ‘High error’, ‘Moderate error’, and ‘Low error’ indicates scenarios when $X_1^{ep}$, $X_2^{ep}$, $X_3^{ep}$, and $C^{ep}$ (including $X_3^{ep^2}$ and $C^{ep^3}$) were used. ‘True model” indicates when $X$, $X_2$, $X_3$, and $C$ (including $X_3^2$ and $C^3$) were used.
Note 2. As the outcome model is log-linear, the effect of all covariates on the outcome in RD scale becomes non-linear. Therefore, biases are presented for estimates of the effect of a covariate at its near-minimum level (Min.) and maximum level (Max.), using the mean level as a reference. Biases for estimates of the effect of a covariate at the mean level +1 (+1) are also presented.

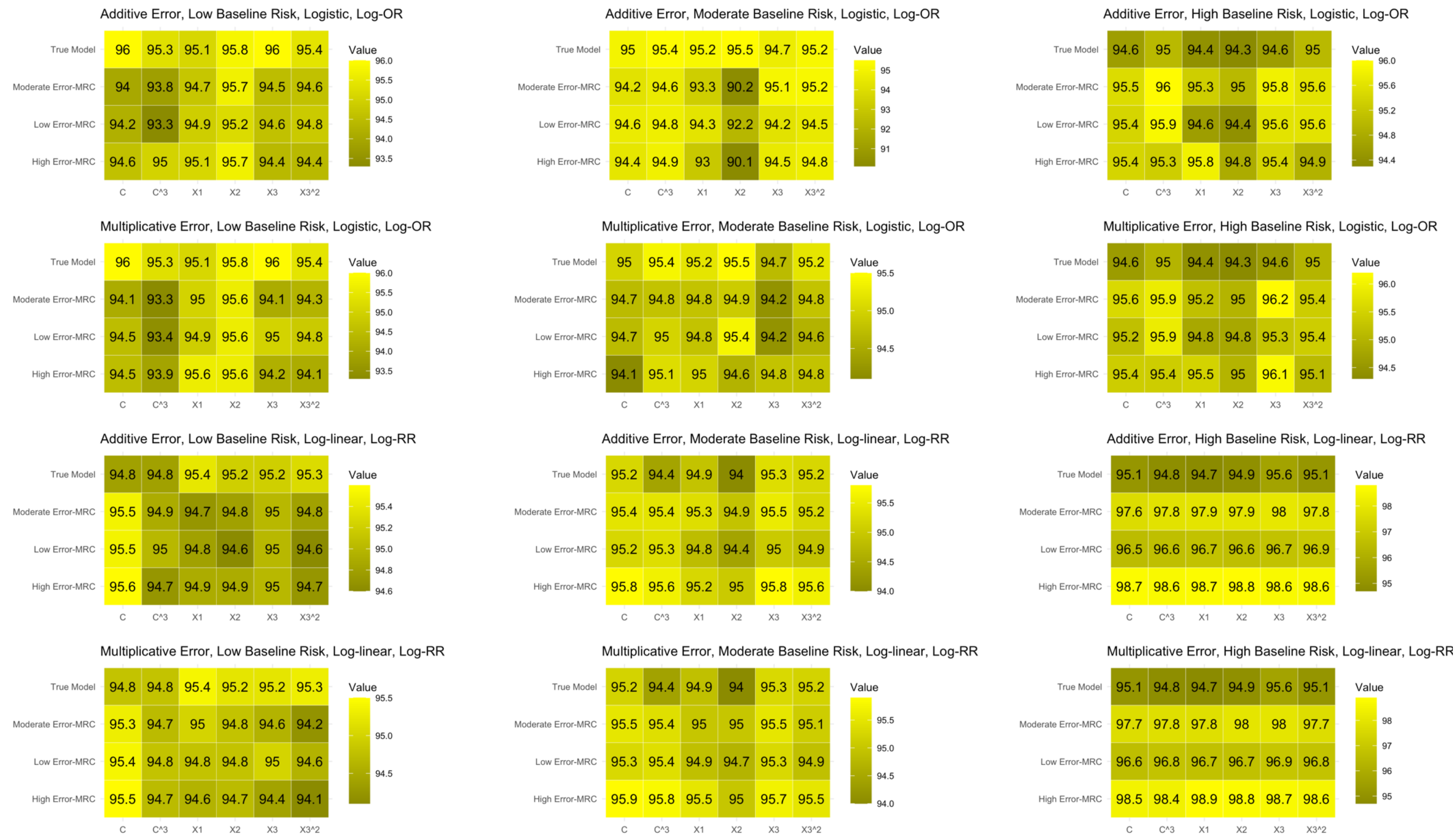

**Figure S8. Nominal coverage of 95% confidence intervals in AEE estimates in the Log-OR scale from logistic models and the Log-RR (rate ratio) scale from log-linear models.**

Note. The coverages were nearly 95% when there were small/negligible biases in AEE estimates.

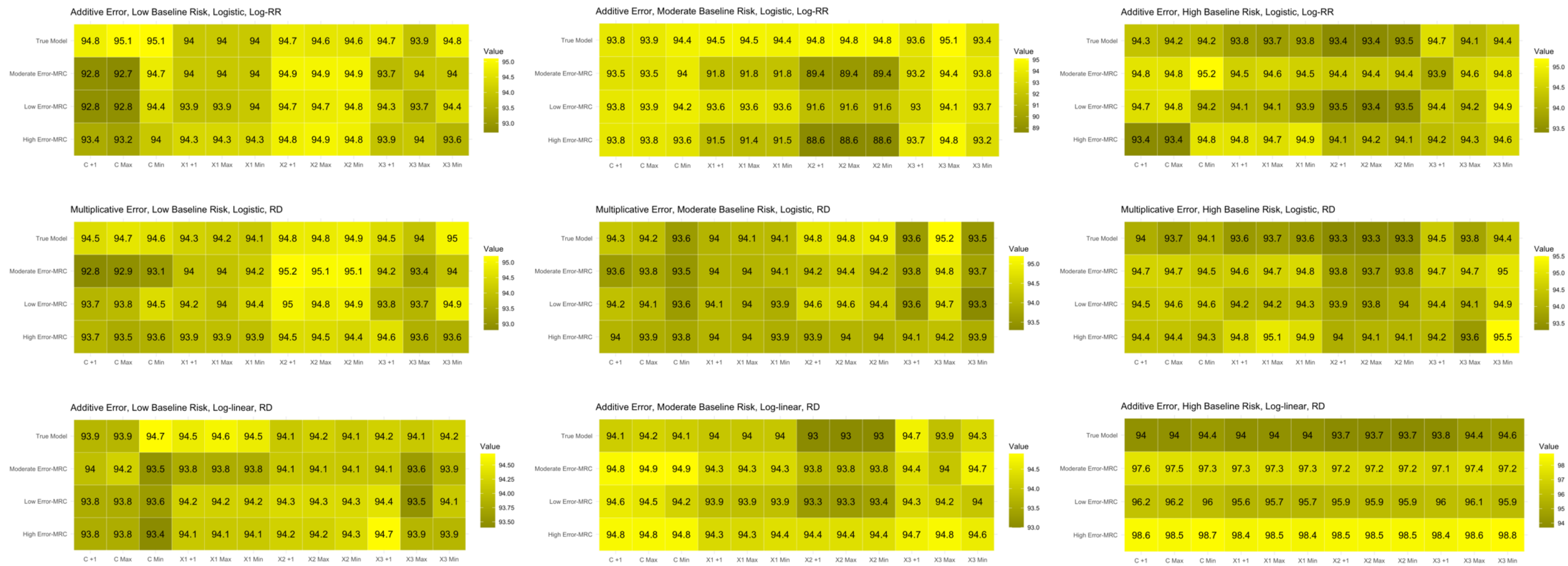

**Figure S9. Nominal coverage of 95% confidence intervals in AEE estimates in the Log-RR (risk ratio) scale from logistic models and the RD scale from logistic and log-linear models.**

Note 1. As the effect of all covariates on the outcome becomes non-linear. Results are presented for estimates of the effect of a covariate at its near-minimum level (Min.) and maximum level (Max.), using the mean level as a reference. Results for the effect of a covariate at the mean level +1 (+1) are also presented.
Note 2. The coverages were nearly 95% when there were small/negligible biases in AEE estimates.

## Appendix S7. Development of MRC-Like ML-based Estimates and Validation Data

### S7.1 Ground-Based Measurements for $PM_{2.5}$, $O_3$, Temperature, and Humidity

We obtained data of the US EPA's 18 $PM_{2.5}$ monitoring stations and 15 $O_3$ monitoring stations in Cook County and the adjacent counties in Illinois for the period, Jan. 2020–Feb. 2021. $PM_{2.5}$ levels were reported as 24-hour averages, while $O_3$ levels were daily 8-hour moving average maximums. 8-hour moving average maximum is a health-relevant exposure metric. A few small negative values were treated as 0.1 $\mu g/m^3$ (slightly below the minimum value) for $PM_{2.5}$ and as 1 ppb (slightly below the minimum value) for $O_3$. The Parameter Occurrence Code in the datasets was evaluated, and daily values included only one POC code. No unrealistic values were found.

For daily temperature and dew temperature, we obtained data from the US EPA's two monitoring stations and US NOAA's 9 Global Summary of Day (GSOD) stations. Relative humidity was estimated using daily temperature and dew temperature. We also obtained 11 Global Historical Climatological Network (GHCN) stations that measure daily temperature (no dew temperature). No unrealistic values were found.

### S7.2. Training and Validation Datasets Construction

The validation dataset was constructed for two purposes: 1) identifying the accuracy of MRC-Like ML-based estimates of environmental conditions; and 2) constructing MRC models.

Figure S10 presents the monitoring locations in Cook County and the neighboring counties in Illinois, of US EPA's Air Quality System and NOAA's GSOD and GHCN. GSOD stations measure temperature and dew (relative humidity were estimated using temperature and dew point). GHCN stations measure temperature but not dew point. The validation data include the six locations, referred to as "co-monitoring locations", which are marked in the figure.

- Location #1 has $PM_{2.5}$, $O_3$, temperature, and relative humidity measurements by US EPA (while "+" symbol is used in the figure for simplicity).
- Location #2 has $PM_{2.5}$ and $O_3$ measurements by US EPA, but there are no meteorological stations. There are four GHCN stations around this location. Temperature measurements from these four stations were highly positively correlated with the measurements from another, with a regression slope close to 1. So, they were averaged. Because relative humidity values from three GSOD stations located at (42.120, -87.905), (42.417, -87.867), and (41.914, -88.246) were highly positively correlated with values from another, with a regression slope close to 1, they were averaged.
- Location #3 has $PM_{2.5}$ and $O_3$ measurements by US EPA, but there are no meteorological stations. Temperature and relative humidity measurements from the three nearest GSOD stations were averaged because values from one location were highly positively correlated with those from another, with a regression slope close to 1.
- Location #4 has $PM_{2.5}$ and $O_3$ measurements by US EPA. The nearby GSOD station was used.

- Location #5 has $PM_{2.5}$ and $O_3$ measurements by US EPA. The nearby GSOD station was used.
- Location #6 has $PM_{2.5}$ and $O_3$ measurements by US EPA. Two GSOD stations located at (41.786, -87.752) and (41.540, -87.532) and the GHCN station located at (41.738, -87.777) was used. Temperature and relative humidity values were averaged because values from one location were highly positively correlated with those from the other, with a regression slope close to 1.

The training datasets for $PM_{2.5}$, $O_3$, relative humidity, and temperature modeling were constructed by excluding the stations used for the validation data. Specifically,

- For $PM_{2.5}$ and $O_3$ modeling, the training datasets included all remaining $PM_{2.5}$ and $O_3$ stations, respectively.
- For relative humidity modeling, the training dataset included all remaining GSOD stations.
- For temperature modeling, the training dataset included all remaining GSOD and GHCN stations, excluding GHCN stations collocated with GSOD stations.

Scatter plots and correlation analyses for each monitoring location pair demonstrated good agreement, with Pearson correlation coefficients ranging from 0.45 to 0.99 for $PM_{2.5}$, 0.67 to 0.97 for $O_3$, 0.98 to 0.999 for temperature, and 0.73 to 0.97 for relative humidity. This suggests that while some heterogeneous spatiotemporal patterns exist in environmental conditions, the overall distribution may not substantially differ in a way that may compromise the representativeness of the validation data.

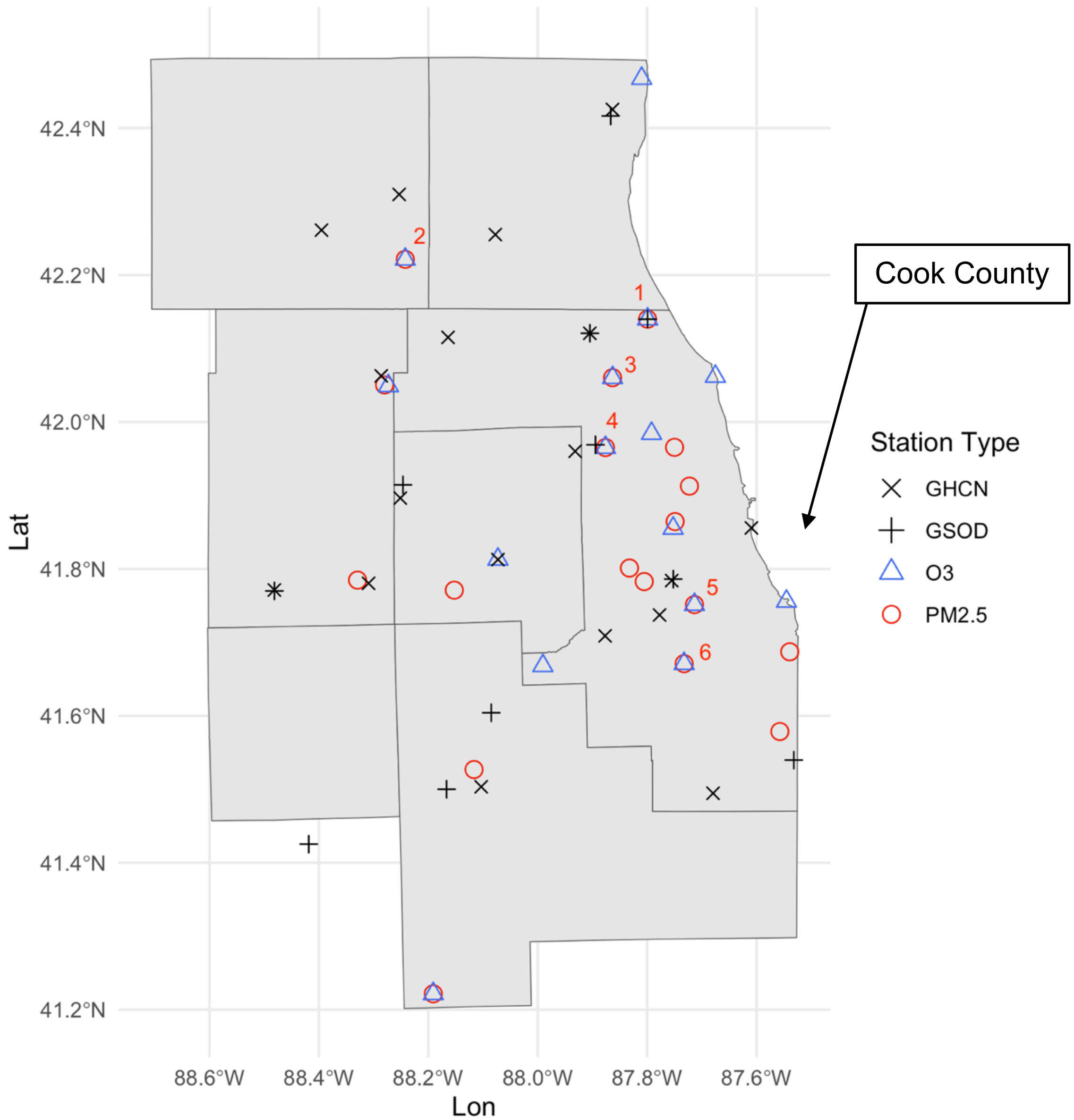

**Figure S10. Locations of ground-based monitoring stations.**

Note. the numbers 1–6 indicate the locations of the monitors used for the validation data

### S7.3. Data Collection and Pre-Processing for Machine Learning-Based Prediction

Table S2 presents datasets and data processing for predictor variables. For traffic data, we obtained traffic volume estimates on each road segment in Illinois from Illinois Department of Transportation (IDOT), which are stationary. Traffic volumes in Cook County and the neighboring counties are constant after 2010, according to IDOT's annual reports, except for the year 2020 and 2021 due to the COVID-19 pandemic.

To obtain hourly traffic volume estimates, we extracted yearly day-of-the-week-specific traffic volume statistics and yearly hour-specific traffic volume statistics from IDOT's annual report and multiplied respective fractions to the traffic volume estimates on each road segment, and then calculated daily averages. Furthermore, to obtain more accurate traffic volume estimates during COVID-19 pandemic, we obtained the Google's COVID-19 Community Mobility Report data for Cook County, which provides county-level daily variation of the individual mobility and is estimated based on cell phone data and multiplied respective fractions to the daily averages of traffic volumes on each road segment.

The Enhanced vegetation Index (EVI) based on MODIS of Terra and Aqua satellites is available every 8-day. EVI values exhibited strong (sub-)seasonality. We estimated daily EVI values with grid-specific generalized additive models including temporal smoothers (i.e., interpolation).

When an original/raw dataset contains missing values at the spatial scale, which was not common, we used focal statistics to fill in these. When the original dataset contained missing values at the temporal scale, which was common for Sentinel satellite-based aerosol indices, we calculated averages using the values from adjacent days for each grid cell.

We obtained 12km gridded, 24-hour daily average estimates for $PM_{2.5}$ and daily 8-hour moving average maximum estimates for $O_3$ from US EPA's CMAQ as predictors. EPA provides downscaled estimates at census tracts, referred to as FAQSD, but we did not use this. The reason is that these downscaled estimates are variables obtained from prediction models without confounders and based on EPA's ground-based measurements including those we will use as validation data. For illustration, according to their methodology, downscaled estimates were created by 'fusing' measurements from ground-based monitors with 12km gridded estimates from CMAQ. For the sake of understanding, this may be expressed as, acknowledging that this might be oversimplification,

$$Measured\ Pollutant_{L,t} = \eta_0 + s\big(CMAQ\ estimate_{L,t}, coordinates\ of\ L, t\big) + residuals_{L,t}$$

and downscaled estimates (DCE) are

$$DCE_{L,t} \coloneqq \widehat{\eta_0} + \hat{s}\big(CMAQ\ estimate_{L,t}, coordinates\ of\ L, t\big)$$

This downscaler model takes a similar form of RC but does not include confounders to adjust for in outcome modeling. This may introduce bias as noted in Appendix S1.3.

For temperature and humidity, we used 4 km-gridded daily estimates of temperature and relative humidity from the PRISM Climate Model for our application analysis, as in Kim et al. (2022)[18]'s work and 9 km-gridded daily estimates from ERA5 Land Hourly. To the best of my knowledge, the details of PRISM and ERA5's methodology are not fully publicly available. Consequently, we cannot confirm whether NOAA's ground-based measurements were directly incorporated in the development of the PRISM Climate/ERA5 Land Hourly models. By visualizing scatter plots and estimating correlations between PRISM/ERA5 temperature estimates and ground-based measurements, we found that measurements from a few monitors, not all, were nearly perfectly correlated with PRISM Climate Model-based estimates and measurements from a few monitors, not all, were nearly perfectly correlated with ERA5 Land Hourly model-based estimates. For relative humidity, nearly perfect correlations were not found.

All datasets were resampled at 500m grids and then spatiotemporally consolidated.

**Table S2. Overview of original data and processing for developing MRC-Like ML-based estimates**

| Datasets | Original temporal resolution | Original spatial resolution | Processing & Final temporal resolution | Processing & Final spatial resolution | Relevance |
|---|---|---|---|---|---|
| **Atmospheric-air pollution-related indices/estimates** | | | | | |
| Sentinel-5P TROPOMI UV index / $O_3$ | Daily | ~5km | Daily | Resampled at 500m | Aerosol index/Correlated with emission sources/secondary particle formation |
| OMI $NO_2$ / $O_3$ | Daily | 0.25° (~27.8km * 20km) | Daily | Resampled at 500m | Aerosol index/Correlated with traffic emission sources/secondary particle formation (for sub-regional patterns) |
| CAMS O3, $PM_{2.5}$, HCHO, $NO_2$, $HNO_3$, NO, $OH^-$, | 3-hourly | 0.75° (~62km * 83km; one grid cell mostly covers Cook County) | Daily | Resampled at 500m | The same as above, but mainly for regional temporal patterns. |
| US EPA's CMAQ $PM_{2.5}$ and $O_3$ estimates | Daily (24-hour average for PM2.5; 8-hour moving average maximum for $O_3$) | 12km | Daily | Resampled at 500m | $PM_{2.5}$ and $O_3$ estimates |
| **Other atmospheric and land variables estimate (e.g., meteorological factors)** | | | | | |
| PRISM Temperature and Relative humidity estimates (Relative humidity were estimated from temperature and dew point temperature estimates) | Daily | 4km | Daily | Resampled at 500m | Temperature and relative humidity estimates |
| ERA5 Land Hourly Temperature, Relative humidity (estimated from temperature and dew), total evaporation, evaporation from the vegetation, total precipitation, wind speeds (U and V components) | Hourly | 0.1° | Daily | Resampled at 500m | Chemical reaction conditions; dry/wet deposition; dispersion; interaction with land features |
| **Land use/cover** | | | | | |
| NLDS Land Cover (developed areas with low/medium/high density; open water, etc.) | One year (2021) | 30m | One year (2021), used to predict spatial heterogeneity | Squared buffers as focal sums (i.e., 100m, 250m, 500m, 1km, 2.5km, 5km, and 10km) and then resampled at 500m | Human activity related to anthropogenic emissions; Pollutants travel long distance (decay over distance); chemical reactions in the atmosphere |
| NLDS Impervious Cover (primary, secondary, tertiary road) | One year (2021) | 30m | One year (2021), used to predict spatial heterogeneity | Squared buffers as focal sums (i.e., 100m, 250m, 500m, 1km, 2.5km, 5km, and 10km) and then resampled at 500m | Traffic-related emissions; Pollutants travel long distance (decay over distance); emissions and activities may differ by road types. |
| NLDS Tree Canopy Cover (TCC) | One year (2021) | 30m | One year (2021), used to predict spatial heterogeneity; TCC is slightly different between years | Squared buffers as focal sums (i.e., 100m, 250m, 500m, 1km, 2.5km, 5km, and 10km) and then resampled at 500m | Dry deposition; spatially correlated with lower anthropogenic emission sources. Biogenic emissions regarding $O_3$. |
| Terra/Aqua MODIS Enhance Vegetation Index | Every 8-day | 250m | Daily (interpolation using grid-specific generalized additive models with temporal smoothers) | Resampled at 500m | Dry deposition; spatially correlated with lower anthropogenic emission sources. Biogenic emissions regarding $O_3$. |
| **Traffic-related information** | | | | | |
| Illinois Department of Transportation Annual Average Daily Traffic | Stationary values, updated every *X* years (varies by road segment) | Line segment | Daily | Weighted sums in different spatial scales (i.e., 100m, 250m, 500m, 1km, 2.5km, 5km, and 10m) and then resampled at 500m | Traffic emissions. Pollutants travel long distance (decay over distance). |
| **Topology** | | | | | |
| USGS Elevation (elevation, slope, and a heat load index estimated) | One year (2023) | 1/3 arc-second (~10m) | One year used to predict spatial heterogeneity | Resampled at 500m | Dispersion, hyper local heat/temperature |
| Distance to Lake Michigan | One year (2023) | Line segment (Calculated using polygon files) | One year used to predict spatial heterogeneity | Resampled at 500m | Different local atmospheric environment (e.g., wind, precipitation, evaporation, condensation) |
| **Other** | | | | | |
| Google's COVID-19 Community Mobility Report | Daily (2020–2021) | County-level (Cook County) | Daily | Identical daily values for all areas | Changes in population-level emission activity during COVID-19 pandemic |

### S7.4. Model Training

Using XGBoost, we fit ML models for $PM_{2.5}$, $O_3$, temperature, and relative humidity, respectively. Using these models, we estimated 500m gridded daily estimates.

We constructed a training dataset for each environmental condition using ground-based monitoring data that were not used to construct the validation dataset (Appendix S7.2). Using this, we trained models. The evaluation metric was root mean squared error. The objective function was tweedie regression with log-link for $PM_{2.5}$ and $O_3$; regression with squared loss for temperature and humidity. Several hyperparameters were selected using a combination of manual search and cross-validation. Specifically, maximum iterations were based on cross-validation. Values for *Max.depth*, *Eta*, *Subsample*, *Colsample_bytree* were selected manually. We selected a subset of the predictor variables using *Gain* and *Frequency* from initially trained models. Using selected variables, we re-trained models.

### S7.5. Model Validation

We evaluated the agreement between modelled estimates, referred to as MRC-Like ML-based estimates and ground-based measurements using the validation dataset. The validation dataset was not used in model training. Linear regressions were fitted, including location- (site-) specific ones.

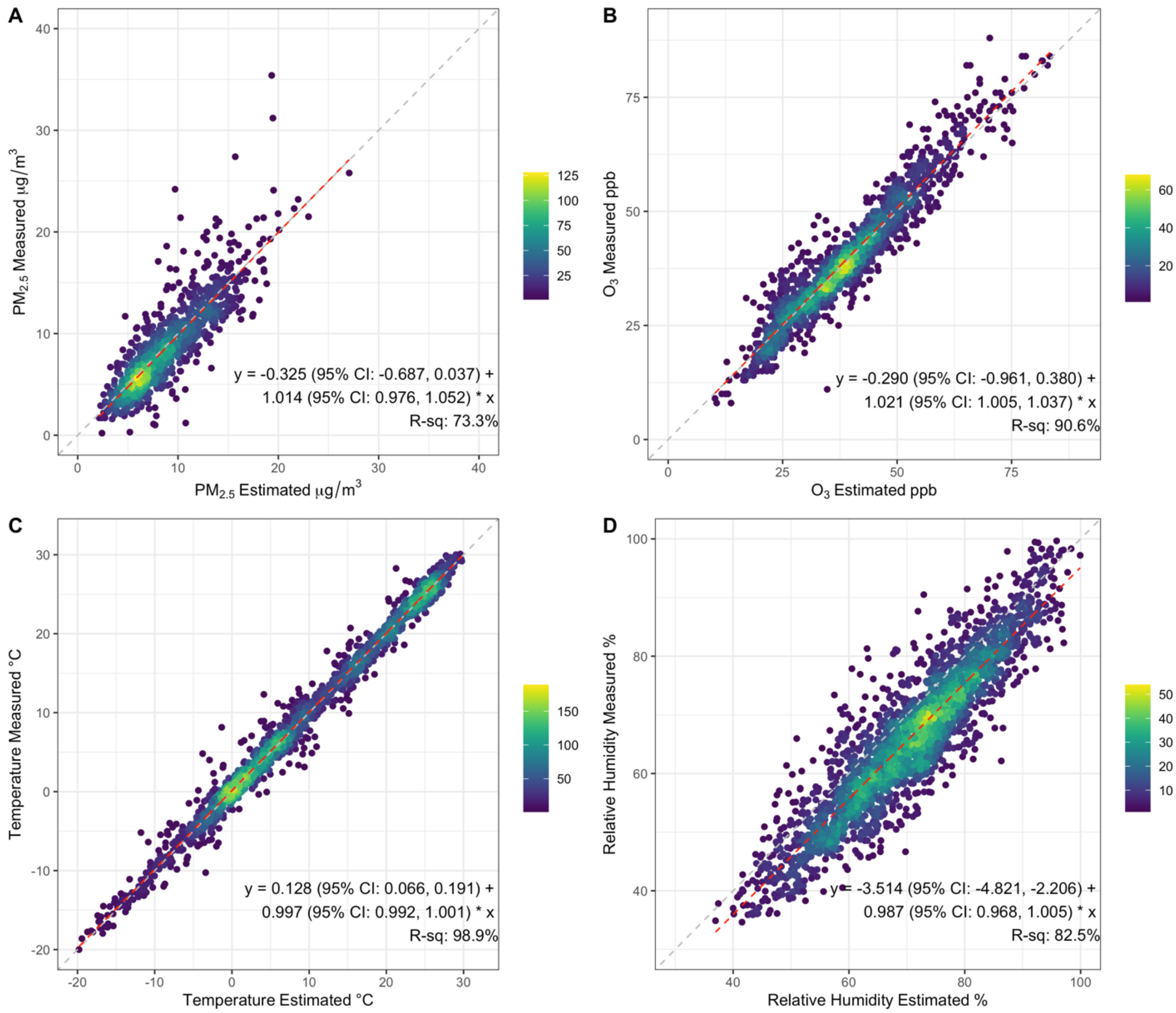

**Figure S11. Scatter plots for MRC-Like ML based estimates of environmental conditions and ground-based measurements in the validation data. A) $PM_{2.5}$; B) $O_3$, C) Temperature, and D) Relative humidity.**

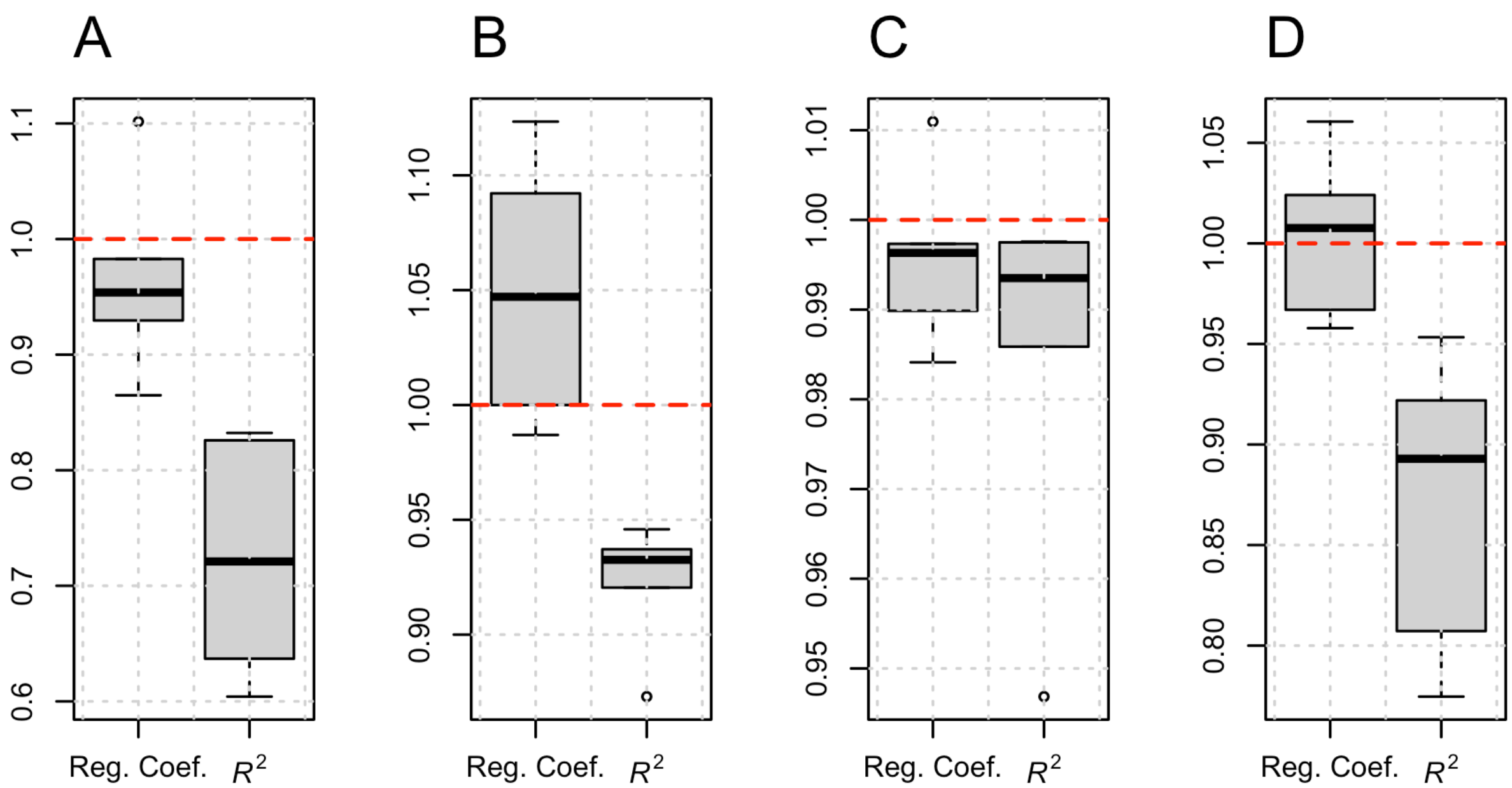

**Figure S12. Monitoring location-specific regression coefficient of ground measurements against MRC-Like ML-based estimates and $R^2$ in the validation data. A) $PM_{2.5}$; B) $O_3$, C) Temperature, and D) Relative humidity.**

Note. The validation data includes six locations; the red dashed lines indicate 1, which is the ideal value.

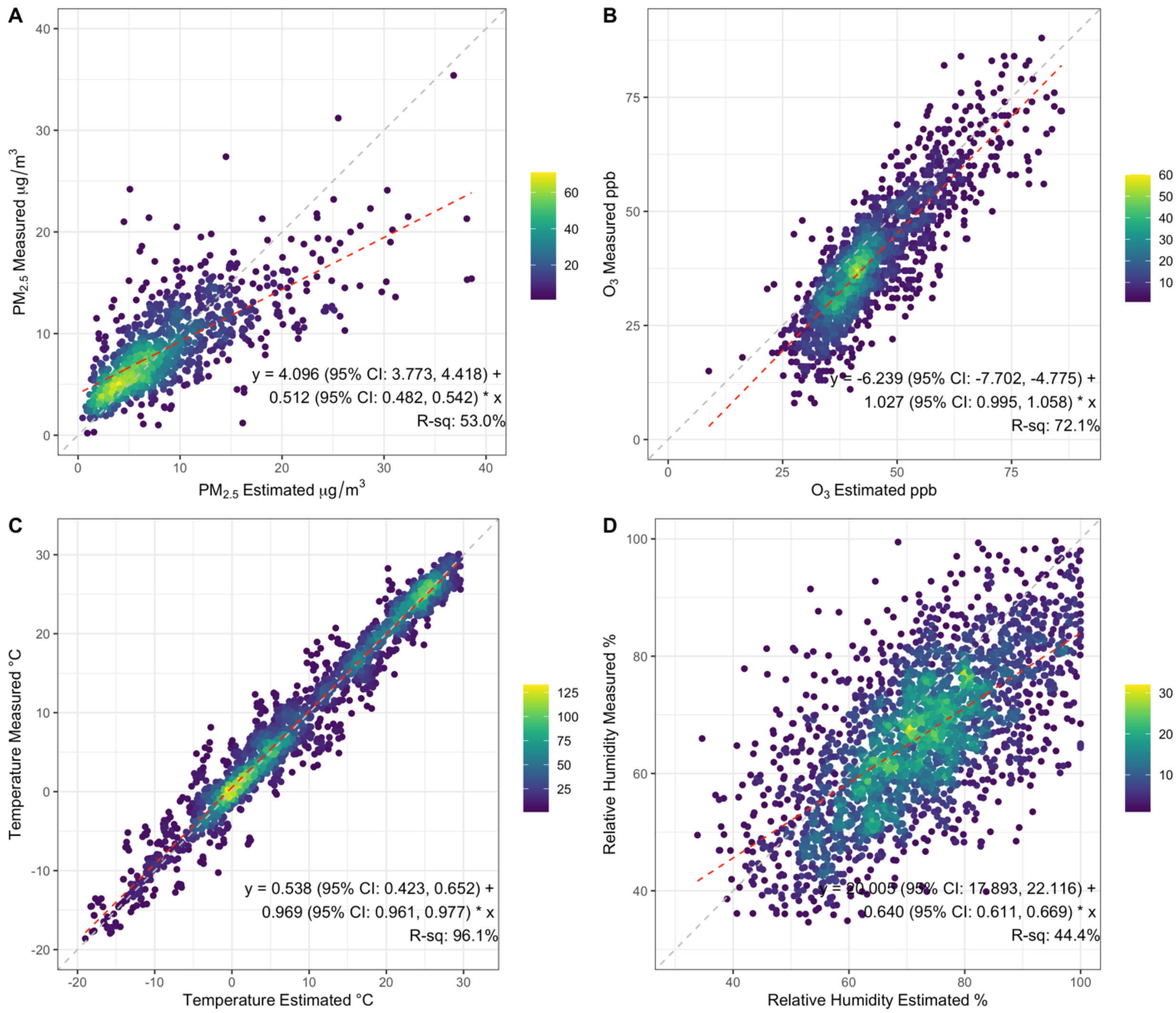

**Figure S13. Scatter plots for CMAQ estimates of $PM_{2.5}$ and $O_3$ and PRISM estimates of temperature and relative humidity and ground-based measurements in the validation data. A) $PM_{2.5}$; B) $O_3$, C) Temperature, and D) Relative humidity.**

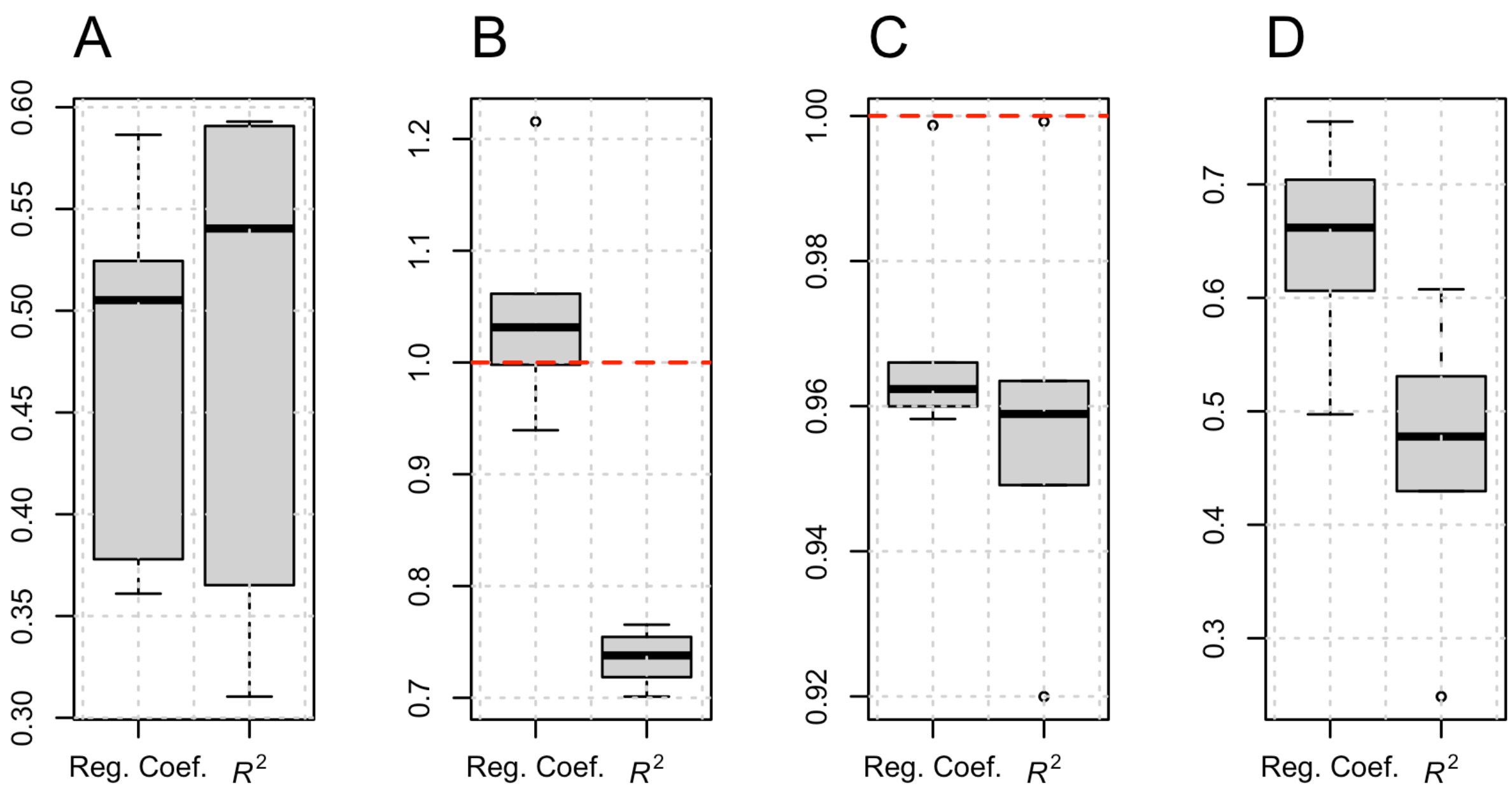

**Figure S14. Monitoring location-specific regression coefficient of ground measurements against CMAQ/PRISM-based estimates and $R^2$ in the validation data. A) $PM_{2.5}$; B) $O_3$, C) Temperature, and D) Relative humidity.**

Note. The validation data includes six locations; The red dashed lines indicate 1, which is the ideal value.

## Appendix S8. Error Sources When Using Modelled Estimates

We provide a conceptualization of EME and CME sources when using modelled estimates of ambient environmental conditions (Figure S15, Go three pages forward), aiming to help readers understand how they were addressed and why the difference between personal exposure and ambient levels may be seen as a component that shapes the population average effect rather than an error component or how this may be viewed as an instrumental variable perspective.

### S8.1. Error from Prediction Modeling

Exposure/confounder prediction models may have inherent uncertainty[19-21] (Source 1 in Figure S15). From a practical standpoint, it may be helpful to characterize exposure/confounder prediction models to help construct validation data. There are two widely used types of air pollution (and meteorological) models: statistical/machine-learning-based models (Source 1A) and assimilation/simulation models (Source 1B). The motivation of this classification is that investigators should understand the data used to generate exposure/confounder estimates in order to construct validation data that is needed to address bias due to measurement error, which was considered in our application analysis. Although this classification is not necessarily precise as there are many variations in the environmental exposure modeling field (e.g., assimilation/simulation model-based estimates used as predictors for statistical/machine learning models), the key distinction lies in whether ground-based measurements that may serve as ground-truth were used for building or training models, (not for validating). Investigators who wish to apply measurement error-correction methods should construct validation data that was not used to train prediction models, meaning that validation data should not include measurements used for model training, as we did in our application analysis.

Statistical or machine learning-based prediction models typically estimate the relationship between a dependent variable–often obtained from ground-based measurements from monitoring stations as the ground truth–and numerous independent variables that may help predict the dependent variable using various algorithms (e.g., traditional statistical algorithms, machine learning algorithms). In the machine learning field, the dependent variable is often referred to as labels, while the independent variables may be called features or classifiers, depending on the context.

Assimilation/simulation models may be constructed based on emission rates, chemical reactions in the atmosphere, dispersion, the fate and transport of air pollutants, and other relevant factors drawing from knowledge environmental engineering and atmospheric science. A series of mathematical equations, with input data plugged into these equations, are used to predict environmental conditions. Constructing these models may not incorporate measurements from ground-based monitors. While it is difficult to fully understand the complexity of atmospheric processes in the real world, these models can provide accurate estimates, albeit with some inevitable errors.

Sometimes, initial predicted values from assimilation/simulation models may be calibrated using statistical methods and measurements (e.g., satellite remote sensing-based measurements and ground-based measurements). Also, once assimilation/simulation models are constructed, the initial modeled estimates may be downscaled to obtain more spatially granular estimates, such as the US Environmental Protection Agency (EPA)'s Fused Air Quality Surface Using Downscaling (FAQSD) models. Downscaling can be performed using statistical modeling, where the dependent variable is often ground-based measurements. The independent variables typically include modeled estimates from assimilation/simulation models, as well as other spatial variables or statistical terms (e.g., spatial smoothers) that can be leveraged to obtain more granular estimates. These calibration and downscaling task may be viewed as a type of statistical modeling (Source 1A). For statistical and machine learning modeling, estimates from assimilation/simulation models may be used as additional predictors in addition to other predictors.

If modeled spatiotemporal estimates are correctly linked to health data, the resulting error may be Berksonian but may not necessarily be purely Berksonian as noted earlier in Appendix S1.3 (e.g., Eq. S1-1, Eq.1-2, Eq.1-3). If modelled estimates are only spatial, which is not uncommon in air pollution exposure models[22,23], non-Berkson error type may additionally arise due to the lack of information about time-varying exposures.

As modelled estimates are inherently error-prone, they should be corrected. We developed spatiotemporal estimates for $PM_{2.5}$, $O_3$, temperature, and humidity using our MRC-Like ML modeling and further corrected these using MRC modeling for our application analysis.

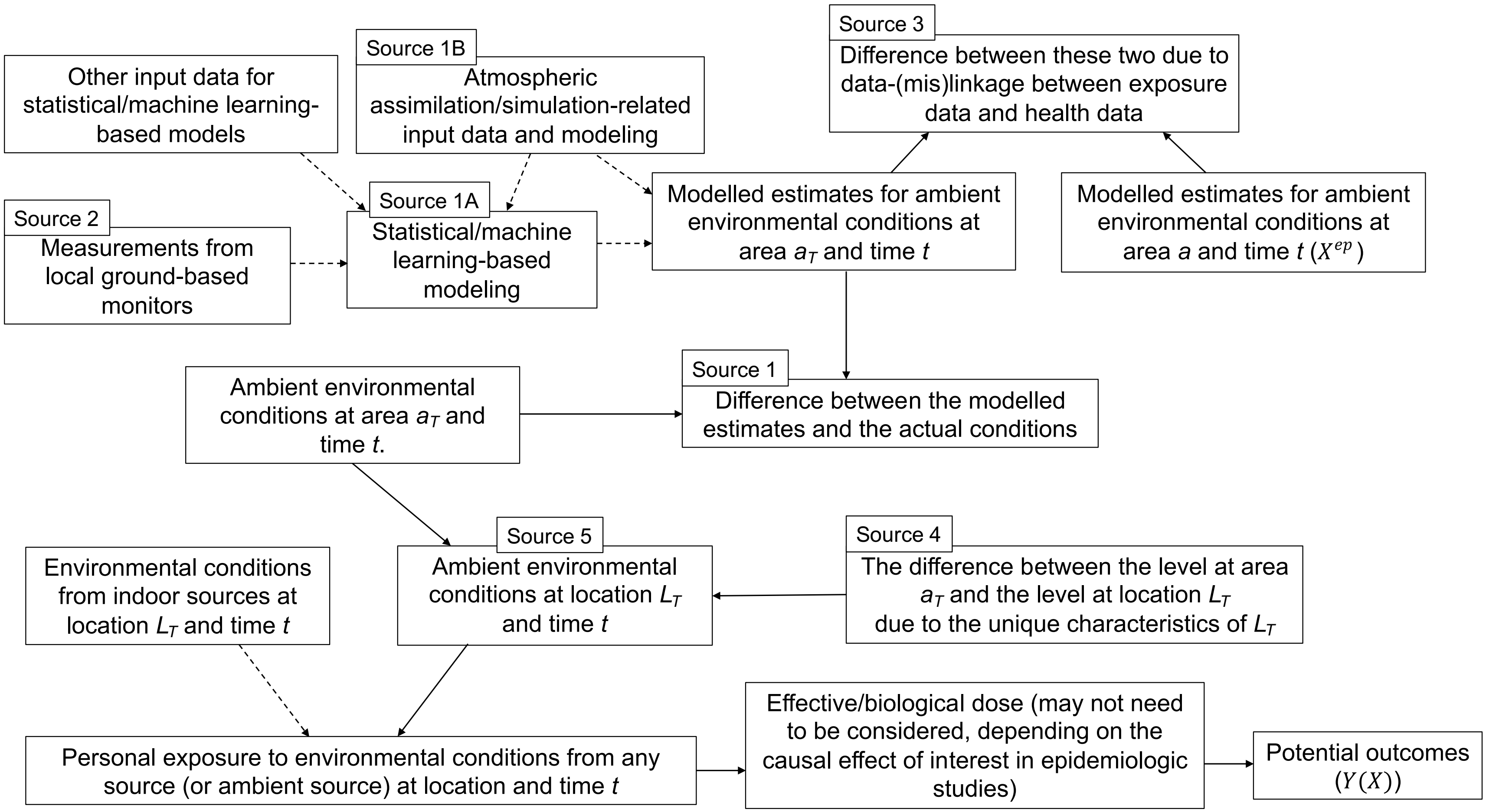

**Figure S15. A conceptualization of sources of EME and CME in epidemiological studies with modelled estimates of ambient air pollutants and meteorological factors.**

Note. Dashed lines denote those that may or may not arise depending on the research setting. The subscript $T$ is used to indicate the "true" (actual) location of an individual at a specific time $t$ so that $L_T$ or $a_T$ may differ from $L$ and $a$ used for data-linkage, where $a$ (or $a_T$) denotes an area that subsumes a location $L$ (or $L_T$)

### S8.2. Instrumental Error in Ground-Based Measurements

In air pollution/temperature epidemiology or air pollution/temperature prediction modeling, ground-based measurements from regulatory monitors (e.g., Air Quality System of US EPA; monitor networks of US National Oceanic Atmospheric Administration (NOAA)) may be considered as ground-truth values. Technically, there may exist instrumental errors (e.g., random noise, sensor error) (Source 2 in Figure S15). We obtained the data that underwent quality control and assurance.

### S8.3. Potential Spatiotemporal Data Mislinkage

Another important error in air pollution epidemiological studies may arise from data linkage with inadequate information of individuals' mobility when exposure/confounder data are linked to health data (Source 3 in Figure S15), instead of direct exposure/confounder assessment. In the real-world, individuals typically move around, meaning their actual exposure to ambient environmental conditions may change over time if the locations they visit have different environmental conditions. It may be challenging to measure the mobility of individuals for ethical and technical constraints. So, investigators may link exposure/confounder data to health data using limited information about the mobility (e.g., residential address).

It is unclear how the error structure evolves as a result of this error. Understanding whether unmeasured mobility could introduce non-negligible bias to causal effect estimates may be complex and is contextual (i.e., study-specific). Depending on the exposure time-windows of interest, unmeasured mobility on a very short time scale (e.g., daily, hourly, daily) may not be concerning in studies for longer-term exposures, as such uncertainty may cancel out over a longer time span, which is an assumption. Sometimes, this assumption may not necessarily result in measurement error that may not introduce bias (e.g., pure Berkson error)[24]. When this assumption does not likely hold, if individual mobility is substantial *and* environmental conditions are heterogeneous across locations that individuals visit, unmeasured mobility may be concerning. The resulting error may be manifest in a different form depending on how investigators design the study[21,25] (e.g., population-level time-series analysis, cohort data analysis), including how exposure variables are created and operationalized[21]. To understand the error structure, investigators need to consider the spatiotemporal distribution of environmental conditions, as the error may vary depending on the condition of interest. While individual mobility may be spatiotemporally autocorrelated, errors may be spatially, temporally, or both spatially and temporally autocorrelated. For example, the spatiotemporal distribution of $PM_{2.5}$ may not be the same as that of $O_3$. Thus, the impact of unmeasured mobility on health effect estimates for $PM_{2.5}$ may not equal those for $O_3$. In addition to EME, because $PM_{2.5}$ or $O_3$ may be a confounder in estimating the effect of the other, CME due to unmeasured mobility may also become complex.

So, a practical question is how to consider the impact of unmeasured mobility on health effect estimates when mobility data are unavailable. For our application analysis, we attempted to identify potential bias due to unmeasured mobility. See Appendix S9.

**S8.4. Difference Between Area-Level Environmental Condition and Personal Exposure to the Condition**

Even if individuals had continued to stay at the given geocoded locations used for the data-linkage, this may not necessarily imply the absence of error. There exist at least two sources to consider: 1) the difference between area-level environmental condition (i.e., true value at an area $a_T$) and location-specific environmental condition (i.e., true value at a location $L_T$ within $a_T$) (Source 4 in Figure S15), where the subscript $T$ indicates the “true” (actual) location of an individual at a specific time, emphasizing that $L_T$ or $a_T$ may differ from $L$ (and $a$) used for data-linkage; and 2) the difference between the location-specific environmental condition and personal exposure to that environmental condition (Source 5).

Source 4 error may arise when investigators do not have spatially granular exposure data for location-specific exposure characterization but area-specific exposure data that is spatially coarse. A common approach is to assign spatially coarse data to specific locations. For example, if investigators use 4-km gridded estimates, then all individuals whose geocoded residential addresses fall within a given grid cell may be assigned the same value from the 4-km gridded estimate. When monitor-based measurements are used, an approach may be to assign the same measurement from a monitor value to individuals whose residential addresses are within *x*-meter of the monitor, referred to as a buffer of x-meter.

The error from Source 4 may be assumed to be purely Berksonian, as the environmental condition at different locations can be perceived as variations from the average level centered on an area that includes those locations. This area-level environmental condition may be perceived to influence the overall trend, while the remaining within-area variations (i.e., inter-location variations) arise from location-specific factors that may be perceived to be spatially distributed with a mean of zero. This purely Berkson error assumption is sometimes referenced by investigators, with many published articles, especially, in air pollution and temperature time-series studies. The error from Source 4 may be assumed to be purely Berksonian as

$$X_{i,L_T,t} = X_{a_T,t} + U_{i,L_T,t}$$

where $X_{i,L_T,t}$ denotes the true exposure of an individual, $i$ who stayed at $L_T$ and at time $t$. In practice, an issue of this assumption is that conceptualizing such an average level centered on an area may sometimes be challenging. For example, (individual-level) environmental epidemiologic studies may use administrative spatial units (e.g., city, ZIP code, census tract) for data-linkage. Administrative spatial units may not be uniformly shaped. Both $X_{a_T,t}$ and $U_{i,L_T,t}$ may exhibit spatiotemporal patterns (so that these two may be correlated), meaning that the error may take a Berksonian form but not be purely Berksonian[21,26]. Consequently, bias in health effect estimates may occur[15,26]. When individual mobility is unmeasured, investigators may wish to conceptualize $a$ that include possible locations around which individuals move, but this is not a simple task. Considering individual mobility, one study found an indication that the purely Berkson error assumption may not necessarily hold[24]. Nevertheless, a violation of this assumption itself may not necessarily compromise observed associations.

Understanding should be contextual and study-specific. The resulting bias depends on several factors: the extent of deviation from the purely Berkson error assumption, exposure time windows, study designs, statistical methods, and possibly net bias due to unmeasured or residual confounding among others. See Appendix S9 for how this error was considered in our application analysis.

Source 5 error arises from the difference between personal exposure to environmental condition and location-specific environmental condition. This difference is likely because the behavioral patterns of individuals (other than the mobility of individuals over different locations) may differ. For example, many people spend most of their time indoors, while others may have different patterns (e.g., outdoor workers). In studies using measurements or estimates of ambient air pollution levels, the difference between ambient and indoor environmental conditions should be noted. This difference may be decomposed into two: 1) personal exposure to environmental conditions from *ambient sources* and 2) that from *indoor sources*. The former is frequently focused as the essence of the causal effect in (ambient) air pollution epidemiologic studies[20,27] and temperature epidemiological studies[28]. For example, effect estimates could provide insight into how ambient air pollution levels relate to health risks (e.g., instrumental variable context[13]), thereby providing insights into regulations on ambient air pollution levels and guidance for individuals on behaviors to minimize exposure to ambient air pollution, although the magnitude of effect estimates based on ambient metrics may not equal what effect estimates would be based on personal exposure metrics. Furthermore, if the causal entity from ambient sources is not different from that from indoor sources, such as temperature and $O_3$–since these are identical physical and chemical entity regardless of the source—then causal effect estimates based on ambient environmental condition variables can reflect the health effects of the exposure from any sources although the magnitude of the effect estimates may differ (See below). $PM_{2.5}$ does not reflect a single physical or chemical substance, which is a composite of various particle sizes and chemical compositions.

Nevertheless, conceptualizing the causal effect of exposure to environmental conditions from ambient sources does not necessarily guarantee no impact from Source 5 error. There is a condition under which Source 5 error can introduce toward-the-null bias[25], which is analogous to what happens in an instrumental variable analysis. For example, for individuals staying indoors or moving between indoor and outdoor environments, a Source 5 error model may be expressed as

$$X_{personal\ exposure} = \gamma_{ambient} X_{ambient\ level} + \gamma_{indoor} X_{indoor\ source}$$

where $\gamma_{ambient}$ and $\gamma_{indoor}$ describe the contribution of ambient and indoor sources to personal exposure, respectively. For individuals staying outdoors, there may be no error from Source 5:

$$X_{personal\ exposure} = X_{ambient\ level}$$

If $\gamma_{indoor} X_{indoor\ source}$ could be seen as an independent random variable and $X_{ambient\ level}$ is independent of $X_{indoor\ source}$, then the resulting error model would take a form of Berkson error (not purely Berksonian):

$$X_{personal\ exposure} = \gamma_{ambient} X_{ambient\ level} + U$$

While this illustrated error model may be an oversimplification in some circumstances, this provides insights on how the error may operate. If individuals are assumed to spend most of time indoors, the infiltration/penetration factor that describes how much outdoor air pollution can penetrate the indoor environment may serve as a surrogate for $\gamma_{ambient}$, which generally <1[29]. Health effect estimates would be underestimated due to this EME if $0 < \gamma_{ambient} < 1$.

If there are unmeasured individuals risk factors (UIR) potentially correlated with $X_{personal\ exposure}$ through $X_{indoor\ source}$, but are independent of $X_{ambient\ level}$, then the use of $X_{ambient\ level}$ may be an instrumental variable analysis as noted by Weisskopf and Webster[13]. This is because of the absence of other pathway through which $X_{ambient\ level}$ can affect $Y$, its independence of any unmeasured confounders such as UIR, and monotonicity. This may further alternatively be viewed as an population average effect of $X_{ambient\ level}$, where this effect may be proportional to the true effect (say $\beta$), by a factor of $\gamma_{ambient}$ ($\beta\gamma_{ambient}$). Note that $\gamma_{ambient}$ may be study-population specific. Figure S1J may describe this, by replacing $X_a$ with $X_{ambient\ level}$ and $X_{personal\ exposure}$ with $X_i$. This implies that the impact of non-differential error in Source 5 on health effect estimates may not need to be considered when the population average effect of an air pollutant from ambient sources is of primary interest (i.e., $\beta\gamma_{ambient}$). This is analogous to the Local Average Treatment Effect (LATE) in causal inference. Regarding CME, if the error component is completely random and the difference between personal exposure to the confounder and the ambient level of the confounder is solely multiplicative similar to $\gamma_{ambient}$, then using the ambient level of the confounder would suffice for confounding adjustment.

**S8.5. Differential Error**

The relationship between ambient environmental condition and personal exposure to this condition (i.e., from ambient sources)–the latter is often unmeasured in large epidemiological studies–may differ across individuals. This difference may not be completely random, as behavioral patterns may systematically differ. For example, some seniors may stay indoors more frequently, while younger adults may be more active for various reasons. The risk of a health outcome may vary by age.

The difference between ambient environmental condition and personal exposure to this condition may exhibit differently across sociodemographic factors and time periods (e.g., seasons because behaviors may change). Since the risk of health outcomes often differs by sociodemographic factors, this suggests that EME/CME may be differential by sociodemographic factors. EME/CME may also be differential over a period of time if the behavioral pattern of individuals varies. Such systematic difference can be encoded as $V$. Recall that $V$ is a set of factors that render the error differential. For example, refining the error form mentioned in Section S8.4.

$$X_{personal\ exposure,L_T,t} = \gamma_{ambient}X_{ambient\ level,L_T,t} + V_{S1,L_T,t}$$

where $V_{S1,L_T,t}$ denotes a numeric factor that makes, if individuals were outside 100% during a time interval $(t-1,t]$, $X_{personal\ exposure,L_T,t} = X_{ambient\ level,L_T,t}$. If individuals were inside 100%, $X_{personal\ exposure,L_T,t} = \gamma_{ambient}X_{ambient\ level,L_T,t} +$ $\gamma_{indoor} X_{indoor\ source,L_T,t}$. And if $X_{indoor\ source,L_T,t}$ is independent of $X_{ambient\ level,L_T,t}$, then $\gamma_{indoor} X_{indoor\ source,L_T,t}$ could be seen as $U$. Although this may be over-simplification, this still provides insights into how $V$ could be controlled. See below.

Differential error may also arise in modeled estimates for environmental condition. This is because the error and the risk of health outcomes may exhibit spatiotemporal patterns and one may be correlated with the other. As environmental conditions often exhibit spatiotemporal patterns, which originates from spatiotemporal patterns of their predictors (e.g., emission sources, atmospheric process), the error in exposure prediction models may also follow such patterns, especially if the prediction model may not fully account for spatiotemporal patterns. Additionally, various spatial statistics in mortality and morbidity across the world demonstrate that the risk of health outcomes is not random and may be spatiotemporally patterned. The error may be encoded as

$$X = g(\hat{X}) + V + U$$

where $\hat{X}$ denotes a modelled estimate, $V$ describes the spatiotemporal pattern, and $U$ denotes the remaining completely random error.

The potential of the existence of differential error does not necessarily mean that effect estimates may be biased. For example, $V$ may be (inadvertently/coincidentally) adjusted for in environmental epidemiologic studies when adjusting for confounders. Henceforth, We say “adjustment for differential error” to emphasize adjustment for $V$, while “adjustment for differential error” may not fully adjust for the bias from EME due to the remaining non-differential component such as $U$. Adjustment for differential error may co-incidentally arise when epidemiologists try to adjust for both measured and unmeasured confounders. They may adjust for measured confounders and consider strategies to address unmeasured confounding. To account for unmeasured confounding that may be related to unmeasured spatiotemporal variables, epidemiologists may use robust study designs (e.g., designs based on matching) and employ spatiotemporal statistical techniques (e.g., smoother in generalized additive models). $V$ may be related to/correlated with measured or unmeasured confounders; $V$ may also be a confounder in some circumstances (See Figure S1). In the example above, recall that the error may differ across sociodemographic factors (subsumed to $V$). Sociodemographic factors are often controlled for in air pollution epidemiologic studies by using matching, covariate adjustment, or (generalized) propensity score-based methods (e.g., propensity score matching, stratification, weighting). As air pollution may exhibit spatiotemporal patterns, investigators may utilize statistical methods to adjust for potential unmeasured confounders that may also have spatial, temporal, or spatiotemporal patterns to reduce the risk of unmeasured confounding. If $V$ exhibits such patterns as well, techniques to adjust for unmeasured confounders may also adjust for $V$ (e.g., spatial/temporal smoother, time-stratification matching for each of the individuals like time-stratified case-crossover analysis).

In our application analysis, we used time-stratified case-crossover design. Effect estimation is based on within-individual variation from each of the matched strata (i.e., between a case at an index time and their self-controls at different time points, referred to as reference time points), which means only temporal, not spatial variation in covariates. Therefore, within-individual differential error (due to time-varying factors) may be of concern. Between-individual (e.g., one case vs. another case; self-controls for one case vs. self-controls for another case) differential error does not influence effect estimates. This means that any differential error related to time-invariant individual-level factors (e.g., sociodemographic factors) would not do so.

Within-individual differential error (due to time-varying factors) may be co-incidentally adjusted for, especially if this differential error is related to slowly time-varying factors (e.g., systematic behavior changes over a long period of time) because time-stratified case-crossover design can control for relatively long-term temporal variation by the case–self-control matching. When time-stratification is based on the same calendar year, month, and day-of-the-week, which is a well-established method in air pollution epidemiology and was used by Kim et al. (2022)[18] and our application analysis, what remains concerning is approximately within-two or three-week temporal variation-related differential error because such variation in covariates remains to be used for effect estimation. This means that $V$ that may vary within this two to three-week span should be controlled for.

Kim et al. (2022)[18] conducted rigorous confounding adjustment using the matching and covariate adjustment, which also applies to our application analysis. In addition to case–self-controls matching, covariate adjustment, regarding the two to three-week span, was conducted to adjust for lag-structures of the co-pollutant ($O_3$, for $PM_{2.5}$ estimation; $PM_{2.5}$ for $O_3$, estimation), daily temperature, and relative humidity. Temperature and relative humidity were adjusted for considering their potential non-linear effects on COVID-19 mortality using distributed lag non-linear terms. Temperature and relative humidity are related to infiltration/penetration factor that may change over time[29-31] (See infiltration/penetration factor in Appendix S8.4), implying, if time-varying infiltration/penetration factor (Source 5) had existed, that adjustment for temperature and humidity may have also adjusted for this differential error. Numerically, time-varying infiltration/penetration factor-considered $X_{ambient}$ could be expressed as

$$\gamma_{ambient,avg}X_{ambient} + \gamma^{*}_{ambient}X_{ambient}$$

where $\gamma_{ambient,avg}$ is a long-term temporal average of the time-varying infiltration/penetration factor and $\gamma^{*}_{ambient}$ is the time-varying infiltration/penetration factor minus $\gamma_{ambient,avg}$. $\gamma_{ambient,avg}$ could be seen as constant, which is a non-differential error component and is not addressed but $\gamma^{*}_{ambient}X_{ambient}$ is the differential error component that may have been adjusted for by adjusting for temperature and humidity. Furthermore, lag-structures of daily COVID-19 infectivity was also adjusted for—this may be related to very short-term behavioral change if the change had ever existed, especially in the early stage of the pandemic, meaning that time-varying $V$ from Sources 3–5 related to behavioral changes may have been adjusted for if such $V$ had existed. In our application analysis, as there was no residual autocorrelation in our MRC models, differential errors are unlikely. Thus, the risk of bias due to differential error from EME/CME in effect estimation is likely low.

Recall that the impact of the remaining non-differential part of EME/CME related to Source 5 (i.e., non-differential difference between ambient levels of environmental conditions and personal exposure to environmental conditions from ambient source) on effect estimates remains unclear, but as illustrated in Appendix S8.4, this issue may alternatively be considered from the population average effect perspective.

**Appendix S9. Identifying Potential bias Due to Unmeasured Mobility.**

By plugging the covariates in the main data into the MRC models, MR-Calibrated variables were obtained and linked to cases and their self-controls using the geocoded locations of the time at which COVID-19 death incidents occurred. We overuse $L$ to denote the geocoded location used for the data-linkage and $a$ to denote the spatial grid (i.e., 500m$^2$) that includes $L$. By rewriting MRC models (Eq.1 in the main text) with a subscript $i$,

$$X_{i,L,t,p} = \hat{X}_{i,L,t,p} + U_{i,L,t,p} \text{ (Eq.2)}$$

$X_{i,L,t,p}$ would represent the true level of ambient environmental condition $p$ for an individual $i$ if they remained at $L$ and day $t$.

As $\hat{X}_{i,L,t,1}$, $\hat{X}_{i,L,t,2}$, $\hat{X}_{i,L,t,3}$, and $\hat{X}_{i,L,t,4}$ are MR-Calibrated variables, the use of these variables can reduce bias due to EME and CME from Source 1 (Figure S15), but may not eliminate EME and CME at all, if individuals moved around far from the geocoded location. $\hat{X}_{i,L,t,p}$ may not equal $X_{i,L_T,t,p}$ where $L_T$ denotes the true location of an individual $i$ at a specific time-point due to the potential data-mislinkage due to unmeasured mobility of individuals (Source 3; $L$ or $a$). Even if $L = L_T$ and $a = a_T$ where $a_T$ is the area including $L_T$, $\hat{X}_{i,L,t,p}$ may not equal $X_{i,L,t,p}$, which may be purely Berksonian (Eq.2)—Source 4 error was replaced by this error through MRC. The remaining question concerns the resulting compound error from this error and the potential data-mislinkage.

When potential of the data-mislinkage exists, investigators may wish to assume that the resulting EME/CME from collectively Sources 3 and 4 may introduce minimal or negligible bias (Appendix S8.4). We refer to this assumption regarding EME as Pure Berkson Error from Unmeasured Mobility Assumption (PBEUMA) and regarding CME as Multi-Dimensional PBEUMA. In our application analysis, the compound error from Sources 3 and 4 were replaced by the difference between $\hat{X}_{i,L,t,p}$ and $X_{i,L_T,t,p}$. By replacing $L$ with $L_T$ in Eq.2,

$$X_{i,L_T,t,p} = \hat{X}_{i,L_T,t,p} + U_{i,L_T,t,p} \text{ (Eq.3)}$$

By integrating Eqs.2 and 3,

$$X_{i,L_T,t,p} = \hat{X}_{i,L,t,p} + U_{i,L,t,p} + \left(\hat{X}_{i,L_T,t,p} - X_{i,L,t,p}\right) + U_{i,L_T,t,p} \text{ (Eq.4)}$$

which informs what the resulting compound error model should be for minimal or negligible bias in the application analysis. Note that Eq.4 may or may not hold in real-world applications.

Minimal bias due to EME would be expected if the regression coefficient of the residuals of $X_{i,L_T,t,p}$ against the residuals of $\hat{X}_{i,L,t,p}$ after adjusting for $\boldsymbol{Z}^{**}$ (following the notation $\boldsymbol{z}^{**}$ in Theorem 3) is 1 and the residuals of $\hat{X}_{i,L,t,p}$ is independent of the residuals of $U_{i,L,t,p}$, $U_{i,L_T,t,p}$, and $\left(\hat{X}_{i,L_T,t,p} - X_{i,L,t,p}\right)$ after adjusting for $\boldsymbol{Z}^{**}$ because the resulting conditional error model would become purely Berksonian in outcome modeling that is conditional on $\boldsymbol{Z}^{**}$. This is PBEUMA for $\hat{X}_{i,L,t,p}$.

The risk of potential residual confounding from CME should also be eliminated. A sufficient condition for the elimination is that

1) PBEUMA for $\hat{X}_{i,L,t,p}$ and

2) The independence of the residuals of $\hat{X}_{i,L,t,p}$ and the residuals of $U_{i,L,t,p'}$, $\left(\hat{X}_{i,L_T,t,p'} - X_{i,L,t,p'}\right)$, $U_{i,L_T,t,p'}$ after adjusting for $\boldsymbol{Z}^{**}$ where we overuse $p' \neq p$.

Both needs to hold, which is Multi-Dimensional PBEUMA for $\hat{X}_{i,L,t,p}$.

If individuals had moved around near the geocoded location rather than moving around far from this, (Multi-Dimensional) PBUEMA would have likely held or deviations from this assumption would have been small because the term, $U_{i,L,t,p} + \left(\hat{X}_{i,L_T,t,p} - X_{i,L,t,p}\right) + U_{i,L_T,t,p}$, collectively, would have been small. However, it is possible that some individuals had moved around over long distances.

Therefore, we identified potential bias due to deviations from PBEUMA in hypothetical circumstances where most individuals had moved around over long distances as extreme scenarios that may be believed to introduce substantial bias. As several monitoring stations in the validation data were located far apart (Figure S10), monitoring locations pair permutation data with different weights over multiple repetitions may serve as such hypothetical scenarios to examine the extent of bias due to deviations from (Multi-Dimensional) PBEUMA. See below for details.

To satisfy PBEUMA for $PM_{2.5}$ and $O_3$, the following form, referred to as the PBEUMA model for $\hat{X}_{i,L,t,p}$ (for *p*=1 and *p*=2) should hold, as explained earlier,

$$X_{i,L_j,t,p} = \gamma_{1,p} + \gamma_{2,p}\hat{X}_{i,L_k,t,p} + \gamma_{3,p}U_{i,L_k,t,p} + \gamma_{4,p}\left(\hat{X}_{i,L_j,t,p} - X_{i,L_k,t,p}\right) + U_{i,L_j,t,p} \quad \text{(Eq. S9-1)}$$

where $\gamma_{2,p} = 1$ and the independence of the residuals of $\hat{X}_{i,L_k,t,p}$ and the residuals of $U_{i,L_k,t,p}$, $\left(\hat{X}_{i,L_j,t,p} - X_{i,L_k,t,p}\right)$, and $U_{i,L_j,t,p}$ after adjusting for confounders.

To examine Multi-Dimensional PBEUMA for $PM_{2.5}$ and $O_3$, the following relationships regarding the residuals of terms after adjusting for $\boldsymbol{Z}^{**}$ were examined to determine whether they were uncorrelated:

- For $PM_{2.5}$, the relationship of $\hat{X}_{i,L_k,t,1}|\boldsymbol{Z}^{**}, \boldsymbol{\hat{X}_{i,L_k,t,1}}$ with $U_{i,L_k,t,1}|\boldsymbol{Z}^{**}, \boldsymbol{\hat{X}_{i,L_k,t,1}}$, $\left(\hat{X}_{i,L_j,t,1} - X_{i,L_k,t,1}\right)|\boldsymbol{Z}^{**}, \boldsymbol{\hat{X}_{i,L_k,t,1}}$, $U_{i,L_j,t,1}|\boldsymbol{Z}^{**}, \boldsymbol{\hat{X}_{i,L_k,t,1}}$, $U_{i,L_k,t,2}|\boldsymbol{Z}^{**}, \boldsymbol{\hat{X}_{i,L_k,t,1}}$, $\left(\hat{X}_{i,L_j,t,2} - X_{i,L_k,t,2}\right)|\boldsymbol{Z}^{**}, \boldsymbol{\hat{X}_{i,L_k,t,1}} U_{i,L_j,t,2}|\boldsymbol{Z}^{**}, \boldsymbol{\hat{X}_{i,L_k,t,1}}$, $U_{i,L_k,t,3}|\boldsymbol{Z}^{**}, \boldsymbol{\hat{X}_{i,L_k,t,1}}$, $\left(\hat{X}_{i,L_j,t,3} - X_{i,L_k,t,3}\right)|\boldsymbol{Z}^{**}, \boldsymbol{\hat{X}_{i,L_k,t,1}}$, $U_{i,L_j,t,3}|\boldsymbol{Z}^{**}, \boldsymbol{\hat{X}_{i,L_k,t,1}} U_{i,L_k,t,4}|\boldsymbol{Z}^{**}, \boldsymbol{\hat{X}_{i,L_k,t,1}}$, $\left(\hat{X}_{i,L_j,t,4} - X_{i,L_k,t,4}\right)|\boldsymbol{Z}^{**}, \boldsymbol{\hat{X}_{i,L_k,t,1}} U_{i,L_j,t,4}|\boldsymbol{Z}^{**}, \boldsymbol{\hat{X}_{i,L_k,t,1}}$ (REL. S9-1)
- For $O_3$, the relationship of $\hat{X}_{i,L_k,t,2}|\boldsymbol{Z}^{**}, \boldsymbol{\hat{X}_{i,L_k,t,2}}$ with $U_{i,L_k,t,1}|\boldsymbol{Z}^{**}, \boldsymbol{\hat{X}_{i,L_k,t,2}}$, $\left(\hat{X}_{i,L_j,t,1} - X_{i,L_k,t,1}\right)|\boldsymbol{Z}^{**}, \boldsymbol{\hat{X}_{i,L_k,t,2}}$, $U_{i,L_j,t,1}|\boldsymbol{Z}^{**}, \boldsymbol{\hat{X}_{i,L_k,t,2}}$, $U_{i,L_k,t,2}|\boldsymbol{Z}^{**}, \boldsymbol{\hat{X}_{i,L_k,t,2}}$, $\left(\hat{X}_{i,L_j,t,2} - X_{i,L_k,t,2}\right)|\boldsymbol{Z}^{**}, \boldsymbol{\hat{X}_{i,L_k,t,2}} U_{i,L_j,t,2}|\boldsymbol{Z}^{**}, \boldsymbol{\hat{X}_{i,L_k,t,2}}$, $U_{i,L_k,t,3}|\boldsymbol{Z}^{**}, \boldsymbol{\hat{X}_{i,L_k,t,2}}$, $\left(\hat{X}_{i,L_j,t,3} - X_{i,L_k,t,3}\right)|\boldsymbol{Z}^{**}, \boldsymbol{\hat{X}_{i,L_k,t,2}}$, $U_{i,L_j,t,3}|\boldsymbol{Z}^{**}, \boldsymbol{\hat{X}_{i,L_k,t,2}} U_{i,L_k,t,4}|\boldsymbol{Z}^{**}, \boldsymbol{\hat{X}_{i,L_k,t,2}}$, $\left(\hat{X}_{i,L_j,t,4} - X_{i,L_k,t,4}\right)|\boldsymbol{Z}^{**}, \boldsymbol{\hat{X}_{i,L_k,t,2}} U_{i,L_j,t,4}|\boldsymbol{Z}^{**}, \boldsymbol{\hat{X}_{i,L_k,t,2}}$ (REL. S9-2)

where $\boldsymbol{\hat{X}_{i,L_k,1}} \coloneqq$

$\{\hat{X}_{i,L_k,t-1,1}, .., \hat{X}_{i,L_k,t-l'',1}, \hat{X}_{i,L_k,t,2}, .., \hat{X}_{i,L_k,t-l'',2}, \hat{X}_{i,L_k,t,3}, .., \hat{X}_{i,L_k,t-l'',3}, \hat{X}_{i,L_k,t,4}, .., \hat{X}_{i,L_k,t-l'',4}\}$;

$\boldsymbol{\hat{X}_{i,L_k,2}}$

$\coloneqq \{\hat{X}_{i,L_k,t,1}, .., \hat{X}_{i,L_k,t-l'',1}, \hat{X}_{i,L_k,t-1,2}, .., \hat{X}_{i,L_k,t-l'',2}, \hat{X}_{i,L_k,t,3}, .., \hat{X}_{i,L_k,t-l'',3}, \hat{X}_{i,L_k,t,4}, .., \hat{X}_{i,L_k,t-l'',4}\}$

$L_j$ and $L_k$ indicate a monitoring location and another monitoring location in the validation data, respectively.

An inherent assumption for using the validation data for (Multi-Dimensional) PBEUMA analyses is that most of the decedents had moved within/around Cook County while they were alive and would have continued to do so had they not died, within the exposure time-window of interest (21 days). It is uncertain whether this is true for all individuals. Nevertheless, deviation from this assumption may not be strong because outbound travel may have been limited during the COVID-19 pandemic. Deviation in the study population would not have been strong if only a small subset of individuals had moved outside Cook County, and even if some had moved, environmental conditions may still be correlated with those in Cook County (Appendix S8.2). Thus, using the validation data may not be optimal but reasonable.

By creating 30 permutations (6×5 from the six locations in the validation data), we tested (Multi-Dimensional) PBUEMA. We added random weights to the 30 permutations. The motivation is that weighted analyses may gauge deviations from (Multi-Dimensional) PBEUMA due to unmeasured mobility patterns of individuals. For illustration, suppose that $L_1$ denotes the location at which most individuals in the study population stayed at time $t$; $L_2$ denotes the location used to link exposure data and the health data so that the mismatch arises. Using $L_1$-$L_2$ paired data, whether the resulting error from this mismatch may not violate PBEUMA can be examined. Since which permutation(s) was more likely to represent most potential mismatches that may have occurred in the data-linkage in our application analysis is unknown, random numbers were drawn from the uniform distribution [0,1] as weights and repeated weighted analyses 100 times. While acknowledging that 30 permutations may not represent all possible mismatches–given that Cook County is geographically large and individuals may have visited outside of it–analyzing this at least increases the confidence in determining the likelihood that (Multi-Dimensional) PBEUMA is valid or to gauge potential bias due to deviations from (Multi-Dimensional) PBEUMA, given that environmental conditions at one location were positively correlated with those at another location (Appendix S7.2).

To increase the likelihood that PBEUMA for $PM_{2.5}$ and $O_3$ holds, the central tendency of the distribution of the estimates of $\gamma_{2,1}$ and $\gamma_{2,2}$ (See Eq. S9-1 for notations) is (approximately) should be nearly 1 and there should not exist extreme values deviated from the central tendency across 100 weighted analyses. Also, the relationships mentioned above regarding Multi-Dimensional PBEUMA (REL. S9-1, REL. S9-2) should have the linear coefficient of zero (i.e., no correlation) without its extreme values deviated from the central tendency, even if the absence of correlation may not equal independence, at least theoretically.

If there exist deviations from (Multi-Dimensional) PBEUMA, potential bias may be calculated, in particular if the effect of environmental conditions is expressed as linear terms in an outcome model. For illustration, suppose a simple true outcome model,

$$Y_{i,L_T,t} = \alpha + \beta_1 X_{i,L_T,t,1} + \beta_2 X_{i,L_T,t,2} + \beta_3 C_{i,L_T,t}$$

This true outcome model can be re-expressed as, by plugging Eq. S9-1 into this,

$$Y_{i,t} = \alpha + \beta_1[\gamma_{1,1} + \gamma_{2,1}\hat{X}_{i,L,t,1} + \gamma_{3,1}U_{i,L,t,1} + \gamma_{4,1}(\hat{X}_{i,L_T,t,1} - X_{i,L,t,1}) + U_{i,L,t,1}] + \beta_2[\gamma_{1,2} + \gamma_{2,p}\hat{X}_{i,L,t,2} + \gamma_{3,p}U_{i,L,t,2} + \gamma_{4,2}(\hat{X}_{i,L_T,t,2} - X_{i,L,t,2}) + U_{i,L,t,2}] + \beta_3 C_{i,t}$$

This can be re-written as

$$Y_{i,t} = \alpha + \beta_1\gamma_{1,1} + \beta_1\gamma_{2,1}\hat{X}_{i,L,t,1} + \beta_1\gamma_{3,1}U_{i,L,t,1} + \beta_1\gamma_{4,1}(\hat{X}_{i,L_T,t,1} - X_{i,L,t,1}) + \beta_1 U_{i,L,t,1} + \beta_2\gamma_{1,2} + \beta_2\gamma_{2,2}\hat{X}_{i,L,t,2} + \beta_2\gamma_{3,2}U_{i,L,t,2} + \beta_2\gamma_{4,2}(\hat{X}_{i,L_T,t,2} - X_{i,L,t,2}) + \beta_2 U_{i,L,t,2} + \beta_3 C_{i,t}$$

Consider that an investigator uses an incorrect model to estimate $\beta_1$ using $\hat{X}_{i,L,t,1}$ and $\hat{X}_{i,L,t,2}$:

$$Y_{i,L_T,t} = \alpha^{ep} + \beta_1^{ep}\hat{X}_{i,L,t,1} + \beta_2^{ep}\hat{X}_{i,L,t,2} + \beta_3^{ep}C_{i,L_T,t}$$

The biased regression coefficient for $X_1$, $\beta_1^{ep}$ is

$$\beta_1^{ep} \cong \beta_1\gamma_{2,1} + \vartheta_{U_{i,L,t,1},\hat{X}_{i,L,t,1}|C}\beta_1\gamma_{3,1} + \vartheta_{(\hat{X}_{i,L_T,t,1} - X_{i,L,t,1}),\hat{X}_{i,L,t,1}|C}\beta_1\gamma_{4,1} + \vartheta_{U_{i,L,t,1},\hat{X}_{i,L,t,1}|C}\beta_1 + \vartheta_{U_{i,L,t,2},\hat{X}_{i,L,t,1}|C}\beta_2\gamma_{3,2} + \vartheta_{(\hat{X}_{i,L_T,t,2} - X_{i,L,t,2}),\hat{X}_{i,L,t,1}|C}\beta_2\gamma_{4,2} + \vartheta_{U_{i,L,t,2},\hat{X}_{i,L,t,1}|C}\beta_2$$

where $\vartheta_{term\ 1,term\ 2|C}$ denotes the regression coefficient of the residuals of $term\ 1$ against the residuals of $term\ 2$ after adjusting for $C$.

Potential bias due to deviations from PBEUMA for $X_1$ is related to:

$$\beta_1\gamma_{2,1} + \vartheta_{U_{i,L,t,1},\hat{X}_{i,L,t,1}|C}\beta_1\gamma_{3,1} + \vartheta_{(\hat{X}_{i,L_T,t,1} - X_{i,L,t,1}),\hat{X}_{i,L,t,1}|C}\beta_1\gamma_{4,1} + \vartheta_{U_{i,L,t,1},\hat{X}_{i,L,t,1}|C}\beta_1$$

Potential bias due to deviations from Multi-Dimensional PBEUMA for $X_2$ (i.e., residual confounding by CME) is related to:

$$\vartheta_{U_{i,L,t,2},\hat{X}_{i,L,t,1}|C}\beta_2\gamma_{3,2} + \vartheta_{(\hat{X}_{i,L_T,t,2} - X_{i,L,t,2}),\hat{X}_{i,L,t,1}|C}\beta_2\gamma_{4,2} + \vartheta_{U_{i,L,t,2},\hat{X}_{i,L,t,1}|C}\beta_2$$

Similarly, these equations can be extended to more than two $\hat{X}$s. For non-linear models, these are approximations.

Figure S16 presents the distribution of the estimates for $\gamma_{2,1}$ and $\gamma_{2,2}$ in the PBEUMA models (Eq. S9.1) conditional on the confounders that were adjusted for in outcome modeling from 100 differently weighted permutation datasets. $\gamma_{2,1}$ and $\gamma_{2,2}$ should be nearly 1 to increase the likelihood of PBEUMA. However, regarding $PM_{2.5}$ effect estimation, $\gamma_{2,1}$ values were around 0.5, meaning that underestimation (50% of the true effect) due to this deviation from PBEUMA may have occurred if most individuals had moved around over long distances. Similarly, regarding $O_3$ effect estimation, $\gamma_{2,2}$ values also around 0.5, meaning that underestimation (50% of the true effect) may have occurred.

Figure S17A presents the distribution of the estimates of the regression coefficient of the other terms in the PBEUMA models for $PM_{2.5}$, $O_3$, temperature, and relative humidity than $\hat{X}_{i,L_k,t,1}$ (Eq. S9-1, $p = 1,2,3, or\ 4$) against $\hat{X}_{i,L_k,t,1}$ to evaluate REL S9-1. This coefficient should be nearly 0 for Multi-Dimensional PBUEMA for $PM_{2.5}$. The coefficient estimates for most of the terms were within ±0.1, meaning that deviations from Multi-Dimensional PBEUMA were small. One term, $\hat{X}_{i,L_j,t,2} - X_{i,L_k,t,2}$ had an estimate up to approximately 0.17. To interpret this, the coefficient of this term, which was $\gamma_{4,2}$ in the PBEUMA model for $O_3$ (Figure S16B) is needed. The product of 0.17 and $\gamma_{4,2} \approx 0.7$ represents an approximate of the upper bound for the potential bias due to residual confounding by $O_3$, which is $0.17 \times 0.7 \times 100\% \approx 11.9\%$. This implies that residual confounding by $O_3$ through $\hat{X}_{i,L_j,t,2} - X_{i,L_k,t,2}$ may account for ~11.9% of the size of the effect of $O_3$ on COVID-19 mortality in the $PM_{2.5}$ effect estimation. As the effect of $O_3$ was found to be positive (See the main manuscript), away-from-the-null bias may have occurred if most individuals had moved around over long distances. Similarly, $U_{i,L_k,t,2}$ had an estimate up to approximately 0.07. The coefficient of this term was approximately up to $-0.6$ ($\gamma_{3,2}$) in the PBEUMA model for $O_3$ (Figure 16B). Thus, residual confounding by $O_3$ through $U_{i,L_k,t,2}$ may account for, as an approximate of the upper bound, $0.07 \times -0.6 \times 100\% \approx -3.5\%$ of the size of the effect of $O_3$ on COVID-19 if most individuals had moved around over long distances. Also note that $\hat{X}_{i,L_j,t,1} - X_{i,L_k,t,1}$ had an estimate up to approximately -0.1 (Figure S17A). The coefficient of this term, which is $\gamma_{4,1}$ in the PBEUMA model for $PM_{2.5}$ (Figure S16A) had an estimate up to approximately 0.8. Thus, $PM_{2.5}$ effect estimates may have been biased by, as an approximate of the upper bound, $-0.1 \times 0.8 \times 100\% \approx -8\%$ of the size of the effect of $PM_{2.5}$ on COVID-19 mortality if most individuals had moved around over long distances, which leads to additional deviation from PBEUMA for $PM_{2.5}$. Collectively, as possibly ~50% underestimation through $\gamma_{2,1} \approx 0.5$ (in the previous paragraph), in summary, for $PM_{2.5}$ effect estimation, deviation from PBEUMA and Multi-Dimensional PBEUMA for $PM_{2.5}$ may have introduced toward-the-null net bias to some degree, if the majority of individuals had moved around over long distances.

For $O_3$, Figure S17B demonstrates deviation from PBEUMA may also have come from $\hat{X}_{i,L_j,t,2} - X_{i,L_k,t,2}$ that had an estimate up to ~-0.1 with $\gamma_{4,2} \approx 0.7$, bias by up to $-0.1 \times 0.7 \times 100\% \approx -7\%$ of the size of the effect of $O_3$ on COVID-19 mortality may have occurred if the majority of individuals had moved around over long distances. Multi-Dimensional PBUEMA may not have been deviated much. Collectively, as possibly ~50% underestimation through $\gamma_{2,2} \approx 0.5$ (in the second preceding paragraph), in summary, for $O_3$ effect estimation, deviation from PBEUMA for $O_3$ may have introduced toward-the-null net bias to some degree, if the majority of individuals had moved around over long distances.

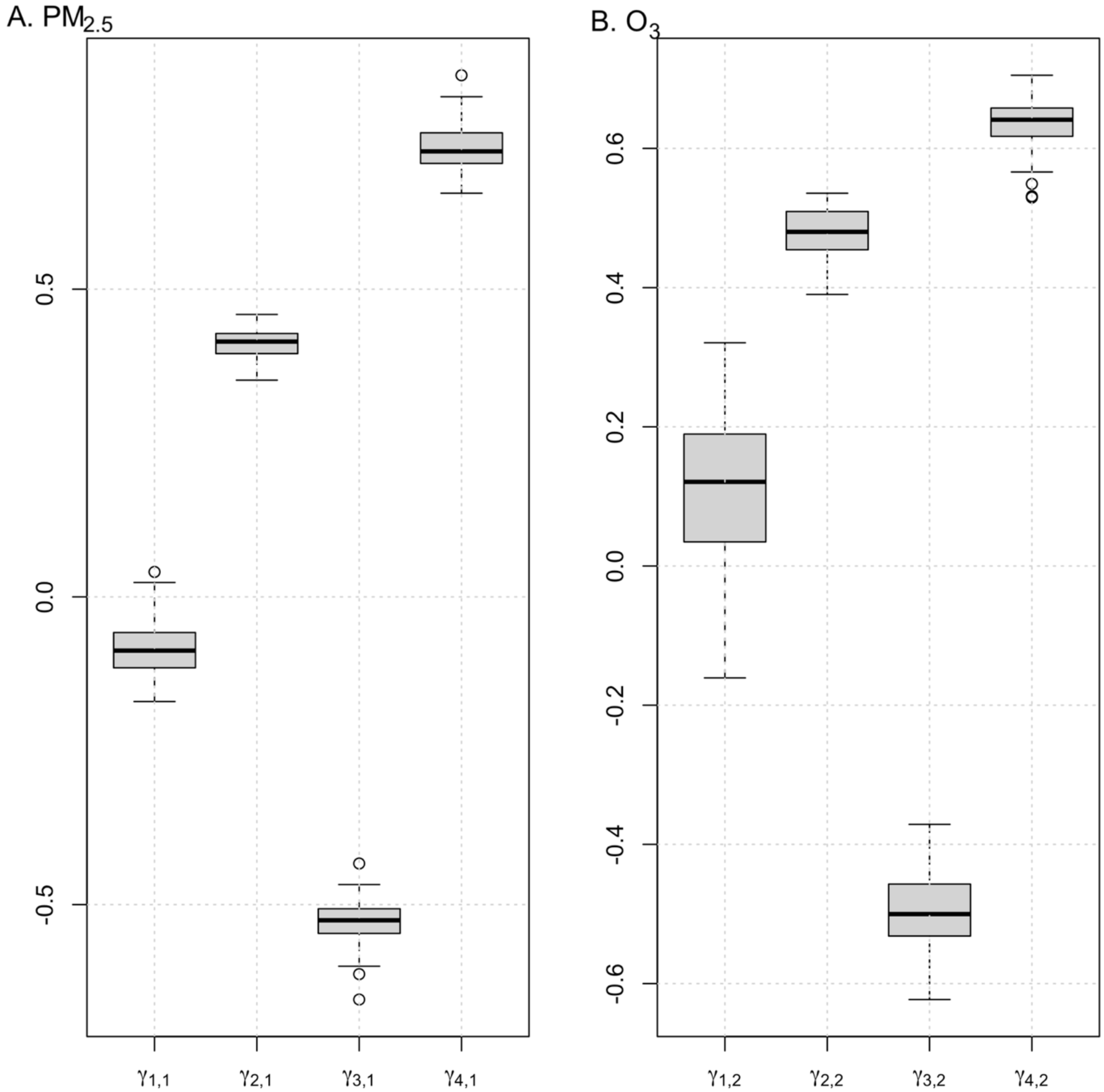

**Figure S16. The distribution of the coefficient estimates of PBEUMA models for $PM_{2.5}$ (A) and $O_3$ (B) conditional on the variables/terms used for confounding adjustment.**

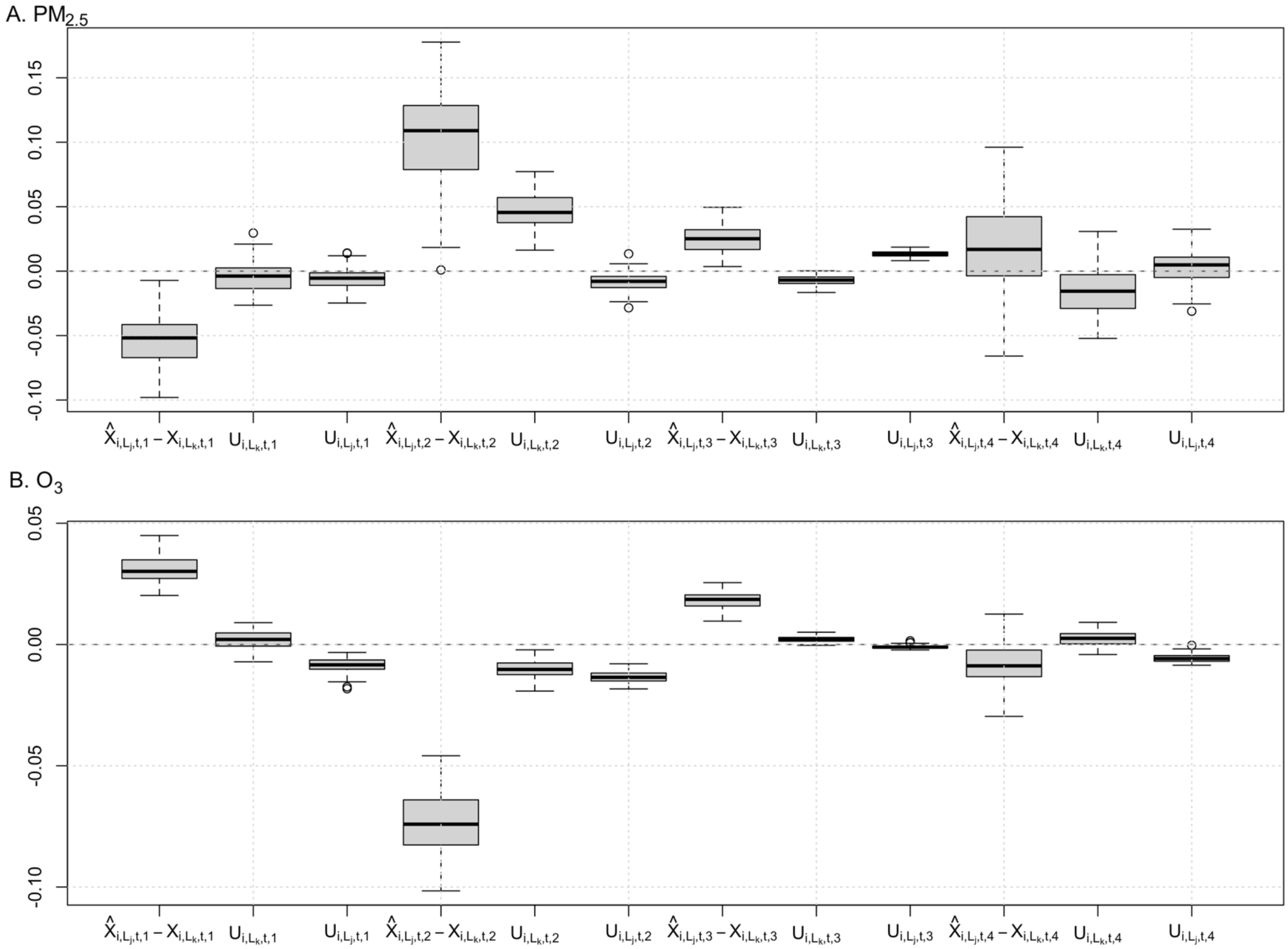

**Figure S17. The regression coefficient estimates of the other terms in the PBEUMA models against MRC-calibrated variables conditional on the variables/terms used for confounding adjustment: $PM_{2.5}$ (A) and $O_3$ (B).**

**Appendix S10. Estimation of Conditional AEE in the RR and RD scales and of Preventable Premature Deaths.**

**S10.1. Estimation of Conditional AEE in the RR and RD scales**

To estimate conditional AEE in the RD and RD scales (Section S6.2 in Online Appendix S6), the probability of *Y*=1 conditional on covariates $P(Y = 1|covariates)$ needs to be estimated using estimated Q-model.

Generally, conditional logistic regression is used to estimate the effects of covariates when the study design is a matched case-control design, including a time-stratified case-crossover design. The intercept of a logistic regression is not estimable in conditional logistic regression due to the conditional likelihood. Furthermore, our analysis is based on case and self-controls data, not a cohort data. Collectively, $P(Y = 1|covariates)$, which is the probability of *Y*=1 for each individual in a population, cannot be directly estimated. External data is needed to estimate this.

A logistic model as Q-model may take the following form,

$$\log\left(\frac{P_i(Y=1)}{P_i(Y=0)}\right) = \beta_0 + \boldsymbol{\beta}\boldsymbol{\mathcal{X}_i}$$

which could have been estimated from cohort data had this existed. By using this model,

$$P_i(Y = 1|\boldsymbol{\mathcal{X}}) = \frac{1}{1 + \exp\left(-(\beta_0 + \boldsymbol{\beta}\boldsymbol{\mathcal{X}_i})\right)}$$

could have been estimated.

However, our application analysis was based on case(–self-controls) data from Cook County Medical Examiner's Office. $\boldsymbol{\beta}$ was estimated using a time-stratified case-crossover analysis with conditional logistic regression. We estimated $\beta_0$ using some assumptions and external data about the empirical COVID-19 case-fatality ratio, which is the total number of observed deaths from COVID-19 divided by the total number of people infected with COVID-19.

For illustration, the conditional likelihood of a case–self-control matched stratum of size $m$, $s$, with the first observation being the case and the others being self-controls, is

$$P\left(Y_{s,1} = 1 \,\&\, Y_{s,j} = 0 \; for \; 1 < j \leq m | X_1, \sum_{\boldsymbol{j=1}}^{\boldsymbol{m}} Y_{s,j} = 1\right) = \frac{\exp(\boldsymbol{\beta}^* \boldsymbol{\mathcal{X}_1})}{\sum_{j=1}^{m} \exp(\boldsymbol{\beta}^* \boldsymbol{\mathcal{X}}_j)}$$

This is the conditional probability of *Y*=1 for each matched stratum that has only one case and 3 or 4 self-controls, (which was concentrated between 0.2 and 0.8, suggesting good approximation conditions to estimate AEE using our developed framework). And the full conditional log-likelihood for our conditional logistic regression is

$$\log\left(\prod_{s=1}^{n} P(Y_{s,1} = 1 \,\& Y_{s,j} = 0 \; for \; 1 < j \leq m | X_1, \sum_{\boldsymbol{j=1}}^{\boldsymbol{m}} Y_{s,j} = 1)\right)$$

It is reasonable to assume that $\boldsymbol{\beta}^* \approx \boldsymbol{\beta}$. This is because Cook County Examiner's Office COVID-19 Case Archive data included nearly all COVID-19 deaths in Cook County from March to February 2021 so that the source population can be seen as the Cook County as a whole. Therefore, the target estimand of our application analysis would be identical to or resemble what would have been targeted in a matched (nested) case-control study within a Cook County-wide full (or representative) cohort data had this data existed. For readers unfamiliar with this, please see the relationship among population-wide time-series design (often, with population mortality registry in environmental epidemiology), open cohort design, and case-crossover design[32,33] and additionally recent case-time-series design[34] (although effect estimates may be subject to varying perturbations of statistical modeling assumptions over analyses) and the relationship of effect measures between different epidemiologic study designs regarding effect measures[35].

As people cannot die from COVID-19 unless they are infected, the probability of interest may become

$$P_i(Y = 1|\boldsymbol{\mathcal{X}}, I = 1) = \frac{1}{1 + \exp\left(-(\beta_0 + \boldsymbol{\beta}\boldsymbol{\mathcal{X}} + \beta_I(I = 1))\right)}$$

where $I = 1$ encodes COVID-infection. $\beta_I(I = 1)$ may be subsumed into $\beta_0$, if the analysis was restricted to infected individuals, so that

$$P_i(Y = 1|\boldsymbol{\mathcal{X}}, I = 1) = \frac{1}{1 + \exp\left(-(\beta_0 + \boldsymbol{\beta}\boldsymbol{\mathcal{X}})\right)}$$

To use regression coefficient estimates from the estimated Q-model for $\widehat{\boldsymbol{\beta}}$, the meaning of these estimates should be considered. Regarding the effect of air pollution on COVID-19 mortality risk (i.e., $\beta$ for air pollutant variables), at least two potential pathways may exist: (1) air pollution increases the risk of death from COVID-19 among individuals already infected; and (2) air pollution increases the risk of COVID-19 infection itself (if this exists). From a mediation analysis perspective, the first pathway may be considered the direct effect of air pollution on COVID-19 mortality, while the second pathway may represent an indirect effect (i.e., air pollution exposures increase the probability of $I = 1$). Since the data used in our application analysis did not include information about the timing of COVID-19 infection in self-controls–which typically precedes the diagnosis–the effect estimates from our conditional regression analysis may fully reflect the direct effect and partially reflect indirect effects, particularly if individuals become infected with COVID-19 within the 21-day period preceding the date of death (e.g., being sentenced to death) but we cannot isolate these two. Considering the incubation period of COVID-19[25], which typically ranges from a few days to a few weeks with a mean of 1-week, it appears that the indirect effect may contribute minimally to our effect estimates. This minimal contribution assumption implies regression coefficient estimates from the estimated Q-model can serve as $\hat{\beta}$. The consequence of lifting this assumption will be discussed below.

The empirical case-fatality ratio can be used to estimate $\beta_0$. For illustration, consider

$$P_i(Y = 1|\boldsymbol{X} = \boldsymbol{0}, I = 1) = \frac{1}{1 + \exp\left(-(\beta_0)\right)}$$

implies that the probability of dying from COVID-19 once an individual who does not have any risk factors of COVID-19 mortality is infected. $P(Y = 1|\boldsymbol{X} = \boldsymbol{0}, I = 1) = E[P_i(Y = 1|\boldsymbol{X} = \boldsymbol{0}, I = 1)]$ would be what the COVID-19 case-fatality ratio would have been had individuals not been exposed to any risk factors. In the real-world, this counterfactual case-fatality ratio is not observable, because people were exposed to various risk factors including environmental hazards. Instead, we can write

$P(Y = 1|I = 1)$,

$$P(Y=1|I=1) = E[P_i(Y=1|I=1)] = E\left[E_{\mathbf{Z}}\left[E_{X_2}\left[E_{X_1}[P_i(Y=1|I=1,X_1,X_2,\mathbf{Z})]\right]\right]\right]$$

$$= E_{\mathbf{Z}}\left[E_{X_2}\left[E_{X_1}\left[\frac{1}{1+\exp\left(-(\beta_0 + f_{X_1}(X_1) + f_{X_2}(X_2) + f_{\mathbf{z}}(\mathbf{Z})\right)}\right]\right]\right]$$

$$\approx Empirical\ (Observed)\ Case - Fatality\ Ratio$$

where $X_1$ is $PM_{2.5}$, $X_2$ is $O_3$, and $\mathbf{Z}$ includes confounders.

Using conditional logistic regression, $f_{X_1}$, $f_{X_2}$, and $f_{\mathbf{z}}(\mathbf{Z})$ can be estimated. To estimate $\beta_0$, we plugged into a certain value to $\beta_0$ and calculate $P(Y=1|I=1)$ and compare this with the empirical case-fatality ratio, iteratively. To calculate $P(Y=1|I=1)$, $P_i(Y=1|I=1,X_1,X_2,\mathbf{Z})$ for all residents in Cook County over the study period should be estimated and averaged. For this, a Cook County-wide full cohort data including all individuals infected with COVID-19 with their values of $X_1$, $X_2$, and $\mathbf{Z}$ could have been used had it existed. Although we do not have information for all residents, it is reasonable to assume that individuals who stayed in an identical location at a specific time shared exposures to ambient $PM_{2.5}$, $O_3$, temperature, and relative humidity. Thus, $P_i(Y=1|I=1,X_1,X_2,\mathbf{Z})$ may be estimated by applying $\hat{f}_{X_1}$, $\hat{f}_{X_2}$, and $\hat{f}_{\mathbf{z}}$ to $X_1$, $X_2$, and $\mathbf{Z}$ in all grid cells over the study period (i.e., the daily time-series of $X_1$, $X_2$, and $\mathbf{Z}$ for all grid cells) and calculating averages weighted by the infected population number of each grid cell. However, we do not have data for the number of infected individuals of each grid cell. We instead conducted this calculation for only the geocoded locations for the decedents due to COVID-19. This approach may not be very accurate because the COVID-19 mortality risk and the infection risk must not have been geographically homogeneous. Nevertheless, obtaining a highly accurate estimate of the intercept may not be necessary, as the marginal COVID-19 mortality risk for Cook County is very low (i.e., with the empirical case-fatality ratio ranging from 2% to 5%, according to the USAFatcts COVID-19 statistics), meaning that AEE estimates in the RR and RD scales would not be sensitive to the point estimate of $\beta_0$ while there are also uncertainties in $\hat{f}_{X_1}$, $\hat{f}_{X_2}$, and $\hat{f}_{\boldsymbol{z}}$. As a result, $\hat{\beta}_0 \approx \log(0.0024)$ was selected such that $E\left[\hat{P}_i(Y=1|I=1)\right] \approx 3.5\%$.

If the minimal contribution assumption does not hold, then, we consider

$$P_i(Y=1|\boldsymbol{\mathcal{X}}) = \frac{1}{1+\exp\left(-(\beta_0 + \boldsymbol{\beta}\boldsymbol{\mathcal{X}})\right)}$$

$\beta_0$ should correspond to the probability of dying from COVID-19 among both those who were infected and those who were not yet infected (i.e., getting infected and then dying from COVID-19). The probability is smaller than that for only those who were infected. This means that $\beta_0$ would be smaller than $\log(0.0024)$. Because this is already a very small quantity, the violation of the minimal contribution assumption would impact very negligibly AEE estimates in the RD and RR scales. But the interpretation of AEE estimates should differ in nuance as $\beta$ for air pollutants here represents the total effect of air pollution on COVID-19 mortality, including both the direct and the indirect effects as noted above. Because the data used did not have information about when individuals were infected, the total effect of air pollution may have been underestimated in our application analysis if the indirect effect had existed (which is unlikely protective).

For AEE estimation, had a Cook County-wide full (or representative) cohort existed, a Q-model could have been applied to everyone in the cohort, and population risks could have been obtained by averaging individuals' estimated probability, which would be marginal AEE. In this application analysis, this is not necessary because our estimand is conditional AEE, not marginal AEE.

### S10.2. Estimation of Preventable Premature Deaths

We estimated preventable premature death estimates as

$$\sum_{i \in D}\left[\frac{P_i\left(Y=1 \middle| \boldsymbol{X_1}=\boldsymbol{X}_{\boldsymbol{RC},\boldsymbol{1}}^{\boldsymbol{ep}}, \boldsymbol{X_2}=\boldsymbol{X}_{\boldsymbol{RC},\boldsymbol{2}}^{\boldsymbol{ep}}\right)-P_i\left(Y=1 \middle| \boldsymbol{X_1}=\boldsymbol{X}_{\boldsymbol{1},\boldsymbol{ref}}, \boldsymbol{X_2}=\boldsymbol{X}_{\boldsymbol{2},\boldsymbol{ref}}\right)}{P_i\left(Y=1 \middle| \boldsymbol{X_1}=\boldsymbol{X}_{\boldsymbol{1},\boldsymbol{ref}}, \boldsymbol{X_2}=\boldsymbol{X}_{\boldsymbol{2},\boldsymbol{ref}}\right)}\right]$$

where $D$ is a set of 7,507 deaths from COVID-19 and $\boldsymbol{X_1}=\boldsymbol{X}_{\boldsymbol{RC},\boldsymbol{1}}^{\boldsymbol{ep}}, \boldsymbol{X_2}=\boldsymbol{X}_{\boldsymbol{RC},\boldsymbol{2}}^{\boldsymbol{ep}}$ indicates natural exposure values (3-week cumulative period). We constructed several sets of the reference values below. $X_{p,ref}$ at a given lag (within 3-week cumulative period) is set to be the following value if the natural exposure value at that lag exceeds this:

1) $X_{1,ref}=5$ μg/m$^3$ and $X_{2,ref}=51$ ppb (Reduction in both pollutants)
2) $X_{1,ref}=5$ μg/m$^3$ and $X_{2,ref}$ remains natural ($PM_{2.5}$ reduction only)
3) $X_{1,ref}$ remains natural and $X_{2,ref}=51$ ppb ($O_3$ reduction only)
4) $X_{1,ref}=2$ μg/m$^3$ and $X_{2,ref}$: 7ppb (Jan.); 18ppb (Feb.); 8ppb (Mar.); 16ppb (Apr.); 25ppb (May); 27ppb (Jun.); 26ppb (Jul.); 17ppb (Aug.); 12ppb (Sep.); 5ppb (Oct.); 10ppb (Nov.); 3ppb (Dec.)
5) $X_{1,ref}=2$ μg/m$^3$ and $X_{2,ref}$ remains natural ($PM_{2.5}$ reduction only)
6) $X_{1,ref}$ remains natural and $X_{2,ref}$: 7ppb (Jan.); 18ppb (Feb.); 8ppb (Mar.); 16ppb (Apr.); 25ppb (May); 27ppb (Jun.); 26ppb (Jul.); 17ppb (Aug.); 12ppb (Sep.); 5ppb (Oct.); 10ppb (Nov.); 3ppb (Dec.) ($O_3$ reduction only)

**Table S3. AEE estimates in the RD and RR scales for cumulative 3-week exposure to $PM_{2.5}$ (per 5.6μg/m$^3$) and $O_3$ (per 22.7ppb).**

| Type of Analysis | Pollutant | Effect Measure | Estimates (95% CI) All | Home/Work | Home/Work & Diff. ≤7d |
|---|---|---|---|---|---|
| Date of Incident | $PM_{2.5}$ | RD* | 0.0027 (0.0012, 0.0039) | 0.0024 (0.0011, 0.0039) | 0.0032 (0.0018, 0.0047) |
| | | psRD* | 0.26 (0.15, 0.32) | 0.24 (0.14, 0.31) | 0.29 (0.20, 0.34) |
| | | RR | 2.11 (1.49, 2.62) | 1.99 (1.45, 2.61) | 2.34 (1.76, 2.97) |
| | $O_3$ | RD* | 0.0009 (0.0001, 0.0023) | 0.0021 (0.0008, 0.0042) | 0.0033 (0.0017, 0.0056) |
| | | psRD* | 0.12 (0.01, 0.24) | 0.23 (0.11, 0.33) | 0.29 (0.19, 0.36) |
| | | RR | 1.37 (1.03, 1.98) | 1.89 (1.34, 2.77) | 2.36 (1.69, 3.35) |
| Date of Death | $PM_{2.5}$ | RD* | 0.0016 (0.0005, 0.0025) | 0.0011 (0.0003, 0.0020) | 0.0011 (0.0004, 0.0022) |
| | | psRD* | 0.18 (0.06, 0.25) | 0.14 (0.05, 0.22) | 0.14 (0.05, 0.23) |
| | | RR | 1.67 (1.19, 2.05) | 1.47 (1.13, 1.84) | 1.46 (1.15, 1.90) |
| | $O_3$ | RD* | 0.0001 (-0.0004, 0.0011) | 0.0009 (0.0001, 0.0023) | 0.0023 (0.0008, 0.0049) |
| | | psRD* | 0.02 (-0.06, 0.13) | 0.11 (0.02, 0.24) | 0.24 (0.10, 0.35) |
| | | RR | 1.05 (0.85, 1.45) | 1.37 (1.05, 1.96) | 1.97 (1.33, 3.03) |

*RD describes the increase in the probability of dying from COVID-19 per one IQR increase in a 3-week exposure at a specific day so that the quantity may be small. As individuals are exposed to air pollution every day, we defined and estimated period-specific RD (psRD) as the increase in the probability of dying from COVID-19 per one IQR increase in a 3-week exposure over March. 2020 to Feburary 2021 (365 days). Since the outcome is death, meaning this can occur only once for each individual, mathematically, psRD was estimated as $psRD \coloneqq 1-(1-p_2)^{365}-(1-(1-p_1)^{365})$ where $p_1$ is the baseline probability (i.e., 0.0024) and $p_2$ is the probability of dying from COVID-19 at one IQR value. The bound of this measure is [0,1]. For example, psRD=0.29 means that the probability of dying from COVID-19 between March 2020 and February 2021 would increase by 0.29 if an individual were exposed to one additional IQR level of an air pollutant every day, compared to a given baseline exposure level.

**References for Online Supplementary Materials (Appendices)**